\documentclass[12pt]{article}
\usepackage{amsmath, calc, amssymb, mathabx}
\usepackage{bm, amsthm}
\usepackage{upgreek}
\usepackage{float}
\usepackage{natbib}
\RequirePackage[colorlinks,citecolor=blue,urlcolor=blue]{hyperref}
\RequirePackage{graphicx}
\usepackage{booktabs}
\usepackage{multirow}
\usepackage{dsfont}
\usepackage[bbgreekl]{mathbbol}
\DeclareSymbolFontAlphabet{\mathbbl}{bbold}
\usepackage{authblk}

\usepackage[disable]{todonotes}
\newcommand{\blind}{1}

\newtheorem{theorem}{Theorem}[section]
\newtheorem{lemma}[theorem]{Lemma}

\newtheorem{corollary}[theorem]{Corollary}
\newtheorem{assumption}{Assumption}

\theoremstyle{remark}\newtheorem{remark}[theorem]{Remark}\theoremstyle{plain}
\newcommand{\E}{\mathrm{E}}
\newcommand{\EIF}{S^{\text{eff, nonpar}}}

\newcommand{\EstFunTraj}{{\psi}^{\text{traj}}}

\newcommand{\EstFunPoi}{{\psi}^{\text{point}}}

\newcommand{\EstestFunTraj}{\widehat{\psi}^{\text{traj-step}}}

\newcommand{\pr}{\mathrm{P}}
\newcommand{\pn}{\mathbb{P}_N}

\newcommand{\pnt}{\mathbb{P}_{NT}}

\newcommand{\var}{\operatorname{var}}

\newcommand{\bs}{\boldsymbol}

\newcommand{\NA}{---}

\newcommand{\bcdot}{\boldsymbol{\cdot}}

\newcommand{\UB}{\operatorname{UB}}
\newcommand{\sm}{\mathbf{sm}}

\newcommand{\vv}{\mathbf{v}}

\newcommand{\va}{\mathbf{a}}

\newlength\oversetwidth
\newlength\underwidth

\begin{document}

\def\spacingset#1{\renewcommand{\baselinestretch}%
{#1}\small\normalsize} \spacingset{1}


\def\spacingset#1{\renewcommand{\baselinestretch}%
{#1}\small\normalsize} \spacingset{1}


\if1\blind
{
  \title{\bf Semiparametric Off-Policy Inference for Optimal Policy Values under Possible Non-Uniqueness}
  \author{
    Haoyu Wei\thanks{Email: h8wei@ucsd.edu \\ We are grateful to Kengo Kato, Runze Li, and the anonymous referees for insightful comments and valuable suggestions.} \hspace{.2cm}\\
    {\small Department of Economics, University of California, San Diego, USA}}
  \maketitle
} \fi

\if0\blind
{
  \bigskip
  \bigskip
  \bigskip
    \title{\bf Semiparametric Off-Policy Inference for Optimal Policy Values under Possible Non-Uniqueness}
  \medskip
  \maketitle
} \fi

\bigskip
\begin{abstract}
Off-policy evaluation (OPE) constructs confidence intervals for the value of a target policy using data generated under a different behavior policy.
Most existing inference methods focus on fixed target policies and may fail when the target policy is estimated as optimal, particularly when the optimal policy is non-unique or nearly deterministic.

We study inference for the value of optimal policies in Markov decision processes.
In an auxiliary augmented transition-sampling experiment, we
characterize the existence of the efficient influence function and
show that non-regularity arises when competing optimal policies have
distinct first-order gradients.  For the actual i.i.d.-trajectory
experiment, we derive the semiparametric efficiency bound and a
uniformly weighted estimator that attains it under a unique optimum,
while the sequential NSAVE procedure trades efficiency for stability
and validity under non-uniqueness.%

Motivated by this analysis, we propose a novel \textit{N}onparametric \textit{S}equenti\textit{A}l \textit{V}alue \textit{E}valuation (NSAVE) method, which yields martingale-based inference and retains a double-robustness property under policy-aligned nuisance estimation.
We further develop a pointwise smoothing-based approximation under
explicit first-stage rates, and a post-selection template with
uniform coverage whenever its stated joint calibration condition is
satisfied.

Simulation studies support the theoretical results.
An application to the Drink Less micro-randomized trial provides confidence intervals for state-adaptive notification policies and their improvement over the randomized behavior policy.
\end{abstract}
\par
\noindent
{\it Keywords:} Doubly robust estimation; non-regular inference; optimal policy; off-policy evaluation.
\vfill

\newpage
\spacingset{1.9}

\maketitle

\section{Introduction}\label{sec:Intro}

Reinforcement learning (RL) is concerned with learning optimal decision rules for sequential decision problems in order to maximize long-term cumulative rewards \citep{sutton2018reinforcement}. 
A fundamental statistical task within RL is off-policy evaluation (OPE), which seeks to estimate the value of a target policy using data generated under a potentially distinct behavior policy. 
OPE plays a pivotal role in offline RL, where new data collection is either costly or ethically constrained, necessitating that inference rely entirely on historical trajectories \citep{luedtke2016statistical,agarwal2019reinforcement,uehara2022review}.

Mobile health interventions provide a concrete example in which this inferential problem is not merely technical.
In the Drink Less micro-randomized trial, users of an alcohol-reduction app were randomized daily to receive no notification, a new notification, or the standard notification, and the primary proximal outcome recorded whether the user opened the app in the hour after the decision point \citep{bell2023notifications}.
The original trial analysis established near-term causal effects of notifications and explicitly motivated further optimization of the notification policy.
However, that analysis focused on marginal proximal effects rather than inference for the value of a learned, state-adaptive policy over a discounted sequential horizon.
This is precisely the setting where offline policy optimization, temporal dependence, possible near-ties among actions, and valid uncertainty quantification interact.

The majority of existing statistical analyses of OPE concentrate on the classical setting in which the evaluation policy is fixed and known \textit{a priori}. 
In this regime, an extensive body of literature has established doubly robust and semiparametrically efficient estimators under various modeling assumptions \citep{jiang2016doubly,kallus2020double,shi2021deeply}. 
However, in many empirical applications, the policy of interest is not pre-specified but is itself estimated from the data as an \emph{optimal} policy. 
This setting introduces a qualitatively different statistical structure:
the maximization can make the optimal value non-smooth when competing
optimal policies are distinguishable to first order, while deterministic
or highly concentrated target policies can create a separate overlap
and numerical-stability problem.

Analogous issues have been extensively studied in the causal inference literature regarding optimal treatment regimes \citep{laber2014dynamic,kosorok2019precision,athey2021policy}.  There, ties typically produce non-regularity when the tied rules are distinguishable to first order; degenerate tie strata with identical first-order behavior are an important exception \citep{luedtke2016statistical}.
Extending such insights to the sequential decision-making framework of Markov decision processes (MDPs) is substantially more challenging due to temporal dependence, the Bellman fixed-point structure, and the complex interaction between policy optimization and value estimation.

Recently, \citet{shi2022statistical} proposed the SAVE estimator for infinite-horizon value inference, establishing semiparametric efficiency for the value of an optimal policy under a linear $Q$-function model and a set of non-degeneracy conditions. 
While SAVE represents a significant step toward principled inference for optimal policy values, its theory relies on stringent structural and regularity assumptions. 
In particular, it requires (i) a low-dimensional linear approximation of the $Q$-function, and (ii) well-conditioned Bellman estimating equations under the target policy. 
When the optimal policy is unique and deterministic, or nearly so, the latter condition often fails: the feature covariance induced by the target policy becomes ill-conditioned, leading to numerical instability and the breakdown of the associated inference. 
Moreover, in such regimes, SAVE no longer admits a doubly robust representation and loses its efficiency guarantees; furthermore, no alternative valid confidence sets are provided once these non-degeneracy conditions are violated.

Our target below is the same infinite-horizon discounted value functional.
The distinction is in the sampling experiment: we observe $N$ independent
stationary behavior episodes, or stationary windows, each containing a fixed
number of transitions.  Under the stationarity condition stated below, this
finite observed window does not change the estimand into a finite-horizon
return; it only determines how many stationary transition samples each
episode contributes to the normalized episode score.  Thus our results are
best viewed as fixed-window, many-episode inference for an infinite-horizon
discounted MDP, rather than as the bidirectional long-trajectory asymptotics
studied by \citet{shi2022statistical}.

This paper develops a unified inferential framework for the value of optimal policies in MDPs that explicitly addresses such non-regular phenomena. 
Our contributions are threefold. 
\begin{itemize}
    \item First, in an auxiliary augmented transition-sampling
    experiment, we characterize the existence of the efficient
    influence function (EIF) for the optimal policy value and derive
    its explicit form when the optimal policy is unique.  More
    generally, under ties the EIF exists exactly when all optimal
    policies share the same fixed-policy EIF; distinguishable optimal
    policies lead to non-regularity.  For the actual i.i.d.-trajectory
    experiment we further compute the semiparametric efficiency bound
    and exhibit a uniformly weighted estimator that attains it under a
    unique optimum.

    \item Second, building on this characterization, we propose a novel \emph{\textit{N}onparametric \textit{S}equenti\textit{A}l \textit{V}alue \textit{E}valuation (NSAVE)}. NSAVE admits martingale-based inference and a doubly robust representation when its nuisance functions target the same sequentially estimated policy, while remaining well-defined in degenerate or near-degenerate regimes where existing methods become unstable.

    \item Third, we study a complementary sample-split softmax approximation.  We state explicit first-stage conditions under which it gives pointwise Gaussian inference, while clarifying that softmax applied to a hard optimal $Q$-function does not generally restore pathwise differentiability at ties.
    Finally, beyond pointwise inference for the optimal value, we also consider a post-selection inference formulation. When the optimal policy is not unique, rather than targeting a single value functional, we give a calibrated template for confidence sets covering the collection of values associated with data-dependent selected policies.  Projection supplies an always-valid asymptotic baseline under the joint Gaussian approximation, while sharper data-dependent calibrations require the stated two-step error guarantee.
\end{itemize}

Through theoretical analysis, simulations, and an application to the
Drink Less micro-randomized trial, we study the stability and coverage of NSAVE and
the pointwise smoothing procedure across regular and non-regular
regimes, including settings where the optimal policy is deterministic
or nearly deterministic.

The remainder of the paper is organized as follows. 
In Sections \ref{sec:setup} and \ref{sec:IFF_for_EIF}, we first
characterize the efficient influence function in the auxiliary
augmented transition-sampling experiment and establish the
non-regularity that arises under first-order distinguishable optimal
policies.
In Sections \ref{sec:est} and \ref{sec:inference}, we then develop the
proposed NSAVE estimator together with martingale-based trajectory
inference and stability properties.
Section~\ref{sec:alternative} introduces the smoothing-based approach and the associated post-selection confidence sets for handling non-unique optimal policies. 
The finite-sample performance of NSAVE, the smoothing approach, and sample-split benchmarks is investigated through extensive simulation studies in Section~\ref{sec:sim}.
In Section~\ref{sec:realdata}, we apply the proposed methods to the Drink Less mobile health micro-randomized trial and conduct off-policy inference for learned notification policies.
Finally, the last section concludes with a discussion of the implications, limitations, and directions for future research.

\section{Problem Formulation}\label{sec:setup}

\subsection{Data Generating Process and Parameter of Interest}

We consider observational data generated from a canonical Markov decision process (MDP).
At any given time $t$, let $(S_t, A_t, R_t)$ denote the state-action-reward triplet.
We observe an offline dataset 
$\big\{ O_{it}: 1 \leq i \leq N, 0 \leq t \leq T\big\}$ with 
$O_{it} = (S_{it}, A_{it}, R_{it},S_{i,t+1})$,
generated by a behavior policy $b(\cdot \mid S)$, where $i$ indexes the episode and $t$ indexes the time point.
For any \textbf{fixed} target policy $\pi(a \mid s)$, OPE generally aims to evaluate the mean return 
$
    \eta(\pi) = \E^{\sim \pi} \left[ \sum_{t = 0}^{+\infty} \gamma^t R_t \right]
$
and construct a valid confidence interval,
where $\E^{\sim \pi}$ denotes the expectation when the system follows policy $\pi$. 
Distinct from existing semiparametric studies on OPE that consider an arbitrary $\pi$, we focus on a specific target policy: the optimal policy $\pi^*$, which maximizes $\eta(\pi)$ over the set of all possible policies $\Pi$. Specifically, the parameter of interest is
\[
     \begin{aligned}
         \eta^* := \eta(\pi^*) = \E^{\sim \pi^*} \bigg[ \sum_{t = 0}^{+\infty} \gamma^t R_t \bigg]\quad \text{ such that } \quad  \pi^* = \arg\max_{\pi \in \Pi} \eta(\pi).
     \end{aligned}
\]

To ensure the value function is identifiable, we adopt standard assumptions in the OPE literature.
For simplicity, we use $f(x \mid y)$ to represent the conditional density of $X$ given $Y = y$.

\begin{assumption}[Data Structure \& Observations]\label{ass:data_obs}
    We observe $N$ i.i.d. episodes
    $O_i=(S_{i0},A_{i0},R_{i0},\ldots,
    S_{iT},A_{iT},R_{iT},S_{i,T+1})$ of a common fixed length
    $T+1$.  Within each episode,
    $(R_t,S_{t+1})\sim f(r_t,s_{t+1}\mid A_t,S_t)$ and
    $A_t\sim b(a_t\mid S_t)$.
\end{assumption}

\begin{assumption}[Markov, Conditional Independence, \& Time-Homogeneity]\label{ass:Mark} 
    Conditional on $(S_t,A_t)$, the pair $(R_t,S_{t+1})$ is
    independent of the prior history, and
    $A_t\mid S_t\sim b(\cdot\mid S_t)$.  The conditional law of
    $(R_t,S_{t+1})$ given $(S_t,A_t)$ and the behavior policy
    $b(a\mid s)$ do not depend on $t$.
\end{assumption}

\begin{assumption}[Stationary Behavior Episodes]\label{ass:stationary}
    Under the behavior policy, $S_t$ has the same marginal law
    $f_b$ for every $t=0,\ldots,T+1$.  The target value uses the
    common initial-state law $\mu_0=f_b$.  The marginal ratio
    $\omega^\pi$ in \eqref{eq:mar_ratio} is square-integrable under
    the stationary behavior state--action law for every policy
    considered below.
\end{assumption}

Assumptions \ref{ass:data_obs} and \ref{ass:Mark} identify
$\eta(\pi)$ for a fixed policy.  Assumption \ref{ass:stationary} is an
additional sampling condition used for the normalized episode scores
and the trajectory-level inference results; it is not needed merely
to define the value functional.
The finite observed length $T+1$ should therefore be read as the length
of the available behavior window, not as a truncation of the target return.
Indeed, under Assumption~\ref{ass:stationary}, if $h(S,A,R,S')$ is any
integrable transition-level term with common stationary expectation
$\E h(S_t,A_t,R_t,S_{t+1})=c$, then
$
    \E\!\left[
    \frac{1}{1-\gamma^{T+1}}\sum_{t=0}^{T}\gamma^t
    h(S_t,A_t,R_t,S_{t+1})
    \right]
    =
    \frac{c}{1-\gamma}.
$
This identity is the reason the normalized finite-window scores used below
can target the infinite-horizon discounted value without truncation bias.
If one instead wished to use one or a few dependent long trajectories, the
same intuition would require additional mixing or blocking assumptions; that
long-trajectory asymptotic regime is not the statistical experiment analyzed
in this paper.  Accordingly, all asymptotic statements below are
many-episode results: $N\to\infty$ with any chosen finite window length
$T$ held fixed.  Larger independent stationary windows can still be
beneficial without changing the estimand; in the efficient
uniform-weight calculation of Section~\ref{sec:traj_eff}, the
transition-residual component of the bound is proportional to
$(T+1)^{-1}$.  A joint regime with $T=T_N\to\infty$ is a natural
extension, but it would require the Lindeberg and nuisance-rate
conditions to be verified uniformly in $T$, and is left outside the
present scope.

We briefly review standard estimation methods. 
The first method involves analyzing the aggregate mean return via the $Q$-function, defined as
\begin{equation}\label{eq:Q_fun}
    Q(a, s; \pi) := \E^{\sim \pi} \bigg[ \sum_{k = 0}^{+ \infty} \gamma^k R_{t + k} \mid A_t = a, S_t = s\bigg].
\end{equation}
The value function can be expressed as
\[
    \begin{aligned}
        & \eta(\pi) = \E^{\sim \pi} \Bigg[ \E^{\sim \pi} \bigg[ \sum_{t = 0}^{+\infty} \gamma^t R_t \mid A_0, S_0 \bigg] \Bigg] = \int Q(a_0, s_0 ; \pi) \pi(a_0 \mid s_0) f(s_0) \mathrm{d} a_0 \, \mathrm{d} s_0, \\
    \end{aligned}
\]
where $f(s_0)$ is the initial state density. We also define the value function $V(s; \pi) =  \int Q(a, s; \pi) \pi(a \mid s) \, \mathrm{d} a$.

The second method is the marginal importance sampling (MIS) estimator, which addresses the curse of horizon. The marginal ratio is defined as
\begin{equation}\label{eq:mar_ratio}
    \omega(a, s ; \pi) := (1 - \gamma) \sum_{t = 0}^{+ \infty}  \frac{\gamma^t f_{\sim \pi, t}(s)}{f_{+ \infty} (s)} \frac{\pi(a \mid s)}{b(a \mid s)} = (1 - \gamma) \sum_{t = 0}^{+ \infty}  \frac{\gamma^t f_{\sim \pi, t}(a, s)}{f_{+ \infty} (a, s)},
\end{equation}
where $f_{\sim \pi, t}$ and $f_{+ \infty}$ denote the time-dependent and stationary densities, respectively. Under stationarity, we have the identity 
\[
    \eta(\pi)
    =\frac{1}{1-\gamma}\E\!\left[\omega(A,S;\pi)R\right]
    =\frac{1}{1-\gamma}
    \E\!\left[\omega(A,S;\pi)r(A,S)\right],
\]
where $r(a,s)=\E(R\mid A=a,S=s)$ and the expectation is under the
stationary behavior-law transition distribution.
Crucially, both the $Q$-function \eqref{eq:Q_fun} and MIS ratio
\eqref{eq:mar_ratio} enter the efficient score in the auxiliary
transition-sampling experiment introduced below.

\subsection{Characterization and Issues}

Let $P$ index the underlying MDP data-generating components.  The same
$P$ induces both a full-trajectory law and, below, an auxiliary
transition-sampling law; these are different statistical experiments.
For any stationary policy $\pi$, define
\[
    \Psi(P;\pi)
    :=\E_{S_0\sim\mu_P}\!\left[V_P^\pi(S_0)\right],
    \qquad
    \Pi^*(P):=\arg\max_{\pi\in\Pi}\Psi(P;\pi),
\]
where $\mu_P$ is the initial-state law.  The optimized-value
functional is
\[
    \Psi^*(P):=\sup_{\pi\in\Pi}\Psi(P;\pi).
\]
When a representative is needed, let
$\pi^*(P)\in\Pi^*(P)$ be any measurable selection.  Thus
$\Psi^*(P_0)=\eta(\pi^*)$ is well-defined even if the optimizing
policy is not unique. We focus on the value functional; discussions
regarding the policy itself can be found in
\cite{kosorok2019precision,athey2021policy,luo2024policy}.
While $\Psi\big(P; \pi^*(P)\big) = \Psi^*(P)$, in general $\Psi\big(P_1; \pi^*(P_2)\big) \neq \Psi^*(P_1)$ if $P_1 \neq P_2$.
To keep population objects distinct from estimated nuisance functions, we henceforth write
\[
    Q_P^\pi(a,s) := Q(P)(a,s;\pi), \qquad
    V_P^\pi(s) := \int Q_P^\pi(a,s)\pi(a\mid s)\,\mathrm{d}a,
\]
and reserve hats, as in $\widehat Q^\pi$, for learned functions.  We use
$Q_P^*(a,s) := \sup_{\pi\in\Pi}Q_P^\pi(a,s)$ only for the optimal
$Q$-function; it should not be confused with $Q_P^\pi$ evaluated at a
generic or estimated policy.%

For the greedy-policy regret and smoothing results below, Assumption
\ref{ass:bellman_complete} imposes the standard Bellman-completeness
condition that $\Pi$ contains the stationary Markov greedy policies
under consideration.  In that case,
$Q_P^*$ is the unique Bellman-optimal action-value function,
$V_P^*(s)=\max_aQ_P^*(a,s)$, and every policy greedy with respect to
$Q_P^*$ is optimal.  The auxiliary envelope results can still be read
for a generic closed policy class, but the later greedy-policy
arguments use this stronger structure.

To separate the pathwise-differentiability calculation from the
actual trajectory-sampling problem, consider the auxiliary
\emph{augmented transition-sampling experiment}
$\mathcal M^\dagger$.  One observation is
$O^\dagger=(S_0,S,A,R,S')$, where $S_0\sim\mu_P$ is independent of
an i.i.d. behavior-law transition $(S,A,R,S')$.  At the truth, the
transition's state marginal is $f_b$, but the auxiliary experiment
does not observe a complete episode.  The primitive MDP components
$(\mu_P,r_P,p_P)$ determine the same fixed-policy and optimized values
as in the trajectory problem.
Write $P_0^\dagger$ for the law of $O^\dagger$ and
$P_{0,b}^{\mathrm{tr}}$ for its behavior-transition marginal.

All EIF and pathwise-differentiability statements in
Section~\ref{sec:IFF_for_EIF} are relative to
$\mathcal M^\dagger$.  They are used to understand the nonsmooth
optimization map, not to identify the efficiency bound for i.i.d.
episodes.  As shown in
\cite{uehara2020minimax,shi2024off,shi2021deeply}, for any fixed
$\pi$ the fixed-policy EIF in $\mathcal M^\dagger$ is
\begin{equation}\label{eq:EIF_for_fixed_policy}
    \begin{aligned}
        & \EIF \{ \Psi(P; \pi)\}\Big|_{P = P_0}(O^\dagger) \\
        = & \frac{1}{1 - \gamma}  \omega(P_0)(A, S; \pi)\big[ R + \gamma V(P_0)(S'; \pi) - Q(P_0)(A, S; \pi)\big] + V(P_0)(S_0; \pi) - \Psi(P_0; \pi), \\
        \text{or } \quad & \EIF_{\eta(\pi)}(O^\dagger; Q, \omega, V, \pi) \\
        =& \frac{1}{1 - \gamma} \omega(P_0)(A, S; \pi) \big[ R + \gamma V(S'; \pi) - Q(A, S; \pi)\big] + V(S_0; \pi) - \eta(\pi),
    \end{aligned}
\end{equation}
If the initial law is known rather than sampled, the
$V_P^\pi(S_0)-\Psi(P;\pi)$ component is omitted.
For $N$ independent full trajectories, the transition component of
the same expression supplies an unbiased score under Assumption
\ref{ass:stationary}, but it is not asserted to be the canonical
gradient of the trajectory experiment.
We denote the estimating functions for a single point $O$ and the trajectory $O_{0:T}$ as
\begin{equation}\label{eq:debiased_terms_for_fixed_policy}
    \begin{aligned}
        & \EstFunPoi_{\eta(\pi)}(O^\dagger; Q, \omega, V, \pi) = \frac{1}{1 - \gamma}  \omega(A, S; \pi) \big[ R + \gamma V(S'; \pi) - Q(A, S; \pi)\big] + V(S_0; \pi) \\
        & \EstFunTraj_{\eta(\pi)}(O_{0:T}; Q, \omega, V, \pi)
        = \frac{1}{1-\gamma^{T+1}}\sum_{t = 0}^T \gamma^t
        \omega(A_t, S_t; \pi)
        \big[ R_t + \gamma V(S_{t + 1}; \pi)
        - Q(A_t, S_t; \pi)\big] + V(S_0; \pi).
    \end{aligned}
\end{equation}
The factor $(1-\gamma^{T+1})^{-1}$ is required because
$\sum_{t=0}^T\gamma^t=(1-\gamma^{T+1})/(1-\gamma)$.
Under stationarity, it makes the trajectory estimating function have
the same expectation as the point-level estimating function.  Without
this normalization, the finite-$T$ score has a non-vanishing
truncation/scale bias even when the density ratio is correct.
The EIF for $\eta(\pi)$ satisfies
\[
     \EIF_{\eta(\pi)}(O^\dagger; Q, \omega, V, \pi)
     =\EstFunPoi_{\eta(\pi)}(O^\dagger;Q,\omega,V,\pi)-\eta(\pi),
\]
and, at the level of expectations, under stationarity,
\[
    \E \big[ \EIF_{\eta(\pi)}(O^\dagger; Q, \omega, V, \pi) \big]
    =\E \big[ \EstFunTraj_{\eta(\pi)}(O_{0:T}; Q, \omega, V, \pi) \big]-\eta(\pi).
\]
For later use, define the centered normalized episode score
\begin{equation}\label{eq:centered_episode_score}
    \phi_T(O_{0:T};\pi)
    :=
    \EstFunTraj_{\eta(\pi)}
    (O_{0:T};Q^\pi,\omega^\pi,V^\pi,\pi)-\eta(\pi).
\end{equation}

Within $\mathcal M^\dagger$, any regular asymptotically linear (RAL)
estimator $\widehat{\Psi}(\widehat{P};\pi)$ has an influence-function
representation \citep[IF,][]{tsiatis2006semiparametric}:
\[
    \widehat{\Psi}(\widehat{P}; \pi) - \Psi(P_0; \pi)
    =\frac{1}{N}\sum_{i=1}^{N}
    \operatorname{IF}(P_0^\dagger;\pi)(O_i^\dagger)
    +o_{P_0^\dagger}(N^{-1/2}).
\]
The EIF $\EIF \{ \Psi(P; \pi)\}\big|_{P = P_0}(O^\dagger)$ is the unique influence function that minimizes the variance $\var_{P_0^\dagger} \big( \operatorname{IF}(P_0^\dagger; \pi)(O^\dagger) \big)$.

Geometrically, the EIF is characterized via the auxiliary-model
tangent space $\mathcal{T}^\dagger$. For any differentiable path
$\{P_{\epsilon}:\epsilon\in\mathbb R\}\subset\mathcal M^\dagger$
passing through $P_0$ at $\epsilon=0$, pathwise differentiability of
$\Psi(P;\pi)$ implies
$
    \frac{\mathrm{d}}{\mathrm{d} \epsilon} \Psi(P_{\epsilon}; \pi) \Big|_{\epsilon = 0} = \E_{P_0^\dagger} \Big[ \widetilde{\psi}(O^\dagger) \dot{\ell}^\dagger (O^\dagger) \Big]
$
by the Riesz representation theorem, where $\dot{\ell}^\dagger(O^\dagger)$ is the score function. When $\widetilde{\psi}(O^\dagger) \in \mathcal{T}^\dagger$, then $\widetilde{\psi}(O^\dagger) = \EIF \{ \Psi(P; \pi)\}\big|_{P = P_0}(O^\dagger)$, and
\begin{equation}\label{eq:path_diff}
    \frac{\mathrm{d}}{\mathrm{d} \epsilon} \Psi(P_{\epsilon}; \pi) \Big|_{\epsilon = 0} = \E_{P_0^\dagger} \Big[ \EIF \{ \Psi(P; \pi)\}\big|_{P = P_0}(O^\dagger) \dot{\ell}^\dagger (O^\dagger) \Big].
\end{equation}

However, the policy map $P\mapsto\pi^*(P)$ need not itself be
differentiable, so a formal chain rule through $\pi^*(P)$ is generally
not justified.  The appropriate tool is an envelope argument: first
obtain a local expansion of $\Psi(P;\pi)$ that is uniform over $\pi$,
and then optimize that expansion.  This yields a directional derivative
determined by the fixed-policy EIFs of the policies that are optimal at
$P_0$.  Section~\ref{sec:IFF_for_EIF} makes this statement precise.

\section{Pathwise Differentiability in an Auxiliary Experiment}
\label{sec:IFF_for_EIF}

For a fixed policy $\pi$, let
\[
    D_{P_0}^{\dagger,\pi}(O^\dagger)
    := \EIF\{\Psi(P;\pi)\}\big|_{P=P_0}(O^\dagger),
    \qquad
    \Pi^*(P_0):=\arg\max_{\pi\in\Pi}\Psi(P_0;\pi).
\]
For the policy metric used below, take
\[
    d_\Pi(\pi_1,\pi_2)
    :=\sup_{s\in\mathcal S}
    \operatorname{TV}\{\pi_1(\cdot\mid s),\pi_2(\cdot\mid s)\}.
\]
The following condition states explicitly the uniformity needed to pass
from fixed-policy differentiability to differentiability of the
optimized value.

\begin{assumption}[Finite-MDP Local Regularity in
$\mathcal M^\dagger$]\label{ass:regularity}
    The state and action spaces are finite,
    $|R|\leq\overline c_R<\infty$, $\gamma<1$, and $\Pi$ is a closed
    class of stationary
    Markov policies.  There is a constant $c>0$ such that the
    stationary behavior state law satisfies $f_b(s)\geq c$ for every
    state reachable under a policy in $\Pi$, and
    $b(a\mid s)\geq c$ whenever $\pi(a\mid s)>0$ for some
    $\pi\in\Pi$.
    For every regular parametric submodel
    $\{P_\epsilon:\epsilon\in(-\epsilon_0,\epsilon_0)\}$ through $P_0$,
    the initial law $\mu_\epsilon$, conditional reward mean
    $r_\epsilon(s,a)$, and transition kernel
    $p_\epsilon(s'\mid s,a)$ satisfy
    \[
        \begin{aligned}
        \|\mu_\epsilon-\mu_0-\epsilon\dot\mu\|_1
        =o(|\epsilon|),\qquad
        \max_{s,a}|r_\epsilon(s,a)-r_0(s,a)
        -\epsilon\dot r(s,a)|
        & =o(|\epsilon|),\\
        \text{and} \qquad \max_{s,a}\sum_{s'}
        |p_\epsilon(s'\mid s,a)-p_0(s'\mid s,a)
        -\epsilon\dot p(s'\mid s,a)|
        &=o(|\epsilon|).
        \end{aligned}
    \]
\end{assumption}

Assumption \ref{ass:regularity} is stated in terms of the primitive MDP
components.  
Lemma \ref{lem:finite_mdp_uniform_expansion} below shows
that it implies, rather than assumes, both the uniform local expansion
and its representation by the fixed-policy canonical gradient in
$\mathcal M^\dagger$.
It does not identify a trajectory-level
canonical gradient or efficiency bound.

The next assumption is imposed only for the regular-case theorem
below.  It is not maintained in the nonregularity result or in the
one-sided NSAVE analysis.
\begin{assumption}[Unique Optimal Policy]\label{ass:deterministic_policy}
    The optimal set is a singleton:
    $\Pi^*(P_0)=\{\pi_0^*\}$.
    The policy $\pi_0^*$ may be deterministic when the action space and
    overlap conditions permit this, but determinism is not required for
    the differentiability conclusion.
\end{assumption}

The overlap conditions in Assumption \ref{ass:regularity} restrict only the
\emph{behavior} law $(f_b,b)$ that generates the data; they place no
probabilistic restriction on the evaluated policies and are therefore logically
independent of Assumption~\ref{ass:deterministic_policy}, which constrains only
the \emph{optimal set} $\Pi^*(P_0)$ and may hold for a stochastic optimizer.

\begin{theorem}[Envelope characterization under uniqueness in
$\mathcal M^\dagger$]
\label{thm:suff_cond_for_EIF_exp}
    Suppose that Assumptions \ref{ass:regularity} and
    \ref{ass:deterministic_policy} hold.
    Then $\Psi^*$ is pathwise differentiable at $P_0$ relative to
    $\mathcal M^\dagger$, and its EIF in that experiment is the
    fixed-policy EIF evaluated at the unique optimizer:
    \[
        \EIF\{\Psi^*(P)\}\big|_{P=P_0}
        =D_{P_0}^{\dagger,\pi_0^*}
        =\EIF\{\Psi(P;\pi)\}\big|_{P=P_0,\,\pi=\pi_0^*}.
    \]
\end{theorem}

Non-uniqueness alone is not sufficient for nonregularity.  If two
optimal policies induce the same first-order behavior, the optimal
value can remain pathwise differentiable.  The relevant condition is
whether optimal policies are distinguishable through their
fixed-policy EIFs.

\begin{corollary}[First-order equivalent tied optima]
\label{cor:equivalent_tied_optima}
    Under Assumption \ref{ass:regularity}, suppose there is a common
    $D_{P_0}^{\dagger,*}\in L_2^0(P_0^\dagger)$ such that
    \[
        D_{P_0}^{\dagger,\pi}
        =D_{P_0}^{\dagger,*}
        \qquad P_0^\dagger\text{-a.s. \quad for every } \quad
        \pi\in\Pi^*(P_0).
    \]
    Then $\Psi^*$ is pathwise differentiable at $P_0$ relative to
    $\mathcal M^\dagger$, with EIF
    $D_{P_0}^{\dagger,*}$, even if $\Pi^*(P_0)$ is not a singleton.
\end{corollary}

\begin{assumption}[First-Order Distinguishable Optimal Policies]
\label{ass:unres_rules}
    There exist $\pi_1,\pi_2\in\Pi^*(P_0)$ such that
    $\|D_{P_0}^{\dagger,\pi_1}
    -D_{P_0}^{\dagger,\pi_2}\|_{P_0^\dagger,2}>0$.
\end{assumption}

\begin{theorem}[Nonregularity under distinguishable optima in
$\mathcal M^\dagger$]
\label{thm:nec_cond_for_EIF}
    Suppose that Assumptions \ref{ass:regularity} and
    \ref{ass:unres_rules} hold in the full
    nonparametric auxiliary model $\mathcal M^\dagger$.  Then the
    right and left pathwise derivatives of $\Psi^*$ disagree along
    some regular submodel of $\mathcal M^\dagger$.  Consequently,
    $\Psi^*$ is not pathwise differentiable at $P_0$ relative to
    $\mathcal M^\dagger$ and has no influence function in that
    experiment.
\end{theorem}

Thus uniqueness is a transparent sufficient condition in
$\mathcal M^\dagger$.  Under non-uniqueness, the sharper
characterization in the full nonparametric auxiliary model is that an
EIF exists if all policies in $\Pi^*(P_0)$ share the same
fixed-policy EIF; otherwise the optimal value is only directionally
differentiable there.  This parallels the exceptional-law distinction
for optimal treatment rules \citep{luedtke2016statistical}.  No
trajectory-level efficiency conclusion is drawn from this
characterization.

\section{Estimation Under Possible Non-Uniqueness}\label{sec:est}

\subsection{The Challenge and Current Gap}

When $\pi$ is fixed, an i.i.d. sample from the auxiliary experiment
$\mathcal M^\dagger$ would yield the augmented-transition one-step
estimator by averaging its estimating function.  For the actual
i.i.d.-episode data, the trajectory version instead averages the
normalized episode score:
\begin{equation*}
    \begin{aligned}
        & \widehat{\eta}_{\text{DR, 1}}(\pi) & := \pnt \EstFunPoi_{\eta(\pi)}(O^\dagger;\widehat{Q},  \widehat{\omega},  \widehat{V},  \pi) \quad
        \text{or} \quad \widehat{\eta}_{\text{DR, 2}}(\pi) & := \pn\EstFunTraj_{\eta(\pi)}(O_{0:T};\widehat{Q},  \widehat{\omega},  \widehat{V},  \pi)
    \end{aligned}
\end{equation*}
with estimated nuisance functions $\widehat{Q},\widehat{\omega},
\widehat{V}$.  The two scores have the same expectation under
stationarity, although their variances need not coincide because
observations within a trajectory are dependent:
\begin{itemize}
    \item Double-Robustness: Assuming either $\widehat{Q}$ (and thus $\widehat{V}$) or $\widehat{\omega}$ is consistent, both estimating functions have the correct expectation.  Their corresponding empirical averages are consistent under the sampling and empirical-process conditions appropriate to their respective experiments.

    \item Semiparametric Efficiency: In the auxiliary augmented
    transition-sampling experiment, if sample splitting or
    cross-fitting is used, both nuisances are $L_2$-consistent, and
    their product error is $o_{P_0}(N^{-1/2})$ (for example, each is
    $o_{P_0}(N^{-1/4})$), the point-level one-step estimator
    achieves semiparametric efficiency, satisfying
    $
        \sqrt{N} \big( \widehat{\eta}_{\text{DR, 1}}(\pi) - \eta(\pi) \big)
        \rightsquigarrow
        \mathcal{N}\!\left(
        0,\,
        \E\!\left[
        \left\{
        \EIF_{\eta(\pi)}
        (O^\dagger;Q,\omega,V,\pi)\right\}^2\right]\right).
    $
    This auxiliary-experiment statement should not be obtained by
    treating dependent transitions within an observed episode as
    i.i.d.  The normalized trajectory estimator instead has the
    variance of its episode-level score.
\end{itemize}

When the optimal policy (or policies) $\pi^*$ is unknown, 
and we have an estimated optimal policy 
$\widehat{\pi}$ computed by
a $Q$-learning type algorithm such that 
\[
    \widehat{\pi}(a \mid s) := \mathds{1} \Big\{ a = \arg \max_{a' \in \mathcal{A}} \widehat{Q}_{\text{opt}} (a, s ) \Big\},
\]
where $\widehat{Q}_{\text{opt}} (a, s ) $ denotes a consistent estimator for the optimal $Q$-function, i.e., $Q(a, s; \pi^* )$,
an intuitive ``plug-in" estimator for the parameter of interest $\eta(\pi^*)$ is $\widehat{\eta}_{\text{DR, j}}(\widehat{\pi})$.  Conditional fixed-policy double robustness remains available, but inference for the data-adaptive policy additionally requires policy stability and control of its regret.
However, as shown in Theorem \ref{thm:nec_cond_for_EIF}, when multiple
optimal policies have distinct fixed-policy EIFs, there is no regular
EIF for the unsmoothed optimal value in
$\mathcal M^\dagger$.  Mere non-uniqueness without first-order
distinguishability does not imply this conclusion, and the theorem is
not a trajectory-efficiency statement.
More importantly, in such cases, the argmax of the optimal $Q$-function may not be uniquely defined; consequently, $\widehat{\pi}(\cdot \mid s)$ might not converge to a \textbf{fixed} quantity for some $s \in \mathcal{S}$. As a result, the plug-in estimator $\widehat{\eta}_{\text{DR, j}}(\widehat{\pi})$ will fluctuate randomly and fail to maintain a stable limiting distribution \citep{shi2022statistical}.

To overcome this issue, \cite{shi2022statistical} proposed a new estimator called \textit{S}equenti\textit{A}l \textit{V}alue \textit{E}valuation (SAVE), denoted as $\widehat{\eta}_{\text{SAVE}}$, assuming that the $Q$-function follows a linear sieve model such that $Q (a, s; \pi) \approx \Phi^{\top}(s) \bm{\beta}_{\pi, a}$, where $\Phi(s)$ is a vector of sieve basis functions. This novel estimator enjoys bidirectional asymptotic normality:
\[
    \begin{aligned}
        \sqrt{NT (K - 1) / K}  \widehat{\sigma}_{\text{SAVE}}^{-1}\big( \widehat{\eta}_{\text{SAVE}} - \eta (\widehat{\pi})\big) \, \rightsquigarrow \, \mathcal{N}(0, 1); \\
        \sqrt{NT (K - 1) / K} \widehat{\sigma}_{\text{SAVE}}^{-1}\big( \widehat{\eta}_{\text{SAVE}} - \eta ({\pi}^*)\big) \, \rightsquigarrow \, \mathcal{N}(0, 1),
    \end{aligned}
\]
as either $N \rightarrow \infty$ or $T \rightarrow \infty$, where $K$ is the number of data partitions and $\widehat{\sigma}_{\text{SAVE}}^2$ is a ``plug-in" type variance estimator. Thus, the readily applicable estimator $\widehat{\eta}_{\text{SAVE}}$ can be used for statistical inference.

Nonetheless, despite its appealing theoretical guarantees, the SAVE estimator $\widehat{\eta}_{\text{SAVE}}$ suffers from several significant limitations. First, SAVE relies critically on a linear structural assumption for the $Q$-function, namely that $Q(s,a;\pi)$ can be well-approximated by a low-dimensional linear form $Q(s,a;\pi)\approx \Phi^\top(s)\bs{\beta}_{\pi,a}$. This assumption may be violated in many realistic sequential decision problems. Second, when the optimal policy is unique and deterministic, $\widehat{\eta}_{\text{SAVE}}$ no longer admits a doubly robust representation and consequently loses both the double robustness property and the associated semiparametric efficiency guarantees. Third, SAVE requires non-degeneracy conditions on the target-policy Bellman equations, which can become difficult to verify for deterministic or highly concentrated target policies. The simulations in Section~\ref{sec:sim} therefore focus directly on the properties established by our theory: fixed-policy double robustness, sequential inference under unique and tied optima, smoothing approximation bias, and sensitivity to weak overlap.

To address these three challenges, we adopt the conceptual framework of \cite{shi2022statistical} while introducing a revised sequential value evaluation procedure and a corresponding estimator, which we term \textit{N}onparametric \textit{S}equenti\textit{A}l \textit{V}alue \textit{E}valuation (NSAVE).

\subsection{Nonparametric SequentiAl Value Evaluation Approach}\label{sec:new_save}

Let $\tau$ be a random permutation independent of the i.i.d.
trajectory observations $\{O_i\}_{i=1}^N$.  Then
$\{O_{\tau(i)}\}_{i=1}^N$ is again an i.i.d. sequence, and
$O_{\tau(j)}$ is independent of
$\mathcal F_{j-1}:=\sigma\{O_{\tau(i)}:i\leq j-1\}$.
Let $\ell_N$ be the number of episodes in the initial training sample.
The first online policy, used to evaluate $O_{\tau(\ell_N+1)}$, is
therefore trained on
$\{O_{\tau(i)}:i\leq\ell_N\}$.

For $j = \ell_N + 1, \ldots, N$, we perform the following steps:
\begin{itemize}
    
    \item \textbf{Policy learning}: Using only
    $\{O_{\tau(i)}:i\leq j-1\}$, estimate the optimal $Q$-function,
    denoted by $\widehat Q_{\tau(j-1)}^*$.  With a fixed measurable
    tie-breaking map, let
    \[
        \widehat a_{\tau(j-1)}(s)
        :=\operatorname{TieBreak}\!\left(
        \arg\max_{a'\in\mathcal A}
        \widehat Q_{\tau(j-1)}^*(a',s)\right),
        \qquad
        \widehat{\pi}_{\tau(j - 1)}^{(Q)}(a \mid s)
        :=\mathds{1}\{a=\widehat a_{\tau(j-1)}(s)\}.
    \]

    \item \textbf{Policy-aligned nuisance training}: On the same past
    data, but treating
    $\widehat\pi_{\tau(j-1)}^{(Q)}$ as the evaluation policy, estimate
    \[
        \widehat Q_{\tau(j-1)}
        (\cdot,\cdot;\widehat\pi_{\tau(j-1)}^{(Q)}),
        \qquad
        \widehat\omega_{\tau(j-1)}
        (\cdot,\cdot;\widehat\pi_{\tau(j-1)}^{(Q)}).
    \]
    Set
    $\widehat V_{\tau(j-1)}(s;\widehat\pi)
      =\sum_a\widehat\pi(a\mid s)
       \widehat Q_{\tau(j-1)}(a,s;\widehat\pi)$.

    \item \textbf{Evaluating}: Using the above nuisance functions, we
    calculate the normalized \textit{episode} estimating function
    \[
        \begin{aligned}
            \EstestFunTraj_{\tau(j)}
            &:=\frac{1}{1-\gamma^{T+1}}\sum_{t = 0}^T \gamma^t \widehat{\omega}_{\tau(j - 1)}\big(A_{\tau(j), t}, S_{\tau(j), t}; \widehat{\pi}_{\tau(j - 1)}^{(Q)}\big) \big[ R_{\tau(j), t} + \gamma \widehat{V}_{\tau(j - 1)} \big(S_{\tau(j), t + 1}; \widehat{\pi}_{\tau(j - 1)}^{(Q)} \big)\\
            &~~~~~~~~~~~~~~  - \widehat{Q}_{\tau(j - 1)}  \big(A_{\tau(j), t}, S_{\tau(j), t}; \widehat{\pi}_{\tau(j - 1)}^{(Q)} \big)\big]  + \widehat{V}_{\tau(j - 1)} \big(S_{\tau(j), 0}; \widehat{\pi}_{\tau(j - 1)}^{(Q)} \big).
        \end{aligned}
    \]
\end{itemize}

For the two-sided result, after the online evaluations are complete,
we may refit the optimal $Q$-function on all $N$ episodes.  Denote
this fit by $\widehat Q_{\mathrm{fin}}^*$ and its fixed-tie-breaking
greedy policy by $\widehat\pi_{\mathrm{fin}}$.  This final refit is
not used inside any martingale increment; it only supplies the policy
value around which the two-sided decomposition is centered.

The density-ratio learner is therefore a fixed-policy nuisance learner
at each sequential step; it is not a second policy optimizer.  This
alignment is what permits the exact double-robustness identity in
Lemma~\ref{lem:exact_dr_identity}.%


Define the online one-step variance as 
\[
    \begin{aligned}
        \widetilde{\sigma}_{\tau(j - 1)}^2
        :=\var\!\left(
        \EstestFunTraj_{\tau(j)}
        \mid\mathcal F_{j-1}\right),
        \qquad j=\ell_N+1,\ldots,N.
    \end{aligned}
\]
Let $\widehat{\sigma}_{\tau(j - 1)}^2$ denote its corresponding consistent estimator. 
In practice, $\widehat{\sigma}_{\tau(j - 1)}^2$ can be estimated using
a sliding window such as
$\{\EstestFunTraj_{\tau(j-m)},\ldots,
\EstestFunTraj_{\tau(j-1)}\}$.  For asymptotic justification, a fixed
window is generally insufficient: one may take $m=m_N\to\infty$ while
requiring the conditional score variances and nuisance targets to
drift by $o_{P_0}(1)$ across the window.  Assumption
\ref{ass:est_var} below records the resulting variance-consistency
requirement directly.
Then, our final estimator is similarly defined as the weighted average:
\[
    \begin{aligned}
        \widehat{\eta}_{\text{NSAVE}} := \bigg\{ \sum_{j = \ell_N + 1}^N \frac{1}{\widehat{\sigma}_{\tau(j - 1)}}\bigg\}^{-1} \sum_{j = \ell_N + 1}^N \frac{\EstestFunTraj_{\tau(j)}}{\widehat{\sigma}_{\tau(j - 1)}}.
    \end{aligned}
\]
The inverse-standard-deviation weights studentize each predictable
martingale increment, giving it approximately unit conditional
variance and leading to the transparent martingale normalization in
Theorem~\ref{thm:CLT_for_R1_CLB}.  They are not the
inverse-\emph{variance} weights that minimize the variance of a
weighted average of independent unbiased observations.  Accordingly,
the weighting is motivated by sequential stabilization, not by a
claim of semiparametric efficiency.%

Intuitively, our novel estimator $\widehat{\eta}_{\text{NSAVE}}$ approximates, but is distinct from, the average \textit{weighted} empirical historical value $\overline{\eta}_w(\widehat{\pi}^{(Q)})$, defined as
\[
    \overline{\eta}_{w}(\widehat{\pi}^{(Q)})
    :=
    \bigg\{\sum_{j=\ell_N+1}^N
    \frac{1}{\widehat{\sigma}_{\tau(j-1)}}\bigg\}^{-1}
    \sum_{j=\ell_N+1}^N
    \frac{\eta(\widehat{\pi}_{\tau(j-1)}^{(Q)})}
    {\widehat{\sigma}_{\tau(j-1)}}.
\]

In the following, we analyze the theoretical properties of our novel estimator $\widehat{\eta}_{\text{NSAVE}}$ to demonstrate its advantages.
\begin{assumption}[Trajectory-Score Regularity]
\label{ass:traj_score_reg}
    The state and action spaces are finite,
    $|R_t|\leq\overline c_R<\infty$, and $\gamma<1$.  The stationary
    behavior state law is
    bounded away from zero on states reachable under the evaluated
    policies, and $b(a\mid s)$ is bounded away from zero whenever an
    evaluated policy assigns positive probability to action $a$ at
    state $s$.  The fitted $Q$-functions and
    marginal ratios used by NSAVE are predictable with respect to
    $\mathcal F_{j-1}$ and are uniformly bounded with probability
    tending to one.  Finally, $n_N=N-\ell_N\to\infty$.
\end{assumption}

This condition concerns the moments of the actual episode scores.  It
is separate from the parametric-submodel conditions used to calculate
the auxiliary-experiment EIF.

\subsubsection{Nuisance Estimation Approaches}

There are various approaches for estimating the nuisance components.
At step $j$, the density-ratio target is
$\omega(\cdot,\cdot;\widehat\pi_j)$ for the already selected
evaluation policy $\widehat\pi_j=\widehat\pi_{\tau(j-1)}^{(Q)}$.
It is not obtained by maximizing reward over occupancy measures.

For example, one may impose the discounted flow constraints
\[
    \begin{aligned}
        & \sum_{a} \omega(a, s') b(a \mid s') f_b(s') = (1 - \gamma) \mu_0(s') + \gamma \sum_{a, s}
        f(s' \mid a, s) \omega(a, s) b(a \mid s) f_b(s),
        \qquad s'\in\mathcal S.
    \end{aligned}
\]
Together with policy compatibility,
\[
    \frac{\omega(a,s)b(a\mid s)}
    {\sum_{a'}\omega(a',s)b(a'\mid s)}
    =\widehat\pi_j(a\mid s),
\]
these equations identify the ratio for $\widehat\pi_j$.  Equivalently,
one can estimate the state-density ratio and set
$\omega(a,s;\widehat\pi_j)
 =\rho(s;\widehat\pi_j)\widehat\pi_j(a\mid s)/b(a\mid s)$.

Another option is fixed-policy Minimax Weight Learning
\citep{nachum2019dualdice,duan2020minimax,uehara2020minimax}.  For a
fixed $\widehat\pi_j$, define
\[
    \begin{aligned}
        \mathcal{L} (Q^{\pi}, \omega^{\pi}; P_0) & : = (1 - \gamma) \E_{\mu_0} \big[ V^\pi(S)\big] \\
        &~~~~~~~~~~~~~ + \E_{P_0} \Big[ \omega(A, S ; \pi)
        \Big\{ R + \gamma V^\pi(S') - Q^\pi(A,S) \Big\} \Big].
    \end{aligned}
\]
Let $\mathbb{P}_{\tau(j - 1) T} :=  \frac{1}{T (j - 1)} \sum_{t = 1}^T \sum_{\iota = 1}^{j - 1} [\bcdot]_{\tau(\iota), t}$ be the empirical distribution measure at Step $j$.
We can then construct policy-aligned nuisance estimators by solving
\begin{equation}\label{eq:empirical_lag_for_nui}
    \begin{aligned}
        \big( \widehat{Q}_{j}^{\widehat\pi_j},
        \widehat{\omega}_{j}^{\widehat\pi_j} \big)
        = \arg\min_{\omega \in \widehat{\Omega}_{\text{flow}}(\widehat\pi_j)}
        \max_{Q \in \mathcal{Q}}
        \mathcal{L} \big(Q^{\widehat\pi_j},
        \omega^{\widehat\pi_j};
        \mathbb{P}_{\tau(j - 1) T} \big).
    \end{aligned}
\end{equation}
An unconstrained reward-maximizing occupancy program instead learns a
second policy and is not covered by the double-robustness theorem.

\section{Inference}\label{sec:inference}

Although it is intuitively plausible that our novel estimator $\widehat{\eta}_{\text{NSAVE}}$ will be consistent as long as the nuisance components are consistently estimated (we will also formally demonstrate its consistency),
such intuition is insufficient for inference. 
The latter typically requires stronger conditions. To characterize the specific requirements, consider that the remainder term can be decomposed into two distinct components corresponding to two different inferential strategies:

\begin{itemize}
    \item For the \textit{Conservative Lower Bound}: 
    \[
        \begin{aligned}
            \widehat{\eta}_{\text{NSAVE}} - \eta(\pi^*) = \underbrace{\widehat{\eta}_{\text{NSAVE}} - \overline{\eta}_w(\widehat{\pi}^{(Q)})}_{=: R_{\text{CLB}, 1N}} + \underbrace{\overline{\eta}_w(\widehat{\pi}^{(Q)}) - \eta(\pi^*)}_{=: R_{\text{CLB}, 2N}}.
        \end{aligned}
    \]
    This is a relatively coarse decomposition, as our primary goal here is to establish a lower bound. It is straightforward to see that $R_{\text{CLB}, 1N}$ represents the empirical error for the average value functions under the estimated policies, while $R_{\text{CLB}, 2N}$ represents the cumulative regret arising from the estimated policies.

    \item For the \textit{Two-Sided Confidence Interval}:
    \[
        \begin{aligned}
            \widehat{\eta}_{\text{NSAVE}}-\eta(\pi^*)
            =
            \underbrace{\widehat{\eta}_{\text{NSAVE}}
            -\eta(\widehat\pi_{\mathrm{fin}})}
            _{=:R_{\text{TCI},1N}}
            +
            \underbrace{\eta(\widehat\pi_{\mathrm{fin}})
            -\eta(\pi^*)}_{=:R_{\text{TCI},2N}}.
        \end{aligned}
    \]
    Here, we use a more refined decomposition consistent with standard analyses:  $R_{\text{TCI}, 1N}$ represents the statistical error, and $R_{\text{TCI}, 2N}$ captures the policy-value error. 
\end{itemize}

\subsection{Conservative Lower Bound}\label{sec:lower_bound}

In both decompositions, the second terms, $R_{\text{CLB}, 2N}$ and $R_{\text{TCI}, 2N}$, are non-positive by the definition of $\pi^*$. Consequently, we have 
\[
    \begin{aligned}
        \eta(\pi^*) = \left\{ \begin{array}{cc}
            \widehat{\eta}_{\text{NSAVE}} -  \big( R_{\text{CLB}, 1N} + R_{\text{CLB}, 2N} \big)  &  \geq \widehat{\eta}_{\text{NSAVE}} - R_{\text{CLB}, 1N}  \\
            \widehat{\eta}_{\text{NSAVE}} -  \big( R_{\text{TCI}, 1N} + R_{\text{TCI}, 2N} \big) &  \geq \widehat{\eta}_{\text{NSAVE}} - R_{\text{TCI}, 1N} 
        \end{array}\right. .
    \end{aligned}
\]
If we can construct a valid $(1 - \alpha)$ upper bound $\UB(R_{1N}; \alpha)$ for either $R_{\text{CLB}, 1N}$ or $R_{\text{TCI}, 1N} $ such that
\[
     \liminf_{N \rightarrow \infty} \pr \big(R_{\text{CLB}, 1N} \leq \UB(R_{1N}; \alpha) \big) \geq 1 - \alpha ~~\text{ or }  ~~\liminf_{N \rightarrow \infty} \pr \big( R_{\text{TCI}, 1N} \leq \UB(R_{1N}; \alpha) \big) \geq 1 - \alpha,
\]
then 
\[
    \begin{aligned}
        \liminf_{N \rightarrow \infty} \pr \big( \eta(\pi^*) \geq \widehat{\eta}_{\text{NSAVE}} - \UB(R_{1N}; \alpha) \big) \geq 1 - \alpha.
    \end{aligned}
\]
This implies that $\widehat{\eta}_{\text{NSAVE}} - \UB(R_{1N}; \alpha)$ serves as a valid lower bound for the optimal value. The following theorem formally states how to construct a valid $\UB(R_{1N}; \alpha)$ and its corresponding estimator $\widehat{\UB}(R_{1N}; \alpha)$. Let
\[
    \begin{aligned}
    \sigma_{R_{1N}}
    &:=
    \sqrt{N-\ell_N}
    \left\{\sum_{j=\ell_N+1}^N
    \frac{1}{\widetilde\sigma_{\tau(j-1)}}\right\}^{-1},\\
    \widehat\sigma_{R_{1N}}
    &:=
    \sqrt{N-\ell_N}
    \left\{\sum_{j=\ell_N+1}^N
    \frac{1}{\widehat\sigma_{\tau(j-1)}}\right\}^{-1},
    \end{aligned}
\]
the standard-error scale of the inverse-standard-deviation weighted
average,
and write
\[
    \widehat\phi_j:=\EstestFunTraj_{\tau(j)},
    \qquad
    M_j:=\widehat\phi_j-
    \E(\widehat\phi_j\mid\mathcal F_{j-1}).
\]
The upcoming theorem, which establishes the asymptotic normality for the first terms in the two types of decompositions, relies on the following assumptions:

\begin{assumption}[Convergence Rates for Nuisance Parameters]\label{ass:nui_rates}
    Let $\widehat\pi_j:=\widehat{\pi}_{\tau(j-1)}^{(Q)}$.
    Conditional on $\mathcal F_{j-1}$, the nuisance estimators target
    the same evaluation policy, i.e.,
    \[
        \widehat Q_j \equiv Q^{\widehat\pi_j},
        \qquad
        \widehat\omega_j \equiv \omega^{\widehat\pi_j}.
    \]
    Here and below, the displayed $L_2(P_0)$ nuisance norms are taken
    under the stationary behavior transition marginal.
    With $n_N=N-\ell_N\to\infty$,
    \[
        \frac{1}{\sqrt{n_N}}
        \sum_{j=\ell_N+1}^N
        \|\widehat Q_j-Q^{\widehat\pi_j}\|_{P_0,2}
        \|\widehat\omega_j-\omega^{\widehat\pi_j}\|_{P_0,2}
        =o_{P_0}(1).
    \]
    A sufficient uniform sequential rate condition is
    \[
        \begin{aligned}
        \max_{\ell_N<j\leq N}
        (j-\ell_N)^{\kappa_Q}
        \|\widehat Q_j-Q^{\widehat\pi_j}\|_{P_0,2}
        &=O_{P_0}(1),\\
        \max_{\ell_N<j\leq N}
        (j-\ell_N)^{\kappa_\omega}
        \|\widehat\omega_j-\omega^{\widehat\pi_j}\|_{P_0,2}
        &=O_{P_0}(1),
        \end{aligned}
    \]
    with $\kappa_Q+\kappa_\omega>1/2$.
\end{assumption}

To see why the displayed uniform rates are sufficient, put
$r=\kappa_Q+\kappa_\omega>1/2$.  On events with probability tending
to one, the product at online index
$k=j-\ell_N$ is bounded by a random $O_{P_0}(1)$ constant times
$k^{-r}$.  Hence
\[
    \frac1{\sqrt{n_N}}\sum_{k=1}^{n_N}k^{-r}
    =
    \begin{cases}
    O(n_N^{1/2-r}),&1/2<r<1,\\
    O(n_N^{-1/2}\log n_N),&r=1,\\
    O(n_N^{-1/2}),&r>1,
    \end{cases}
\]
and every case converges to zero.

If the implementation estimates an occupancy ratio through a separate
policy $\widehat\pi_j^{(\omega)}$, the assumption requires its resulting
ratio to converge to $\omega^{\widehat\pi_j}$; convergence merely to
the ratio of a different policy does not give the usual
double-robustness identity.

\begin{assumption}[Nondegenerate Predictable Scales]
\label{ass:non_zero_var}
    There are constants $0<\underline\sigma<\overline\sigma<\infty$
    such that $\widehat\sigma_{\tau(j-1)}$ and
    $\widetilde\sigma_{\tau(j-1)}$ are
    $\mathcal F_{j-1}$-measurable and, with probability tending to
    one,
    \[
        \underline\sigma
        \leq
        \inf_{\ell_N<j\leq N}
        \{\widehat\sigma_{\tau(j-1)}
        \wedge\widetilde\sigma_{\tau(j-1)}\}
        \leq
        \sup_{\ell_N<j\leq N}
        \{\widehat\sigma_{\tau(j-1)}
        \vee\widetilde\sigma_{\tau(j-1)}\}
        \leq\overline\sigma.
    \]
\end{assumption}

\begin{assumption}[Conditions for Estimated Variances]\label{ass:est_var}
    \[
        \frac{1}{N-\ell_N}
        \sum_{j=\ell_N+1}^N
        \left(
        \frac{\widehat\sigma_{\tau(j-1)}}
        {\widetilde\sigma_{\tau(j-1)}}-1
        \right)^2
        =o_{P_0}(1).
    \]
\end{assumption}

\begin{assumption}[Lindeberg Condition]\label{ass:Lindeberg}
    For any $\epsilon > 0$,
    \[
        \frac{1}{N-\ell_N}
        \sum_{j=\ell_N+1}^N
        \E\!\left[
        \left.
        \frac{M_j^2}{\widetilde\sigma_{\tau(j-1)}^2}
        \mathds{1}\!\left\{
        \frac{|M_j|}
        {\sqrt{N-\ell_N}\widetilde\sigma_{\tau(j-1)}}
        >\epsilon\right\}
        \right|\mathcal F_{j-1}\right]
        =o_{P_0}(1).
    \]
\end{assumption}

\begin{theorem}\label{thm:CLT_for_R1_CLB}
    Suppose that Assumptions \ref{ass:data_obs},
    \ref{ass:Mark}, \ref{ass:stationary},
    \ref{ass:traj_score_reg}, and
    \ref{ass:nui_rates}--\ref{ass:Lindeberg} hold. Then
    \[
        \begin{aligned}
            \sigma_{R_{1N}}^{-1} R_{\text{CLB}, 1N}
            &\rightsquigarrow \mathcal{N}(0,1),\\
            \widehat\sigma_{R_{1N}}/\sigma_{R_{1N}}
            &\overset{P_0}{\to}1.
        \end{aligned}
    \]
    Consequently,
    $\widehat\sigma_{R_{1N}}^{-1}R_{\text{CLB},1N}
    \rightsquigarrow\mathcal N(0,1)$.
\end{theorem}

Let
\[
    {\UB}(R_{1N};\alpha)
    :=z_{1-\alpha}\sigma_{R_{1N}},
    \qquad
    \widehat{\UB}(R_{1N};\alpha)
    :=z_{1-\alpha}\widehat\sigma_{R_{1N}}.
\]
The two scales are defined above.  Theorem
\ref{thm:CLT_for_R1_CLB} implies that
$\widehat{\eta}_{\text{NSAVE}}-\widehat{\UB}(R_{1N};\alpha)$
provides a readily applicable conservative lower bound for
$\eta(\pi^*)$.

\begin{corollary}\label{cor:lower_bound_method}
    Under the conditions in Theorem \ref{thm:CLT_for_R1_CLB}, we have that 
    \[
        \liminf_{N \rightarrow \infty} \pr_{P_0} \Big( \eta(\pi^*) \geq \widehat{\eta}_{\text{NSAVE}} - \widehat{\UB}(R_{1N}; \alpha) \Big) \geq 1 - \alpha.
    \]
\end{corollary}

Here, we only establish asymptotic normality for the first term in the coarse decomposition. The reason is that ensuring asymptotic normality for $R_{\text{TCI}, 1N}$ requires regularity conditions for the \textit{Estimated Optimal Policies}. In contrast, as can be seen from Theorem \ref{thm:CLT_for_R1_CLB} and Corollary \ref{cor:lower_bound_method}, we do not impose any conditions on the estimated policy sequence $\{ \widehat{\pi}_{\tau(j - 1)} \}_{j > \ell_N}$. Therefore, compared with the conditions in \cite{shi2022statistical}, which require regularity and so-called margin conditions for both the estimated and true optimal policies,
our novel estimator $\widehat{\eta}_{\text{NSAVE}}$ admits valid inference procedures without such restrictions.

\subsection{Two-Sided Confidence Interval}\label{sec:2_sided_CI}

Under stronger conditions, more accurate inference via a Two-Sided Confidence Interval for $\eta^*$ is possible. To achieve this, we first need to establish the asymptotic properties of the two terms $R_{\text{TCI}, 1N}$ and $R_{\text{TCI}, 2N}$.
Here, $R_{\text{TCI},1N}$ is the fluctuation of NSAVE around
the value of the explicitly defined final refitted policy
$\widehat\pi_{\mathrm{fin}}$.  The martingale CLT applies to its
online evaluation component, while the difference between the
historical online policy values and
$\eta(\widehat\pi_{\mathrm{fin}})$ is controlled by the uniform
sequential regret rate.  The term $R_{\text{TCI},2N}$ is the regret
of the final refitted policy.

\begin{assumption}[Bellman-Complete Policy Class]
\label{ass:bellman_complete}
    The class $\Pi$ contains every deterministic stationary Markov
    policy produced by the fixed tie-breaking greedy rule used below.
    The optimal action-value function is the Bellman-optimal solution
    \[
        Q_{P_0}^*(a,s)
        =r_0(s,a)+\gamma
        \E\!\left[V_{P_0}^*(S')\mid S=s,A=a\right],
        \qquad
        V_{P_0}^*(s)=\max_aQ_{P_0}^*(a,s).
    \]
    Consequently, every policy greedy with respect to
    $Q_{P_0}^*$ is optimal.
\end{assumption}

\begin{assumption}[Estimated-Policy Stability]\label{ass:est_pol}
    Abbreviate
    $\widehat Q_j^*:=\widehat Q_{\tau(j-1)}^*$ and
    $\widehat\pi_j:=\widehat\pi_{\tau(j-1)}^{(Q)}$.
    Let $\widehat\pi_j$ be greedy with respect to
    $\widehat Q_j^*$.
    For some $\kappa_*>0$, and with $\alpha$ denoting the margin
    exponent in Assumption \ref{ass:margin_con},
    \[
        \begin{aligned}
        \max_{\ell_N<j\leq N}
        (j-\ell_N)^{\kappa_*}
        \|\widehat Q_j^*-Q_{P_0}^*\|_\infty=O_{P_0}(1), \quad
        n_N^{\kappa_*}
        \|\widehat Q_{\mathrm{fin}}^*-Q_{P_0}^*\|_\infty
        =O_{P_0}(1),
        \quad \kappa_*(1+\alpha)>1/2.
        \end{aligned}
    \]
\end{assumption}

This is a uniform sequential rate, rather than a separate pointwise
$O_{P_0}$ statement for every $j$; the uniform form is what justifies
summing the policy regrets over the online steps.  The separate final
rate applies to the all-sample refit used in the TCI decomposition.
These rates are not deduced
from a reverse occupancy-ratio inequality.  The sup-norm rate is
strong and is most plausible in tabular or correctly specified
finite-dimensional models, or under nonparametric conditions that
deliver uniform rather than merely integrated error control.  It is
needed only for the two-sided optimal-value result; the conservative
lower bound does not require it.  No convergence of the auxiliary
transition EIF is imposed here: the trajectory CLT studentizes each
episode-score increment by its own predictable conditional standard
deviation.%

Define the optimal action gap and the suboptimal-action set by
\[
    \begin{aligned}
    V_{P_0}^*(s)
    &:=\max_{a\in\mathcal A}Q_{P_0}^*(a,s),\quad
    \Delta^*(a,s)
    &:=V_{P_0}^*(s)-Q_{P_0}^*(a,s),\quad
    \mathcal A_{\mathrm{sub-opt}}(s)
    &:=\{a\in\mathcal A:\Delta^*(a,s)>0\}.
    \end{aligned}
\]
We use the convention that the minimum over an empty set is $+\infty$.

\begin{assumption}[Margin-Type Condition]\label{ass:margin_con}
    There exist $\alpha>0$ and $C<\infty$ such that, for
    all sufficiently small $\delta>0$,
    \[
        \sup_{\pi\in\Pi}
        \pr_{S\sim d_{P_0}^{\pi}}
        \left\{0<
        \min_{a\in\mathcal A_{\mathrm{sub-opt}}(S)}
        \Delta^*(a,S)
        \leq\delta\right\}
        \leq C\delta^\alpha.
    \]
\end{assumption}

\begin{theorem}\label{thm:CLT_for_R1_TCI}
    Under the conditions in Theorem \ref{thm:CLT_for_R1_CLB} and
    Assumptions \ref{ass:bellman_complete}--\ref{ass:margin_con}, we have
    \[
        \begin{aligned}
            \sigma_{R_{1N}}^{-1} R_{\text{TCI}, 1N}
            \rightsquigarrow \mathcal{N}(0,1),
            \qquad
            \widehat\sigma_{R_{1N}}^{-1} R_{\text{TCI}, 1N}
            \rightsquigarrow \mathcal{N}(0,1)
        \end{aligned}
    \]
    as $N \to \infty$. 
\end{theorem}

\begin{theorem}\label{thm:o_P_for_R2_TCI}
    Under the conditions in Theorem \ref{thm:CLT_for_R1_CLB} and
    Assumptions \ref{ass:bellman_complete}--\ref{ass:margin_con}, we have
    $
        R_{\text{TCI}, 2N} = o_{P_0}\big( (N- \ell_N)^{-1 / 2} \big)
    $
    under the rate condition in Assumption \ref{ass:est_pol}.
\end{theorem}

\begin{corollary}\label{cor:CLT_for_NSAVE}
    Under the conditions in Theorem \ref{thm:o_P_for_R2_TCI}, we have 
    \[
        \widehat\sigma_{R_{1N}}^{-1}
        \big(\widehat{\eta}_{\text{NSAVE}}-\eta(\pi^*)\big)
        \rightsquigarrow\mathcal N(0,1).
    \]
\end{corollary}

As partially shown in Corollary \ref{cor:CLT_for_NSAVE}, compared with the results in \cite{shi2022statistical}, our estimator $\widehat{\eta}_{\text{NSAVE}}$ demonstrates several advantages. Specifically:
\begin{itemize}
    \item Double-Robustness: Under policy alignment, our estimator
    retains the usual fixed-policy double-robustness identity.  We
    state this property formally below.
    \item Weaker restrictions on the estimated optimal policies: The convergence rate requirement for the estimated optimal policies is the same as that in \cite{shi2022statistical}, yet we do not require the associated effective sample size to be larger than a specific number inversely proportional to the convergence rate.
    \item Weaker constraints for the margin conditions: We only require that the probability, rather than the Lebesgue measure, satisfies the margin condition.
\end{itemize}

\subsection{Double-Robustness under Policy Alignment}

We next state the double-robustness property.  Semiparametric
efficiency is a separate question: for $N$ independent trajectories,
one must first derive the canonical gradient for the entire trajectory
experiment.  The martingale CLT above identifies the limiting variance
of the proposed score, but by itself does not prove that this variance
equals the trajectory-level efficiency bound.

\begin{theorem}\label{thm:DR_under_policy_alignment}
    Suppose Assumptions \ref{ass:data_obs}, \ref{ass:Mark},
    \ref{ass:stationary}, \ref{ass:traj_score_reg}, and
    \ref{ass:non_zero_var} hold.  Suppose the normalized trajectory
    scores use nuisances aligned with the same evaluation policy and
    \[
        \frac1{n_N}\sum_{j=\ell_N+1}^N
        \|\widehat Q_j-Q^{\widehat\pi_j}\|_{P_0,2}
        \|\widehat\omega_j-\omega^{\widehat\pi_j}\|_{P_0,2}
        =o_{P_0}(1).
    \]
    If, in addition,
    $\overline\eta_w(\widehat\pi)\overset{P_0}{\to}\eta^*$, then
    $\widehat\eta_{\mathrm{NSAVE}}\overset{P_0}{\to}\eta^*$.

    A sufficient double-robustness condition for the displayed
    product average is that one nuisance error converges uniformly to
    zero over $j>\ell_N$, while the other is uniformly
    $L_2(P_0)$-bounded.  The roles of $Q$ and $\omega$ may be
    interchanged.
\end{theorem}

Theorem \ref{thm:DR_under_policy_alignment} concerns consistency.  The
remaining question --- whether the limiting variance equals the
efficiency bound of the i.i.d.-episode experiment --- is settled in the
next subsection.  The answer has two parts: the \emph{discounted}
episode score underlying NSAVE is in general \emph{not} efficient,
because it weights the $T+1$ within-episode transitions unequally even
though Assumption~\ref{ass:stationary} makes them identically
distributed; but a \emph{uniformly} weighted variant is efficient, and
under a unique optimum it attains the trajectory-level bound.%

\subsection{Trajectory-Level Efficiency under Uniqueness}
\label{sec:traj_eff}

Throughout this subsection we work in the \emph{enlarged i.i.d.-episode
model} $\mathcal M^{\mathrm{ep}}$, in which the initial law $\mu_0$, the
behavior policy $b$, and the transition--reward law $g(r,s'\mid s,a)$
are variation-independent nonparametric components, evaluated at a true
law satisfying Assumption~\ref{ass:stationary}.  That assumption enters
only through the resulting fact that every transition tuple
$(S_t,A_t,R_t,S_{t+1})$ shares the stationary law
$P_{0,b}^{\mathrm{tr}}$, so the $T+1$ within-episode transitions are
identically distributed.  We do \emph{not} impose the stationary
fixed-point restriction $\mu_0=f_b(b,g)$ as a model constraint when
forming the tangent space; the effect of doing so is recorded in
Remark~\ref{rem:constrained_bound}.  Write
$\varepsilon^\pi:=R+\gamma V^\pi(S')-Q^\pi(A,S)$ for the Bellman
residual, so $\E(\varepsilon^\pi\mid S,A)=0$.

Define the \emph{uniformly weighted episode score}
\begin{equation}\label{eq:eff_score}
    \Phi^{\mathrm{eff}}_T(O_{0:T};Q,\omega,V,\pi)
    := \frac{1}{(1-\gamma)(T+1)}\sum_{t=0}^{T}
       \omega(A_t,S_t;\pi)\big[R_t+\gamma V(S_{t+1};\pi)
       -Q(A_t,S_t;\pi)\big]+V(S_0;\pi),
\end{equation}
and let
$\phi^{\mathrm{eff}}_T(O_{0:T};\pi)
:=\Phi^{\mathrm{eff}}_T(O_{0:T};Q^\pi,\omega^\pi,V^\pi,\pi)-\eta(\pi)$
be its centered population version.  It differs from the discounted
score $\phi_T$ of \eqref{eq:centered_episode_score} only in replacing
the within-episode weights $\gamma^t/(1-\gamma^{T+1})$ by the uniform
weights $1/\{(1-\gamma)(T+1)\}$; both average to $\eta(\pi)$, and by the
stationarity step in the proof of Lemma~\ref{lem:exact_dr_identity}
both carry the \emph{same} product-form double-robustness remainder, so
the weighting affects only the variance.  Concretely, if $c$ denotes
the expectation of the unscaled transition remainder in
Lemma~\ref{lem:exact_dr_identity}, stationarity makes every
within-episode transition have expectation $c$, and the uniform weights
sum to $(T+1)/\{(1-\gamma)(T+1)\}=1/(1-\gamma)$, the same factor
obtained from the normalized discounted weights.

\begin{theorem}[Episode efficient influence function and bound]
\label{thm:traj_eif}
    Fix $\pi\in\Pi$.  Under Assumptions~\ref{ass:data_obs},
    \ref{ass:Mark}, \ref{ass:stationary}, \ref{ass:regularity}, and
    \ref{ass:traj_score_reg}, $\eta(\pi)$ is pathwise differentiable in
    the enlarged i.i.d.-episode model $\mathcal M^{\mathrm{ep}}$, and its
    efficient influence function is $\phi^{\mathrm{eff}}_T(\cdot;\pi)$.
    Consequently the semiparametric efficiency bound (in
    $\mathcal M^{\mathrm{ep}}$) is
    \[
        \sigma^2_{\mathrm{eff},T}(\pi)
        =\var\{V^\pi(S_0)\}
        +\frac{\var_{P_{0,b}^{\mathrm{tr}}}
        \{\omega^\pi(A,S)\,\varepsilon^\pi\}}
        {(1-\gamma)^2\,(T+1)}.
    \]
\end{theorem}

\begin{corollary}[An efficient estimator for a fixed policy]
\label{cor:eff_est}
    Let
    $\widehat\eta_{\mathrm{eff}}(\pi)
    =\pn\,\Phi^{\mathrm{eff}}_T
    (O_{0:T};\widehat Q,\widehat\omega,\widehat V,\pi)$
    use cross-fitted, policy-aligned nuisances with
    $\|\widehat Q-Q^\pi\|_{P_0,2}
    +\|\widehat\omega-\omega^\pi\|_{P_0,2}=o_{P_0}(1)$ and
    $\|\widehat Q-Q^\pi\|_{P_0,2}\,
    \|\widehat\omega-\omega^\pi\|_{P_0,2}=o_{P_0}(N^{-1/2})$.  Then
    $\sqrt N\{\widehat\eta_{\mathrm{eff}}(\pi)-\eta(\pi)\}
    \rightsquigarrow\mathcal N(0,\sigma^2_{\mathrm{eff},T}(\pi))$, so
    $\widehat\eta_{\mathrm{eff}}(\pi)$ is semiparametrically efficient
    for $\eta(\pi)$.
\end{corollary}

\begin{corollary}[Efficiency for the optimal value under uniqueness]
\label{cor:eff_opt}
    Work in the enlarged model $\mathcal M^{\mathrm{ep}}$ and suppose the
    uniqueness Assumption~\ref{ass:deterministic_policy} holds with
    optimal policy $\pi^*_0$, together with
    Assumptions~\ref{ass:bellman_complete}--\ref{ass:margin_con}.
    Assume \emph{(i) sample splitting:} the final greedy policy
    $\widehat\pi_{\mathrm{fin}}$ and the policy-aligned nuisances
    $(\widehat Q,\widehat\omega,\widehat V)$ used to evaluate a given
    fold are fit on the complementary fold, hence are independent of the
    evaluation fold (equivalently, $K$-fold cross-fitting with fixed
    $K$, fold sizes proportional to $N$);
    \emph{(ii) a product nuisance rate for the random policy:} writing
    $\pi=\widehat\pi_{\mathrm{fin}}$, conditionally on the training
    fold,
    \[
        \|\widehat Q-Q^{\pi}\|_{P_0,2}
        +\|\widehat\omega-\omega^{\pi}\|_{P_0,2}=o_{P_0}(1),
        \qquad
        \|\widehat Q-Q^{\pi}\|_{P_0,2}\,
        \|\widehat\omega-\omega^{\pi}\|_{P_0,2}=o_{P_0}(N^{-1/2});
    \]
    and \emph{(iii) a $\sqrt N$-negligible regret:} the final-fit rate
    of Assumption~\ref{ass:est_pol} holds, under which
    Lemma~\ref{lem:greedy_regret_margin} gives
    $\sqrt N\{\eta(\pi^*_0)-\eta(\widehat\pi_{\mathrm{fin}})\}
    =o_{P_0}(1)$.  Then the efficiency bound for the optimal value
    $\Psi^*$ in $\mathcal M^{\mathrm{ep}}$ equals
    $\sigma^2_{\mathrm{eff},T}(\pi^*_0)$, and the cross-fitted plug-in at
    the final greedy fit attains it:
    \[
        \sqrt N\{\widehat\eta_{\mathrm{eff}}
        (\widehat\pi_{\mathrm{fin}})-\Psi^*(P_0)\}
        \rightsquigarrow
        \mathcal N\big(0,\sigma^2_{\mathrm{eff},T}(\pi^*_0)\big).
    \]
\end{corollary}

The discounted NSAVE increment has per-episode transition variance
$\var_{P_{0,b}^{\mathrm{tr}}}\{\omega^\pi\varepsilon^\pi\}\,
(1+\gamma^{T+1})/\{(1-\gamma^2)(1-\gamma^{T+1})\}$, whereas the efficient
transition variance is
$\var_{P_{0,b}^{\mathrm{tr}}}\{\omega^\pi\varepsilon^\pi\}/
\{(1-\gamma)^2(T+1)\}$.  The two agree at $T=0$, and the discounted one
is strictly larger for every $T\ge1$: as $T\to\infty$ the efficient term
vanishes at rate $(T+1)^{-1}$, while the discounted term stays bounded
below by
$\var_{P_{0,b}^{\mathrm{tr}}}\{\omega^\pi\varepsilon^\pi\}/(1-\gamma^2)$.
Thus the discounted weighting of NSAVE buys its predictable martingale
structure --- hence the sequential, assumption-light lower bound of
Section~\ref{sec:lower_bound} and its validity under non-uniqueness ---
at a non-vanishing transition-variance cost.  In relative terms, the
gap between the discounted and efficient transition variances grows
with the horizon because the efficient transition component is averaged
over the $T+1$ stationary decision points.  When a unique optimum is
the target and efficiency is the priority, the uniformly weighted estimator
$\widehat\eta_{\mathrm{eff}}$ of \eqref{eq:eff_score} should be used.

\begin{remark}[Constrained stationary model]
\label{rem:constrained_bound}
    Imposing the stationarity relation $\mu_0=f_b(b,g)$ as a model
    restriction can only shrink the tangent space, and hence can only
    lower the efficiency bound; thus $\sigma^2_{\mathrm{eff},T}(\pi)$
    remains an attainable regular asymptotic variance under that
    restriction.  Characterizing the \emph{sharp} constrained bound
    requires projecting the gradient onto the smaller tangent space,
    which we leave to future work (see Section~\ref{sec:final}); we state
    the efficiency results in the enlarged model
    $\mathcal M^{\mathrm{ep}}$, where $\sigma^2_{\mathrm{eff},T}(\pi)$ is
    sharp.
\end{remark}

\section{Alternative Inference Approaches}\label{sec:alternative}
\subsection{Smoothing}\label{sec:smoothing}

Softmax policies provide a useful approximation to greedy policies.
However, applying softmax to an estimated \emph{hard} optimal
$Q$-function does not by itself make the map
$P\mapsto Q_P^*$ differentiable at ties.  We therefore treat the
construction below as a sample-split, pointwise approximation method,
not as a general restoration of pathwise differentiability.  A
genuinely smooth regularized target would instead require a soft
Bellman operator and its corresponding influence function.

For any real-valued \textit{function} $h: \mathbb{R}^p \rightarrow \mathbb{R}$ and vector $\vv \in \mathbb{R}^d$, define the softmax smoothing approximation $\varphi_{\beta} \{ \cdot \}$ and the multiple softmax operator $\sm_{\beta} \{ \cdot \}$ such that 
\begin{equation}\label{eq:smoothing_operators}
    \varphi_{\beta} \{ h(\cdot)\}: = h(\cdot) \times \frac{\exp\{\beta h(\cdot)\}}{1 + \exp\{\beta h(\cdot)\}} \quad \text{and} \quad \sm_{\beta} \{ \vv \} := \frac{\sum_{j = 1}^d v_j \exp\{\beta v_j\}}{\sum_{j = 1}^d \exp\{\beta v_j\}},
\end{equation}
where $\beta > 0$ denotes the degree of smoothing.

In the static case, the $Q$-function under a policy $\pi$ reduces to $Q (a, s ;\pi) = Q(a, s)$, as $\pi$ simply selects an action $a$ given $x$ to maximize the $Q$-function.
As pointed out in \cite{whitehouse2025inference}, by smoothing the \textit{value} function $\eta^*(P) = \E_P [\max_{a} Q(a, S)]$ with 
$
    \eta_{\beta}^* = \E_P \big[\sm_{\beta} \big\{ Q(\va, S) \big\} \big],
$
one can differentiate $\eta_{\beta}^*(P_{\epsilon})$ with respect to $\epsilon$ and then use the first-order condition to construct a Neyman orthogonal score (or the estimating equation for the pathwise derivative at $\epsilon = 0$).
However, challenges arise when extending their framework to the dynamic setting inherent in MDPs: for any fixed $\pi$, the (dynamic) $Q$-function is derived from the fixed point of the Bellman equation, such that 
$
    V_P^\pi=r_{P,\pi}+\gamma K_{\pi,P}V_P^\pi,
$
where $K_{\pi,P}$ is the transition kernel defined in
\eqref{eq:transition_kernel}.
A static soft maximum applied after estimating $Q_P^*$ does not, by
itself, identify the value of a fixed policy or supply the
policy-evaluation nuisances needed by the dynamic one-step score.
We therefore use the policy optimization perspective found in
entropy-regularized MDPs \citep{neu2017unified}: instead of treating a
smoothed statewise maximum as the target value, we smooth the greedy
policy and then evaluate that policy through its own Bellman equation.
Specifically, define the \textit{smoothing-greedy} policy
$\pi_{\beta}(P)$ as
\begin{equation}\label{eq:smoothing_greedy_pol}
    \pi_{\beta}(P)(a \mid s) := \frac{\exp\{\beta Q(P)(a, s; \pi^*(P))\}}{\sum_{a' \in \mathcal{A}} \exp\{\beta Q(P)(a', s; \pi^*(P))\}}.
\end{equation}
As $\beta\to\infty$, $\pi_\beta(P)$ converges to the uniform mixture
over maximizing actions.  This limit is an optimal policy, but it need
not equal an arbitrarily selected representative $\pi^*(P)$ when ties
are present.

Our procedure for smoothed nuisance estimation proceeds sequentially:
\begin{itemize}
    \item First, we estimate the \textit{optimal} Q-function $Q^*$ using \textit{any} off-policy algorithm (e.g., Fitted Q-Iteration), yielding $\widehat{Q}_{\text{opt}}(\cdot, \cdot)$;
    \item Second, using this estimate and under a chosen smoothing sequence $\beta_N$, we construct the plug-in policy $\widehat{\pi}_{\beta_N}$ using $\widehat{Q}$ as
    $
        \widehat{\pi}_{\beta_N}(a \mid s) := \frac{\exp\{\beta_N \widehat{Q}_{\text{opt}}(a, s)\}}{\sum_{a' \in \mathcal{A}} \exp\{\beta_N \widehat{Q}_{\text{opt}}(a', s)\}};
    $
    \item Finally, on the training sample we separately estimate the
    policy-evaluation nuisances
    $Q^{\widehat\pi_{\beta_N}}$,
    $V^{\widehat\pi_{\beta_N}}$, and
    $\omega^{\widehat\pi_{\beta_N}}$.  These are denoted by
    $\widehat Q_{\mathrm{eval}}$,
    $\widehat V_{\mathrm{eval}}$, and
    $\widehat\omega_{\mathrm{eval}}$.
\end{itemize}
These nuisance estimates are then plugged into the one-step estimator, leading to
\begin{equation}\label{eq:smoothed_est}
    \widehat{\eta}_{\beta_N}
    :=
    \pn\EstFunTraj_{\eta(\pi)}
    (O_{0:T};
    \widehat Q_{\mathrm{eval}},
    \widehat\omega_{\mathrm{eval}},
    \widehat V_{\mathrm{eval}},
    \widehat\pi_{\beta_N}).
\end{equation}
Using $\widehat Q_{\mathrm{opt}}$ itself in the evaluation score would
generally be incorrect, because
$Q^*\neq Q^{\pi_\beta}$ for finite $\beta$.
The smoothed estimator \eqref{eq:smoothed_est} can be regarded as a modified version of our NSAVE estimator, where we replace the sequential evaluation with the smoothing technique. Specifically, we decompose the difference between $\widehat{\eta}_{\beta_N}$ and the true value $\eta(\pi^*)$ as follows:
\[
    \begin{aligned}
        \widehat{\eta}_{\beta_N} - \eta(\pi^*) = \underbrace{\widehat{\eta}_{\beta_N} - \eta(\widehat\pi_{\beta_N})}_{:= R_{\text{SM}, 1N}} + \underbrace{\eta(\widehat\pi_{\beta_N}) - \eta(\pi^*)}_{:= R_{\text{SM}, 2N}}.
    \end{aligned}
\]
As shown in the proof of Theorem \ref{thm:CLT_for_R1_TCI}, provided consistent nuisance estimates with appropriate convergence rates are selected, the statistical error $R_{\text{SM}, 1N}$ will converge to a normal distribution. Meanwhile, the policy-value error $R_{\text{SM}, 2N}$, controlled by the smoothing parameter $\beta_N$, becomes $o_{P_0}(N^{-1 / 2})$ via the smoothing mechanism rather than sequential evaluation.
The following theorem gives a transparent sufficient condition for a
pointwise Gaussian limit.  Its deliberately strong first-stage rate
condition isolates policy-learning error from policy-evaluation error.

\begin{theorem}\label{thm:smoothing_RAL}
    Suppose Assumptions \ref{ass:data_obs}, \ref{ass:Mark},
    \ref{ass:stationary}, \ref{ass:traj_score_reg}, and
    \ref{ass:bellman_complete}, and \ref{ass:margin_con} hold.
    Suppose the policy and nuisance
    estimators are trained on a sample independent of the $N$ i.i.d.
    evaluation episodes, and
    \[
        r_N:=\|\widehat Q_{\mathrm{opt}}-Q_{P_0}^*\|_\infty
        \quad\text{satisfies}\quad
        \beta_N r_N=o_{P_0}(N^{-1/2}).
    \]
    Assume also
    \[
        \beta_N\to\infty,\qquad
        \beta_N^{-(1+\alpha)}=o(N^{-1/2}),
    \]
    and the policy-evaluation nuisances satisfy
    \[
        \|\widehat Q_{\mathrm{eval}}
        -Q^{\widehat\pi_{\beta_N}}\|_{P_0,2}
        +
        \|\widehat\omega_{\mathrm{eval}}
        -\omega^{\widehat\pi_{\beta_N}}\|_{P_0,2}
        =o_{P_0}(1),
    \]
    together with
    \[
        \|\widehat Q_{\mathrm{eval}}
        -Q^{\widehat\pi_{\beta_N}}\|_{P_0,2}
        \|\widehat\omega_{\mathrm{eval}}
        -\omega^{\widehat\pi_{\beta_N}}\|_{P_0,2}
        =o_{P_0}(N^{-1/2}).
    \]
    Suppose that the population and estimated softmax policies and
    their limit belong to the evaluated policy class $\Pi$, and let
    $\pi_\infty$ be the uniform mixture over optimal actions.  Then
    \[
        \sqrt N\{\widehat\eta_{\beta_N}-\eta^*\}
        =
        \frac1{\sqrt N}\sum_{i=1}^N
        \phi_T(O_i;\pi_\infty)+o_{P_0}(1),
    \]
    where $\phi_T(\cdot;\pi_\infty)$ is the centered normalized
    trajectory score.  Hence the estimator has a pointwise Gaussian
    limit with variance $\var\{\phi_T(O;\pi_\infty)\}$.
\end{theorem}

The same conclusion holds for a fixed number of cross-fitting folds
if the displayed rates and boundedness conditions hold uniformly over
folds and every fold size is asymptotically proportional to $N$.
Indeed, the proof below applies conditionally within each fold, and
the foldwise influence-function sums recombine into the full-sample
sum.

The two smoothing-bias conditions are deliberately restrictive.  For
example, choosing $\beta_N$ just above
$N^{1/\{2(1+\alpha)\}}$ already requires the first-stage sup-norm
error to be smaller than
$N^{-1/2-1/\{2(1+\alpha)\}}$ up to slack.  Thus this theorem is a
transparent sufficient result, not a claim that generic
nonparametric learners automatically satisfy the required rates.
Computationally, the procedure estimates the optimal $Q$-function
once to construct the policy and then estimates the two
policy-evaluation nuisances for that policy.
Theorem \ref{thm:smoothing_RAL} is pointwise.  When optimal policies
have different fixed-policy gradients, it does not imply regularity
under local alternatives or semiparametric efficiency for the
nonregular optimal-value parameter.%

\subsection{Post-Selection Inference }\label{sec:PSI}

When $\mathcal{A} \times \mathcal{S}$ is finite and small, or more generally when the candidate policy class $\Pi = \{ \pi_1, \dots, \pi_K \}$ is finite, we can employ \textit{Post-Selection Inference} (PSI) techniques to address the non-regularity. Unlike the smoothing approach, which modifies the target parameter to a smooth approximation $\eta(\pi_{\beta})$, PSI aims to construct a valid confidence interval for the \textit{value of the empirically selected policy} itself, denoted as $\eta(\widehat{\pi}_N)$, where $\widehat{\pi}_N = \arg\max_{\pi \in \Pi} \widehat{\eta}(\pi)$.

Standard inference that treats $\widehat{\pi}_N$ as a fixed policy fails to account for the \textit{winner's curse}: the selection process systematically favors policies with positive estimation noise, leading to an upward bias in the naive estimator. 

To rigorously correct for this bias while accounting for the potential non-uniqueness of optimal policies (ties) and the high correlation between OPE estimates, we adopt the \textbf{Two-Step Inference on Multiple Winners} framework proposed by \cite{petrou2024inference}.
Our PSI procedure proceeds sequentially as follows:
\begin{itemize}
    \item \textbf{Step 0: OPE estimation and selection.} We estimate the values for all candidate policies.
    Let $\widehat{\bm{\eta}}=(\widehat{\eta}(\pi_1),\dots,\widehat{\eta}(\pi_K))^\top$ be the OPE estimates (e.g., doubly robust for each $\pi_k$). We also estimate the asymptotic covariance matrix $\widehat{\boldsymbol{\Sigma}}$, for example from the empirical covariance of their jointly asymptotically linear episode scores, such that
    $
        \sqrt{N}(\widehat{\bm{\eta}} - \bm{\eta}) \rightsquigarrow \mathcal{N}(\bm{0}, \boldsymbol{\Sigma}) \quad \text{with} \quad \widehat{\boldsymbol{\Sigma}} \overset{\pr}{\to} \boldsymbol{\Sigma} \quad \text{and} \quad \widehat{\mathcal{A}}_{\text{opt}} = \arg\max_{k \in [K]} \widehat{\bm{\eta}}$;

    \item \textbf{Step 1: A $(1-\delta_1)$ confidence region for the nuisance governing selection.} Construct a confidence region $\mathcal{C}_\eta(\widehat{\bm{\eta}};\delta_1)\subseteq \mathbb{R}^K$ such that
        $\liminf_{N\to\infty} \pr\Big(\bm{\eta}\in \mathcal{C}_\eta(\widehat{\bm{\eta}};\delta_1)\Big)\ \ge\ 1-\delta_1$.
    We take the confidence region to contain its center
    $\widehat{\bm\eta}$, as is the case for the usual centered
    rectangular or ellipsoidal regions.  Using
    $\mathcal{C}_\eta$, define the \emph{plausible optimal set}
        $\widehat{\mathcal{A}}^{+} := \bigcup_{\bm{\eta} \in \mathcal{C}_\eta(\widehat{\bm{\eta}};\delta_1)} \arg\max_{k\in[K]} \penalty 0 \eta_k$.
    By construction,
    $\widehat{\mathcal A}_{\mathrm{opt}}\subseteq
    \widehat{\mathcal A}^{+}$ because
    $\widehat{\bm\eta}\in\mathcal C_\eta$, and, on
    ${\bm{\eta}\in
    \mathcal{C}_\eta(\widehat{\bm{\eta}};\delta_1)}$, the true
    optimal set
    $\mathcal A_{\mathrm{opt}}(\bm\eta):=
      \arg\max_{k\in[K]}\eta_k$
    is also contained in $\widehat{\mathcal A}^{+}$.

    \item \textbf{Step 2: Calibrate a simultaneous critical value over $\widehat{\mathcal{A}}^{+}$.} Let the standardized errors be
    $Z_k := \sqrt N\,\widehat{\boldsymbol{\Sigma}}_{kk}^{-1 / 2}
    \big(\widehat{\eta}(\pi_k)-\eta(\pi_k)\big)$.
    With
    $\mathcal E_1
      :=\{\bm\eta\in
      \mathcal C_\eta(\widehat{\bm\eta};\delta_1)\}$,
    choose a possibly data-dependent critical value
    $q_N(\widehat{\mathcal A}^{+})$ satisfying the joint
    two-step error guarantee
    \[
        \limsup_{N\to\infty}
        \pr\left(
        \mathcal E_1\cap
        \left\{\max_{k\in\widehat{\mathcal A}^{+}}|Z_k|
        >q_N(\widehat{\mathcal A}^{+})\right\}\right)
        \leq\delta_2-\delta_1.
    \]
    Because $\widehat{\mathcal A}^{+}$ depends on the same data as
    $\bm Z$, this condition does not follow merely by inserting the
    realized set into an unconditional Gaussian max-quantile.  It
    must be verified by the chosen two-step, locally simultaneous, or
    projection construction.

    \item \textbf{Final PSI confidence set (reported on the selected set).} Define
        $\mathcal{C}_{\text{PSI}} := \bigtimes_{k\in \widehat{\mathcal{A}}_{\text{opt}}} \Bigl[ \widehat{\eta}(\pi_k)\ \pm\ q_N(\widehat{\mathcal{A}}^{+})\cdot\sqrt{\widehat{\boldsymbol{\Sigma}}_{kk}/N}\Bigr]$.
\end{itemize}

The conservative projection implementation always satisfies the
required condition: take $\delta_1=0$,
$\widehat{\mathcal A}^{+}=[K]$, and use the
$(1-\delta_2)$ Gaussian max-quantile over all $K$ coordinates.
Sharper two-step or locally simultaneous implementations may replace
projection when their own calibration theorem verifies the displayed
joint error bound; see Appendix~\ref{sec:psi_other_methods} and
\citet{petrou2024inference}.

\begin{corollary}\label{cor:PSI}
    Suppose
    $0\leq\delta_1\leq\delta_2\leq1$,
    $\liminf_N\Pr(\mathcal E_1)\geq1-\delta_1$, the confidence region
    contains $\widehat{\bm\eta}$, and the Step~2 joint error guarantee
    holds.  Then
    $
        \liminf_{N \to \infty} \pr_{P_0} \Big\{  \big(\eta(\pi_k)\big)_{k \in \widehat{\mathcal{A}}_{\text{opt}}} \in  \mathcal{C}_{\text{PSI}} \Big\} \geq 1 - \delta_2.
    $
\end{corollary}

If the optimal policy is uniquely separated,
$\Pr(\widehat{\mathcal A}^{+}=\{k^*\})\to1$, and the first-step error
allocation satisfies $\delta_{1,N}\to0$, then the resulting critical
value can converge to the ordinary two-sided Gaussian critical value,
so the interval has the same first-order length as the oracle interval.
Without these additional conditions, Corollary~\ref{cor:PSI} gives
coverage but does not by itself imply oracle-length efficiency.  With
ties, the simultaneous set remains valid for all reported winners
under the stated joint calibration condition.

\section{Simulations}\label{sec:sim}

We design the Monte Carlo experiments to isolate the claims made by the
theory.  Episodes are independent, each episode starts from the stationary
distribution of the behavior chain, rewards are Bernoulli, and
$\gamma=0.7$.  Thus the maintained stationarity, bounded-reward, and overlap
conditions hold by construction.  We consider
$N\in\{100,250,500,1000\}$ episodes and $T\in\{25,50\}$ observed
transitions, with 500 replications for every design point.  All trajectory
scores use the finite-window normalization $1-\gamma^T$ and the marginal
discounted occupancy ratio, rather than a cumulative product of action
ratios.  Appendix~\ref{sec:sim_setup} gives the complete data-generating
mechanisms and implementation details.

The experiments have four complementary purposes.  First, a fixed-policy
experiment varies whether $Q^\pi$ and $\omega^\pi$ are correct, misspecified,
or estimated, thereby testing the exact double-robust remainder.  Second, a
six-state MDP with a uniformly separated optimal action examines the regular
unique-optimum regime.  Third, a one-state MDP with exactly equal action
values creates first-order distinguishable tied policies.  Finally, near-tie
and weak-overlap designs stress the smoothing approximation and importance
weighting; these results are reported in Appendix~\ref{sec:sim_impl}.%

We report two NSAVE intervals: the one-sided conservative lower bound
(NSAVE--CLB) and the two-sided interval that additionally uses the
policy-regret condition (NSAVE--TCI).  We also report a 50--50 sample-split
smoothed estimator with
$\beta_N=\{\log(\lfloor N/2\rfloor)\}^2$, and a sample-split hard-greedy DR
benchmark.  For both split estimators, the evaluation score uses
$Q^{\widehat\pi}$ and $\omega^{\widehat\pi}$ for the policy actually learned
on the training half; in particular, the smoothed score does not insert a
hard-optimal $Q$-function.

\paragraph{Fixed-policy double robustness.}
Table~\ref{tab:mc_dr} reports representative results for $T=50$.  Bias is
small and coverage is close to nominal whenever either nuisance is correct.
The estimated-both estimator behaves similarly.  In contrast, deliberately
setting both $Q$ and $\omega$ incorrectly produces a persistent bias of about
$-0.61$ and zero coverage.  Across all $N$ and both horizons, coverage ranges
from 0.918 to 0.964 for the four valid configurations, whereas it is zero
throughout when both nuisances are wrong.  This pattern is the finite-sample
signature of the product remainder in the fixed-policy DR identity.

\begin{table}[htbp]
\centering
\caption{Fixed-policy double-robustness experiment for $T=50$; coverage is
for nominal 95\% two-sided intervals.  ``Wrong
$Q$'' sets $Q=V=0$, and ``wrong $\omega$'' sets $\omega=1$.  Each entry is
based on 500 replications.}
\label{tab:mc_dr}
\begin{tabular}{lrrrr}
\toprule
& \multicolumn{2}{c}{$N=100$} & \multicolumn{2}{c}{$N=1000$} \\
\cmidrule(lr){2-3}\cmidrule(lr){4-5}
Nuisance configuration & Bias & Coverage & Bias & Coverage \\
\midrule
Oracle $Q$, oracle $\omega$       & $-0.001$ & 0.942 & $-0.001$ & 0.944 \\
Oracle $Q$, wrong $\omega$        & $-0.004$ & 0.946 & $\phantom{-}0.001$ & 0.944 \\
Wrong $Q$, oracle $\omega$        & $-0.015$ & 0.918 & $-0.001$ & 0.964 \\
Estimated $Q$, estimated $\omega$ & $-0.008$ & 0.946 & $\phantom{-}0.002$ & 0.940 \\
Wrong $Q$, wrong $\omega$         & $-0.614$ & 0.000 & $-0.611$ & 0.000 \\
\bottomrule
\end{tabular}
\end{table}

\begin{figure}[!htbp]
\centering
\begin{minipage}[t]{0.62\linewidth}
\centering
\vspace{0pt}
\includegraphics[width=\linewidth]{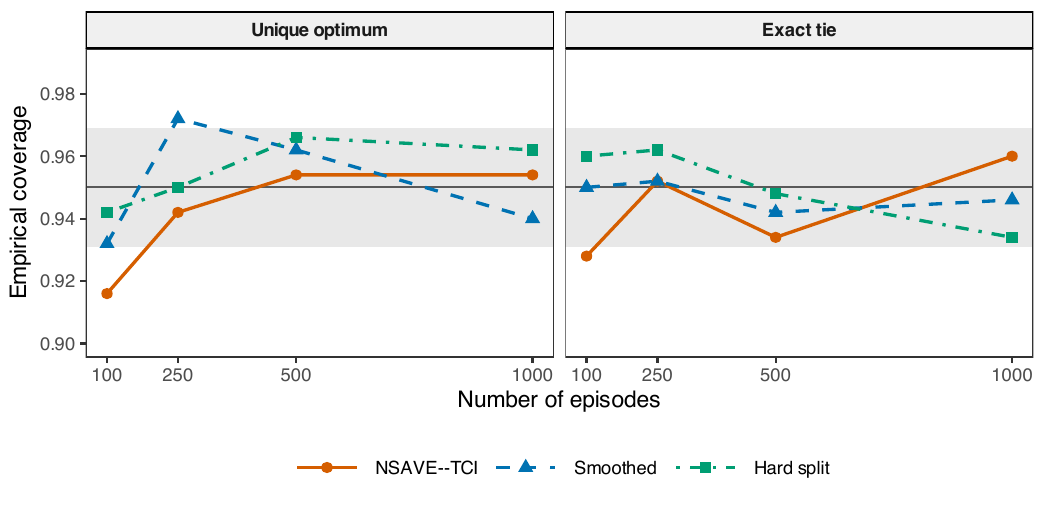}
\par\smallskip
{\small (a) Empirical coverage of nominal 95\% two-sided intervals.}
\par\medskip
\includegraphics[width=\linewidth]{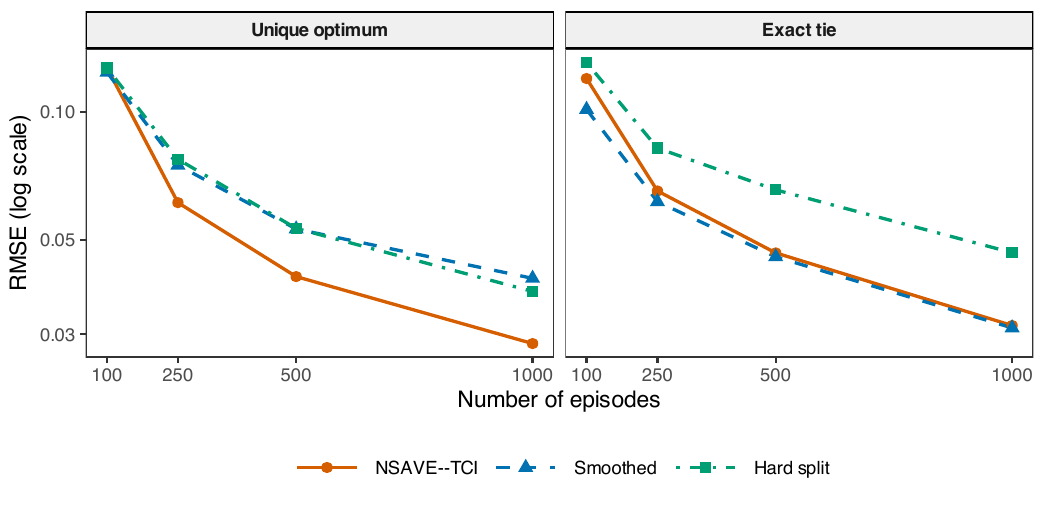}
\par\smallskip
{\small (b) Root mean squared error.}
\end{minipage}\hfill
\begin{minipage}[t]{0.34\linewidth}
\centering
\vspace{0pt}
\includegraphics[width=\linewidth]{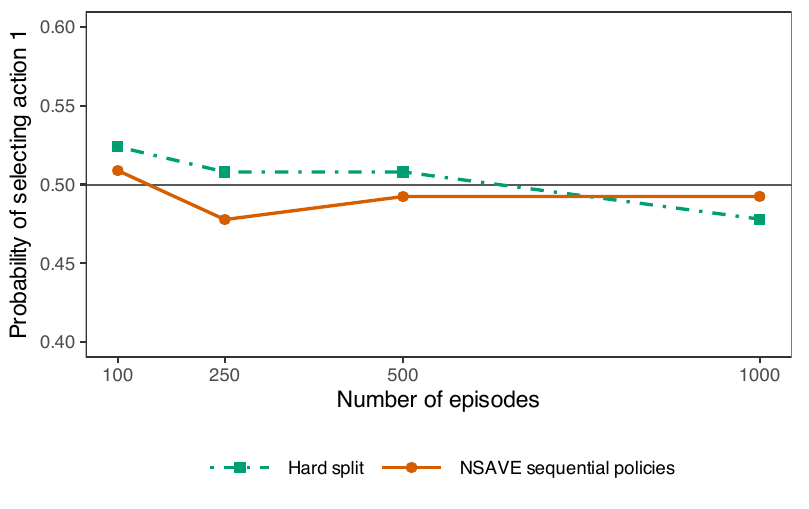}
\par\smallskip
{\small (c) Probability of selecting action~1 in the exact-tie design.}
\end{minipage}
\caption{Monte Carlo performance for $T=50$.  Panels (a)--(b) show
coverage and RMSE under a unique optimum and an exact tie; the gray band in
panel (a) is the 95\% Monte Carlo uncertainty band around 0.95 based on 500
replications.  Panel (c) shows that neither hard policy selector converges
to a fixed action in the exact-tie design, although all selected policies
have the same value.}
\label{fig:mc_main}
\label{fig:mc_tie_selection}
\end{figure}

\paragraph{Unique and tied optima.}
Panels (a)--(b) of Figure~\ref{fig:mc_main} report two-sided coverage and RMSE for $T=50$.
The shaded region in the coverage panel is the 95\% Monte Carlo uncertainty
band around 0.95, whose half-width is approximately 0.019 with 500
replications.  In the unique-optimum design, NSAVE--TCI has some
small-sample undercoverage at $N=100$ (0.916), but its coverage is
0.942--0.954 for $N\geq250$ and its RMSE decreases from 0.127 to 0.029.
The hard-split and smoothed procedures also show decreasing RMSE and
generally calibrated pointwise coverage.  The weighted historical-policy
regret of NSAVE decreases from $2.68\times10^{-3}$ at $N=100$ to
$2.10\times10^{-4}$ at $N=1000$; the smoothing regret decreases from 0.028
to 0.003.

Under the exact tie, the RMSE of all procedures again decreases at the
expected root-$N$ scale and empirical coverage lies between 0.928 and 0.962
over the displayed sample sizes.  At the same time,
panel (c) of Figure~\ref{fig:mc_tie_selection} shows that the hard learned policy does not
converge to a fixed action: both the sequential and split procedures choose
action~1 with probability close to one half even at $N=1000$.  Since the
tied policies have the same value, their policy regret is zero.  This
experiment is pointwise; it illustrates stable value inference despite
nonconvergent policy selection, but it does not assert regularity under
local alternatives.

The one-sided NSAVE--CLB coverage ranges from 0.918 to 0.964 over all
unique, tied, near-tied, and overlap designs.  Appendix~\ref{sec:sim_impl}
shows that weak overlap primarily widens intervals, while the finite
temperature used here can leave non-negligible smoothing bias in the
near-tie design.  We report this latter behavior explicitly because it
demonstrates the practical content of the smoothing-bias condition rather
than contradicting the pointwise smoothing result.

\section{Application to the Drink Less MRT}
\label{sec:realdata}

We apply the proposed methods to the Drink Less micro-randomized trial (MRT), a mobile health study of push notifications designed to increase engagement with an alcohol-reduction app \citep{bell2023notifications}.
At each daily decision point, participants were randomized to receive no notification, a new notification, or the standard notification.
The original MRT analysis showed that notifications can increase near-term engagement and motivated further optimization of the notification policy.
Our goal is complementary: using the randomized longitudinal data, we evaluate learned state-adaptive policies and construct confidence intervals for their discounted values and their improvements over the randomized behavior policy.

The analysis uses the shared Drink Less dataset with $N=349$ complete 30-day participant trajectories.
The action set is
$\mathcal A=\{\text{no notification},\text{new notification},\text{standard notification}\}$, with known behavior probabilities $(0.40,0.30,0.30)$.
The reward is the binary primary outcome indicating whether the participant opened the app between 8pm and 9pm.
We use days 15--30 after a 14-day burn-in period, so that each participant contributes a length-15 trajectory, and set $\gamma=0.7$.
The primary state is an 18-level history summary recording current pre-decision engagement, recent three-day engagement, and recent three-day notification burden.
The target initial law is set to the pooled behavior stationary state law, in accordance with Assumption~\ref{ass:stationary}.
Further preprocessing and sensitivity analyses are given in Appendix~\ref{sec:drinkless_setup}.

We report three learned-policy analyses.
The primary analysis is NSAVE with predictable sequential policy fitting: after an initial training sample, the policy used for a participant is learned only from prior participants, and the resulting normalized episode scores are aggregated with predictable inverse-scale weights.
For comparison, we also report five-fold cross-fitted hard-greedy and finite-temperature smoothed policies.
All improvement intervals are based on paired participant-level scores relative to the randomized behavior policy, rather than subtracting two separately estimated values.

Table~\ref{tab:drinkless_main} summarizes the main estimates.
The randomized behavior policy has estimated value $0.259$.
The NSAVE historical learned policies have estimated value $0.337$ and estimated improvement $0.095$ over behavior, with a 95\% two-sided confidence interval $[0.024,0.167]$ and a 95\% one-sided lower confidence bound $0.035$.
Because every historical learned-policy value is bounded above by the optimal policy value, the one-sided NSAVE lower bound also gives a conservative lower bound for the optimal improvement over behavior.
The cross-fitted hard-greedy and smoothed analyses yield qualitatively similar improvements, with the smoothed analysis giving the largest and most stable estimate.

\begin{table}[!htp]
\centering
\small
\caption{Drink Less main real-data results.  Values are discounted policy values or paired improvements over the randomized behavior policy.  Parentheses give 95\% two-sided confidence intervals; the final column gives the 95\% one-sided lower confidence bound for the improvement.}
\label{tab:drinkless_main}
\begin{tabular}{lccc}
\toprule
Method & Policy value & Improvement & One-sided lower \\
\midrule
Behavior & $0.259\ (0.212,0.307)$ & \NA & \NA \\
NSAVE & $0.337\ (0.251,0.423)$ & $0.095\ (0.024,0.167)$ & $0.035$ \\
Hard greedy & $0.334\ (0.244,0.424)$ & $0.075\ (0.007,0.144)$ & $0.018$ \\
Smoothed & $0.382\ (0.297,0.468)$ & $0.124\ (0.067,0.181)$ & $0.076$ \\
\bottomrule
\end{tabular}
\end{table}

\begin{figure}[!htbp]
    \centering
    \begin{minipage}[t]{0.48\linewidth}
    \centering
    \vspace{0pt}
    \includegraphics[width=\linewidth]{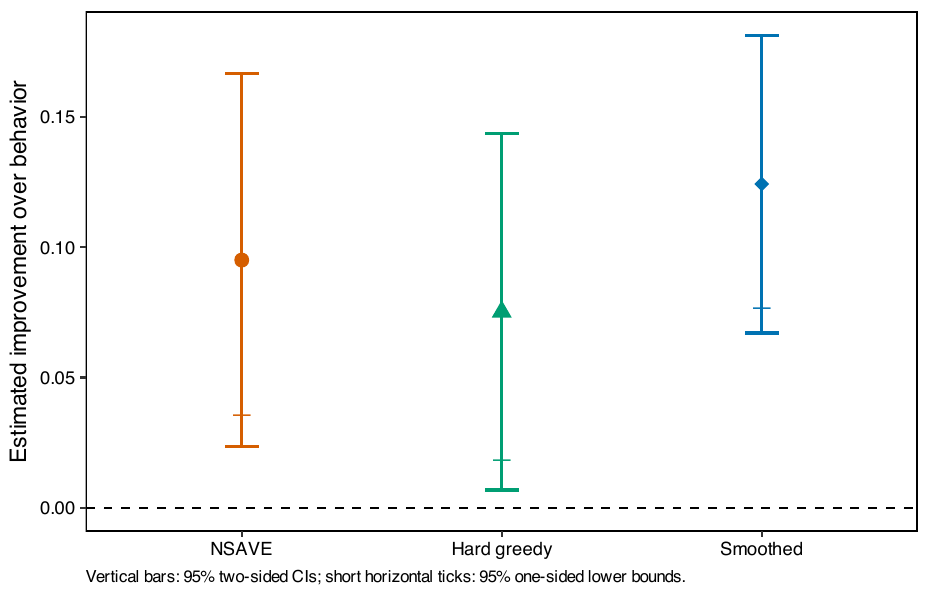}
    \par\smallskip
    {\small (a) Estimated improvement over the randomized behavior policy.}
    \end{minipage}\hfill
    \begin{minipage}[t]{0.48\linewidth}
    \centering
    \vspace{0pt}
    \includegraphics[width=\linewidth]{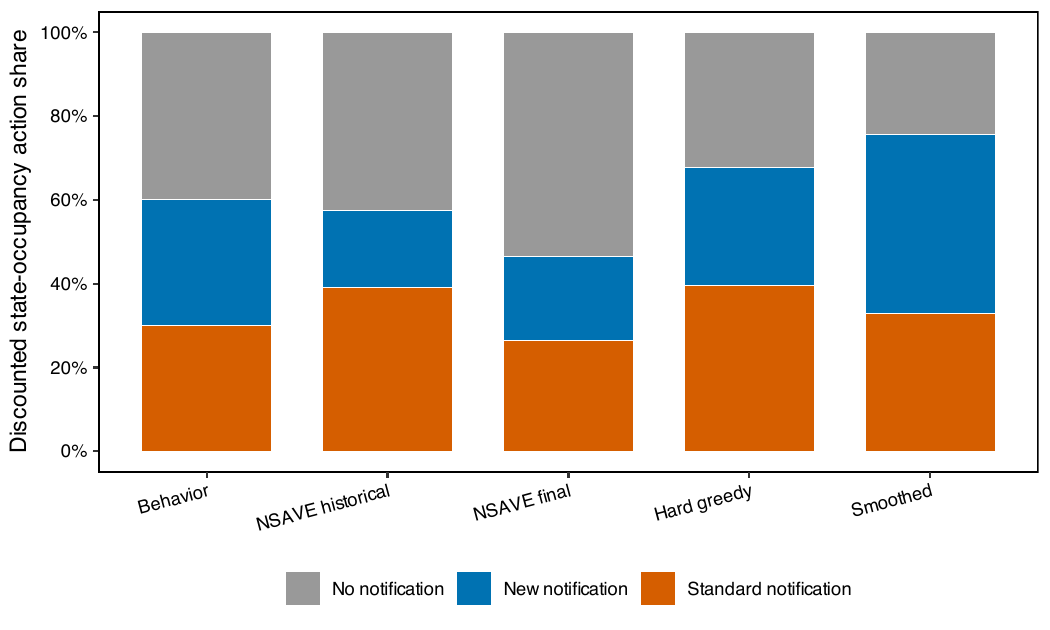}
    \par\smallskip
    {\small (b) Discounted state-occupancy action composition.}
    \end{minipage}
    \caption{Drink Less real-data diagnostics.  Panel (a) shows estimated
    improvement over the randomized behavior policy; vertical bars are 95\%
    two-sided confidence intervals, and short horizontal ticks are 95\%
    one-sided lower confidence bounds.  Panel (b) shows the discounted
    state-occupancy action composition for the randomized behavior policy
    and learned notification policies.  The action composition is
    descriptive; near-ties in fitted action values mean that it should not
    be interpreted as identifying a unique optimal rule.}
    \label{fig:drinkless_improvement}
    \label{fig:drinkless_policy}
\end{figure}

Panel (a) of Figure~\ref{fig:drinkless_improvement} visualizes the paired improvement intervals.
All three analyses produce positive point estimates and positive one-sided lower bounds.
The result should be interpreted as evidence for improved policy value, not as evidence that a unique deterministic notification rule has been identified.
Indeed, the fitted top-two action-value gaps are small for a substantial part of the discounted state occupancy, and different regularized learning procedures select somewhat different action mixtures.
Panel (b) of Figure~\ref{fig:drinkless_policy} reports these action-composition diagnostics.
The empirical picture is therefore consistent with the non-regular setting emphasized by the theory: several near-tied policies can have similar values, while value inference remains stable when formulated through policy-aligned episode scores.

\section{Final Remarks}\label{sec:final}

Finally, we emphasize that, in the auxiliary augmented
transition-sampling experiment, existence of an efficient influence
function for the optimal value hinges critically on regularity of the
policy optimization map.  When the optimal policy is unique, the
problem reduces locally to fixed-policy inference in that experiment
\citep{uehara2022review,shi2025statistical}.  The trajectory
procedures in this paper are instead justified by episode-score
martingale arguments.
When optimal policies are non-unique and have distinct fixed-policy
gradients, the value functional is nonregular and standard root-$N$
regular inference can fail, a phenomenon closely related to
nonregular parameters in optimal treatment regimes and post-selection
inference \citep{laber2014dynamic,whitehouse2025inference}.

Our NSAVE procedure provides stable martingale-based inference under
the stated policy-alignment and rate conditions.  The softmax
approximation provides a pointwise alternative under stronger
first-stage rates, while the post-selection confidence sets offer
calibrated simultaneous coverage for data-selected policy values.
On the efficiency side, we establish the trajectory-level bound and an
estimator attaining it in the enlarged episode model; deriving the
sharp efficiency bound under the structural stationarity restriction
$\mu_0=f_b$, which would require projecting onto the constrained
tangent space, is a natural direction for future work.

These results clarify both the scope and the limitations of existing approaches such as SAVE \citep{shi2022statistical}, and suggest that non-regularity is an intrinsic feature of optimal policy inference rather than a technical artifact.

\section*{Data Availability Statement}

The Drink Less data analyzed in this study were made available with the micro-randomized trial of \citet{bell2023notifications}.  The original trial article is available at \url{https://mhealth.jmir.org/2023/1/e38342}.

\section*{Disclosure Statement}

The authors report there are no competing interests to declare.

\section*{Supplementary Material}

The online Supplementary Material contains detailed configurations for the simulations and real data application (Appendix A), alternative confidence set constructions for post-selection inference (Appendix B), proofs of the theoretical results (Appendices C--E), and auxiliary lemmas (Appendix F).

\bibliographystyle{agsm}
\bibliography{ref}
\newpage 

\appendix

\section{Supplementary Material for Simulation and Real Data Application}

This appendix provides the full simulation specifications and additional
Monte Carlo diagnostics, followed by implementation details for the real
data application.

\subsection{Simulation Setup and Data-Generating Processes}\label{sec:sim_setup}

\paragraph{General setup.}
For every design, $\gamma=0.7$ and rewards are conditionally Bernoulli
with mean $r(s,a)$.  We first generate the behavior transition kernel
$K_b(s,s')=\sum_a b(a\mid s)P(s'\mid s,a)$ and set the initial law to its
stationary distribution $f_b$.  Consequently, the episodes are i.i.d. and
each within-episode state process is stationary under the behavior policy.
The exact finite-window score used in the computation is
\[
 \widehat V^\pi(S_0)
 +\frac{1}{1-\gamma^T}\sum_{t=0}^{T-1}\gamma^t
 \widehat\omega^\pi(S_t,A_t)
 \{R_t+\gamma\widehat V^\pi(S_{t+1})
                    -\widehat Q^\pi(S_t,A_t)\},
\]
where
$\omega^\pi(s,a)=d^\pi(s)\pi(a\mid s)/\{f_b(s)b(a\mid s)\}$.
Thus no cumulative trajectory importance ratio is used.  We use
$N\in\{100,250,500,1000\}$, $T\in\{25,50\}$, and 500 replications.

\paragraph{Fixed-policy DR and unique-optimum designs.}
Both designs use six states and two actions with a uniform behavior policy.
For every $(s,a)$, independent Gamma$(0.8,1)$ draws are normalized and mixed
with 0.08 times the uniform distribution to form $P(\cdot\mid s,a)$.
The reward means are generated independently from
Uniform$(0.15,0.85)$.  We retain a draw only if the minimum statewise gap
between the best and second-best actions is at least 0.06 in the DR design
and 0.08 in the unique-optimum design.  Fixed seeds are used so that all
sample sizes and horizons compare estimators on the same underlying MDP.

\paragraph{Exact-tie design.}
There is one state and two actions, with
$P(1\mid1,a)=1$, $b(1\mid1)=b(2\mid1)=1/2$, and
$r(1,1)=r(1,2)=0.5$.  The two deterministic policies therefore have the
same value.  Their fixed-policy scores nevertheless depend differently on
the observed action and reward, so the design represents
first-order distinguishable tied policies.

\paragraph{Near-tie design.}
This is again a one-state, two-action MDP with a uniform behavior policy,
but the reward means are 0.485 and 0.515.  The fixed action gap is 0.03.
This is a finite-sample stress test, not a local-to-zero triangular array.

\paragraph{Overlap designs.}
The strong- and weak-overlap designs use the same one-state MDP, with reward
means 0.425 and 0.575.  They differ only in the behavior probability assigned
to the optimal action: 0.50 under strong overlap and 0.15 under weak overlap.
The initial law and controlled dynamics are therefore identical, so changes
in variance and interval length can be attributed to overlap.

\subsection{Implementation Details and Additional Results}\label{sec:sim_impl}

\paragraph{Nuisance estimation and policy learning.}
All nuisance quantities are estimated tabularly.  We add 0.25 to every
transition cell, 0.50 to every behavior-policy action count, and one
pseudo-observation with reward mean 0.50 to every state-action reward cell.
The discounted occupancy ratio is then computed from the estimated
transition and behavior kernels.

For NSAVE, the episodes are randomly permuted and the initial training size
is
$\ell_N=\max(30,4|\mathcal S||\mathcal A|)$, truncated if necessary to leave
at least ten evaluation episodes.  The policy and nuisance functions are
refitted after every five evaluation episodes using only the preceding
data.  Within the preceding sample, the most recent 25\% is used to estimate
the predictable score scale and the remainder is used to fit the model.
The resulting scores are aggregated with inverse-scale weights.  The same
point estimate yields a one-sided 95\% lower bound (NSAVE--CLB) and a
two-sided 95\% interval (NSAVE--TCI).

The smoothed and hard-greedy benchmarks randomly split the episodes in half.
The first half estimates the policy and nuisance functions and the second
half evaluates the learned policy.  The smoothed policy is
\[
 \widehat\pi_{\beta_N}(a\mid s)
 \propto \exp\{\beta_N\widehat Q^*(s,a)\},
 \qquad
 \beta_N=\{\log(\lfloor N/2\rfloor)\}^2.
\]
After this policy is formed, its own
$\widehat Q^{\widehat\pi_{\beta_N}}$ and
$\widehat\omega^{\widehat\pi_{\beta_N}}$ are used in the DR score.

\begin{figure}[!htbp]
\centering
\begin{minipage}[t]{0.48\linewidth}
\centering
\vspace{0pt}
\includegraphics[width=\linewidth]{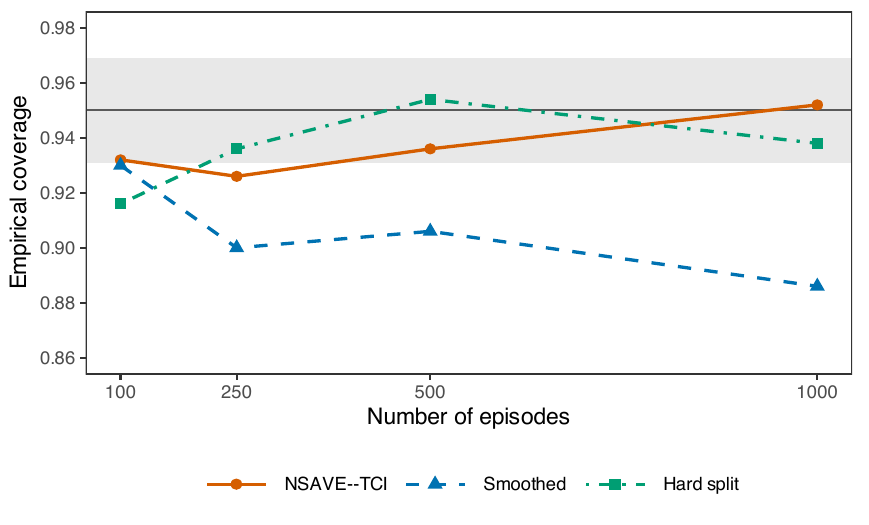}
\par\smallskip
{\small (a) Empirical coverage.}
\end{minipage}\hfill
\begin{minipage}[t]{0.48\linewidth}
\centering
\vspace{0pt}
\includegraphics[width=\linewidth]{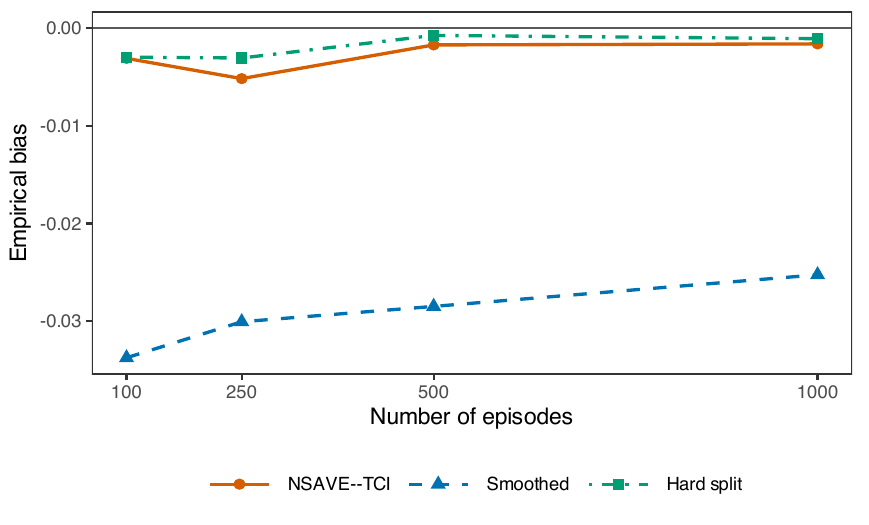}
\par\smallskip
{\small (b) Empirical bias.}
\end{minipage}
\caption{Near-tie stress test for $T=50$.  The gray band in panel (a) is
the 95\% Monte Carlo uncertainty band around 0.95.}
\label{fig:mc_near_tie}
\end{figure}

\paragraph{Near-tie results.}
Figure~\ref{fig:mc_near_tie} separates inferential error from smoothing
approximation bias.  At $T=50$, NSAVE--TCI coverage increases from 0.932 at
$N=100$ to 0.952 at $N=1000$, and its absolute bias is below 0.006 at every
sample size.  The hard-split benchmark has similarly small bias.  By
contrast, the chosen finite-temperature schedule leaves smoothing bias
between $-0.034$ and $-0.025$, and smoothing coverage falls from 0.930 to
0.886.  The associated policy regret remains 0.024 at $N=1000$.  This is
not a contradiction to pointwise smoothing theory: it shows that, at these
sample sizes and this small action gap, the smoothing approximation bias is
not yet negligible relative to the standard error.  Faster annealing would
reduce this bias but can make policy learning more sensitive to estimation
noise.

\begin{figure}[!htbp]
\centering
\includegraphics[width=0.78\linewidth]{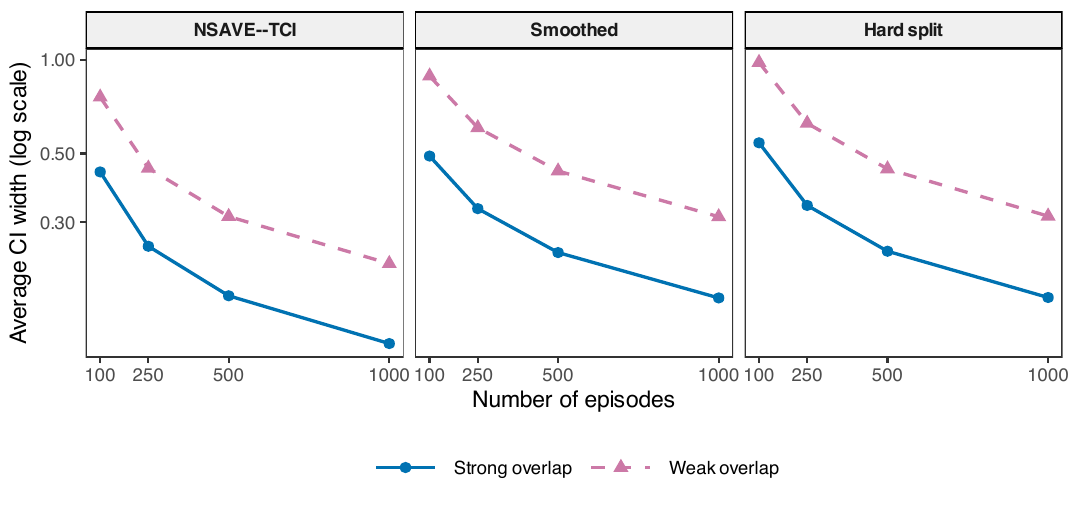}
\caption{Average two-sided confidence-interval width under strong and weak
overlap for $T=50$.}
\label{fig:mc_overlap}
\end{figure}

\paragraph{Overlap results.}
Figure~\ref{fig:mc_overlap} shows that reducing the behavior probability of
the optimal action from 0.50 to 0.15 increases the average interval width by
roughly a factor of 1.8 for all three two-sided procedures.  For example,
the NSAVE--TCI width at $N=1000$ increases from 0.122 to 0.221.
Coverage is largely preserved as $N$ grows, although weak overlap produces
visible small-sample undercoverage for NSAVE--TCI: 0.908 at $N=100$,
improving to 0.944 and 0.940 at $N=500$ and $N=1000$, respectively.  Thus
the experiment displays the variance cost of weak overlap without
confounding it with a change in the state distribution.

\begin{figure}[!htbp]
\centering
\begin{minipage}[t]{0.48\linewidth}
\centering
\vspace{0pt}
\includegraphics[width=\linewidth]{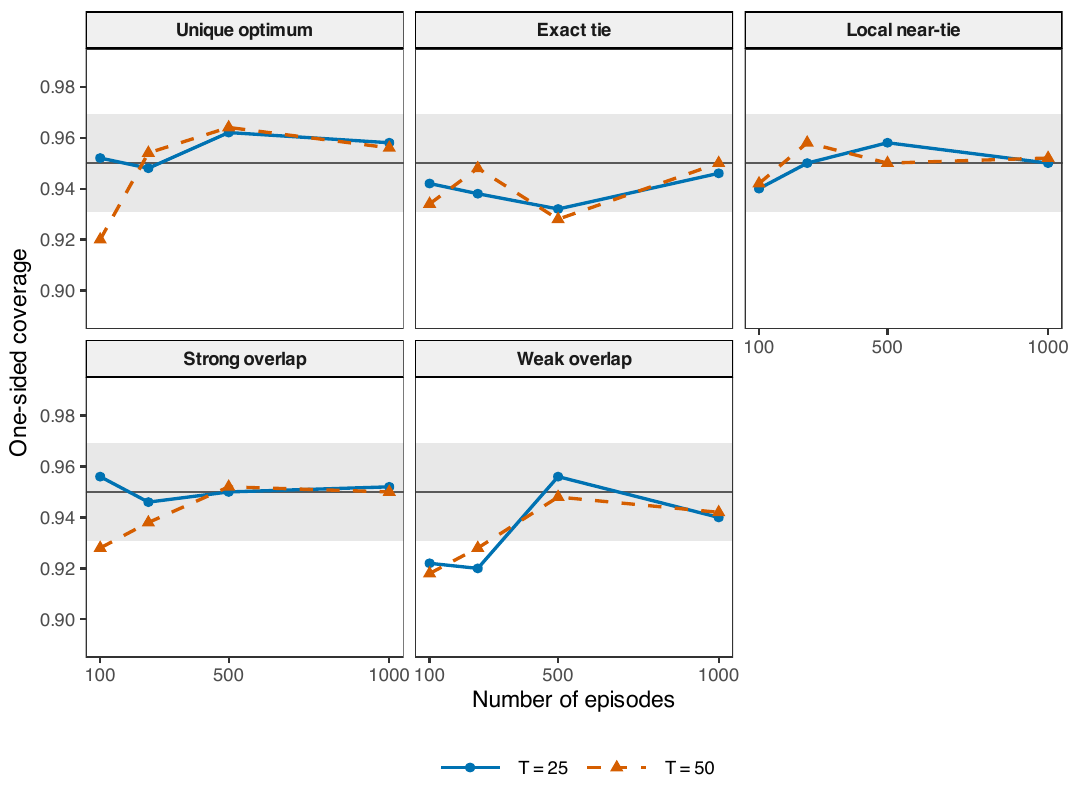}
\par\smallskip
{\small (a) One-sided coverage of NSAVE--CLB.}
\end{minipage}\hfill
\begin{minipage}[t]{0.48\linewidth}
\centering
\vspace{0pt}
\includegraphics[width=\linewidth]{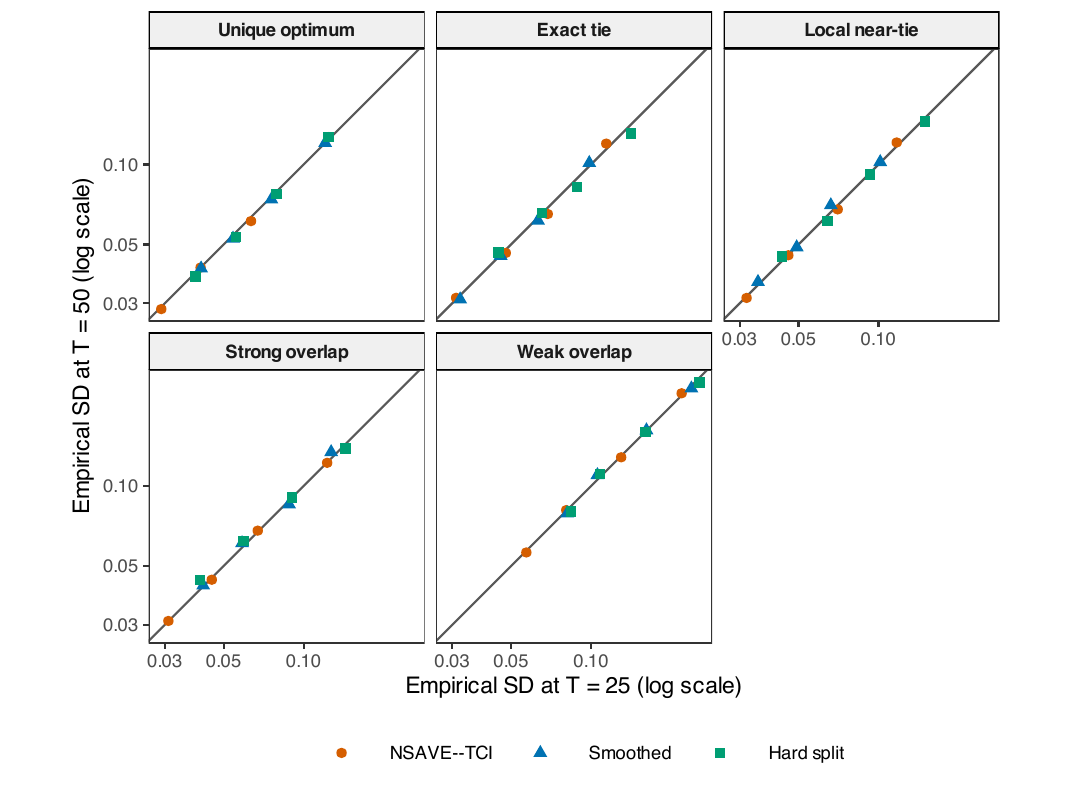}
\par\smallskip
{\small (b) Horizon sensitivity of the empirical standard deviation.}
\end{minipage}
\caption{Additional Monte Carlo diagnostics.  Panel (a) shows one-sided
coverage of NSAVE--CLB; the gray band is the 95\% Monte Carlo uncertainty
band around 0.95.  Panel (b) compares empirical standard deviations across
horizons; points on the diagonal have equal empirical standard deviations
at $T=25$ and $T=50$.}
\label{fig:mc_clb}
\label{fig:mc_horizon}
\end{figure}

\paragraph{One-sided coverage and horizon sensitivity.}
Panel (a) of Figure~\ref{fig:mc_clb} reports NSAVE--CLB coverage for both horizons.  The
minimum coverage over all displayed designs is 0.918 and the maximum is
0.964.  The small differences between $T=25$ and $T=50$ are expected here:
because the return and score are normalized and $\gamma=0.7$, both horizons
already exceed the effective discount horizon.  Across the primary
two-sided procedures, the mean ratio of empirical standard deviations at
$T=50$ and $T=25$ is 0.996, with range 0.923--1.085; the average absolute
coverage difference is 0.010.  Panel (b) of Figure~\ref{fig:mc_horizon} gives the
corresponding pointwise comparison.

\begin{figure}[!htbp]
\centering
\includegraphics[width=0.74\linewidth]{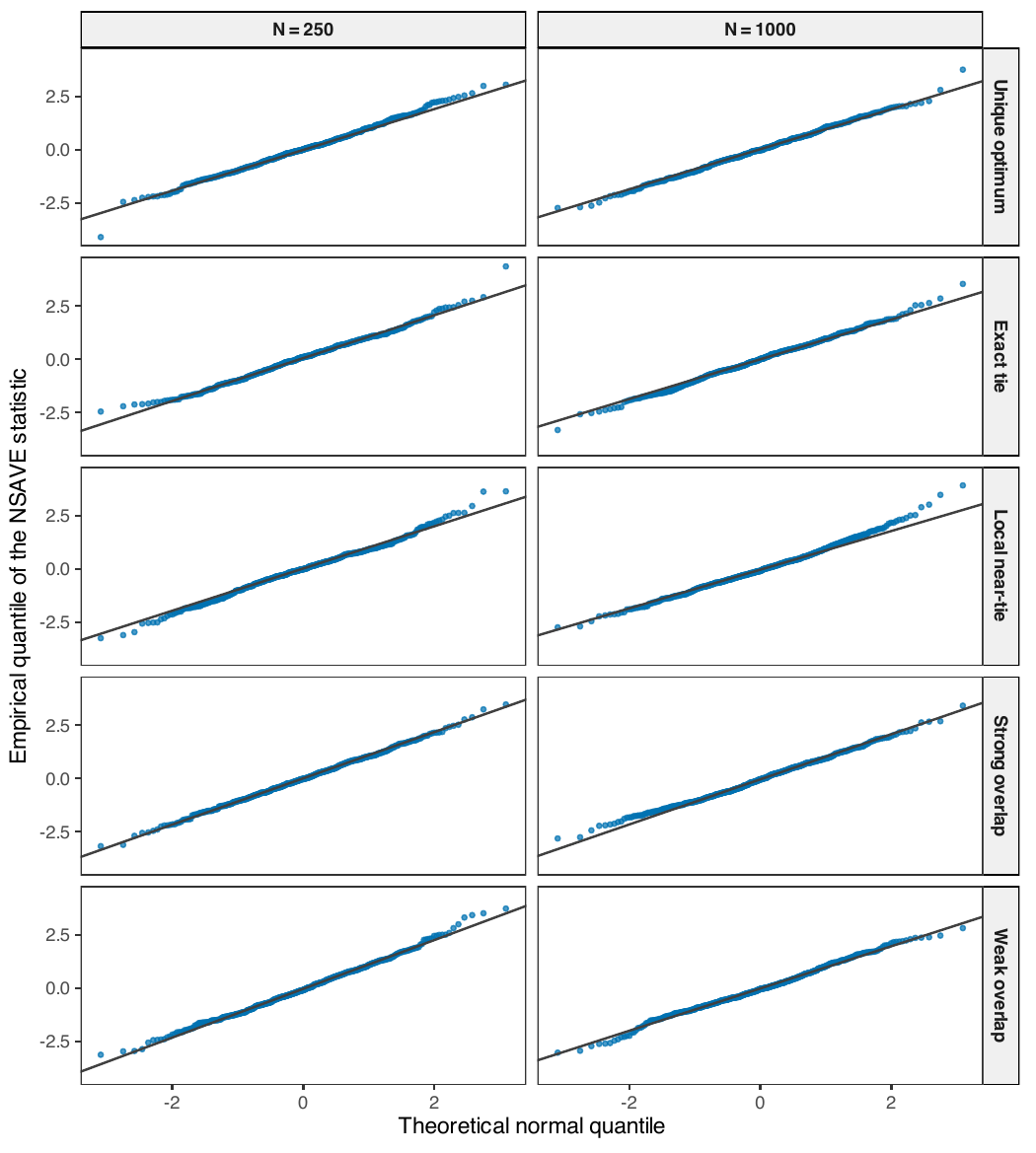}
\caption{Normal Q--Q plots for the studentized NSAVE statistic with $T=50$.}
\label{fig:mc_qq}
\end{figure}

\paragraph{Studentized NSAVE statistic.}
Finally, Figure~\ref{fig:mc_qq} compares the empirical quantiles of
\[
 Z_{\mathrm{hist}}
 =
 \frac{\widehat\eta_{\mathrm{NSAVE}}
       -\overline\eta_{\mathrm{historical}}}
      {\widehat{\mathrm{se}}_{\mathrm{NSAVE}}}
\]
with standard normal quantiles.  The historical target is the
inverse-scale-weighted value of the predictable policies actually evaluated
by NSAVE.  The plots are close to linear by $N=250$ and improve further at
$N=1000$; the largest finite-sample departure occurs under weak overlap.

\subsection{Drink Less Data Construction and Additional Results}
\label{sec:drinkless_setup}

\paragraph{Data source and action coding.}
We use the shared data from the Drink Less MRT of \citet{bell2023notifications}.
The analysis is based on Dataset A, which contains $349$ participants observed for all $30$ daily decision points.
The decision point is 8pm.
The three actions are coded from the randomized notification version:
\texttt{Z} corresponds to no notification, \texttt{X} to a new notification, and \texttt{Y} to the standard notification.
The known behavior probabilities are $0.40$, $0.30$, and $0.30$, respectively, and are used directly in all importance-weight and marginal-occupancy calculations.
We do not estimate the behavior policy.

\paragraph{MDP construction.}
The reward $R_t$ is the primary binary outcome indicating whether the participant opened the app from 8pm to 9pm on day $t$.
The primary analysis uses days 15--30, after a 14-day burn-in period.
Transitions are formed from days 15--29, with day 30 supplying the final next state, so the effective horizon is $L=15$.
The target initial-state distribution is the pooled behavior stationary distribution estimated over the analysis window, rather than the empirical day-15 distribution.
This matches the stationary experiment used in Assumption~\ref{ass:stationary}.

The primary state representation, denoted \texttt{history18}, has 18 levels of the form \texttt{B\#\_E\#\_N\#}.
Here \texttt{B} is whether the participant engaged with the app before 8pm on the current day; \texttt{E} is the number of previous three days with any recorded engagement, capped at two; and \texttt{N} is the number of previous three days on which any notification was sent, also capped at two.
The engagement count uses app openings before 8pm, during the primary 8pm--9pm window, and after 9pm on the corresponding previous day.
As a transparent sensitivity analysis, we also use a 16-level \texttt{core16} state that records current pre-8pm engagement, the previous primary outcome, whether a notification was sent on the previous day, and previous post-9pm engagement.

\paragraph{Nuisance fitting and policy learning.}
For a fixed policy, we fit a finite-state MDP with hierarchical tabular smoothing.
Each reward mean is shrunk toward the action-specific global reward rate, and each transition row is shrunk toward the action-specific global next-state distribution.
The primary shrinkage strengths are $8$ for rewards and $5$ for transitions.
The marginal state-action occupancy ratio is then computed from the fitted transition kernel and the policy being evaluated.
The resulting episode score is the normalized finite-window version of the marginal-occupancy doubly robust score, with normalization $1-\gamma^L$.
No cumulative trajectory importance ratios are used.

The primary NSAVE run uses an initial training size $\ell_n=100$, a rolling scale-estimation window of 30 participants, and policy refitting every five participants.
The participant ordering is fixed by the reproducibility seed used in the analysis.
All reported NSAVE intervals are therefore conditional on this fixed predictable ordering.
The cross-fitted hard-greedy and smoothed analyses use five folds.
For the smoothed policy, the primary inverse-temperature is $\beta=\{\log(N_{\mathrm{train}})\}^2$.

\begin{figure}[!htbp]
    \centering
    \begin{minipage}[t]{0.48\linewidth}
    \centering
    \vspace{0pt}
    \includegraphics[width=\linewidth]{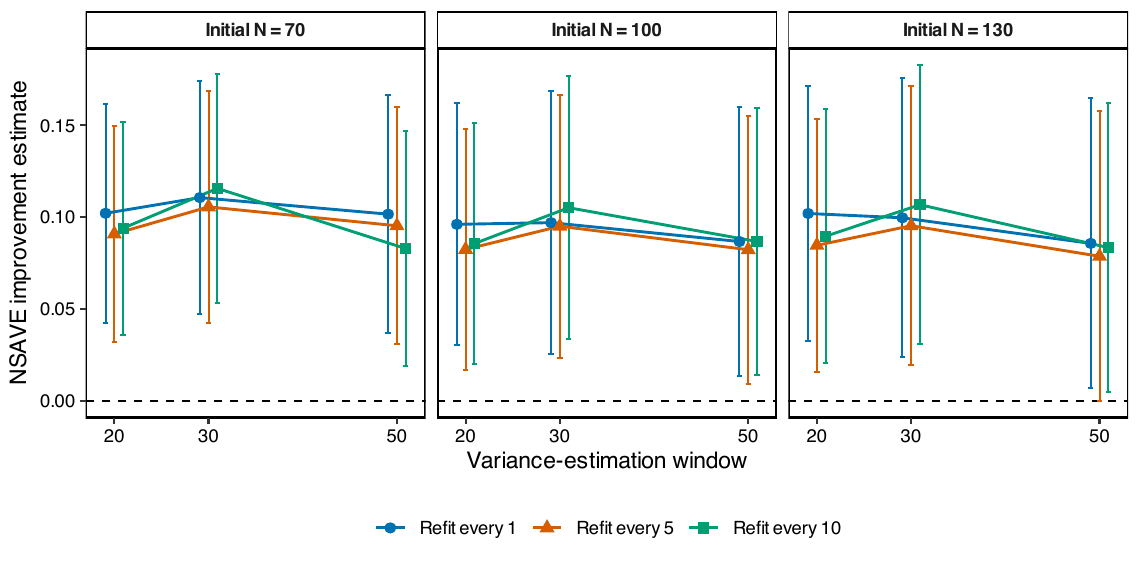}
    \par\smallskip
    {\small (a) Tuning sensitivity.}
    \end{minipage}\hfill
    \begin{minipage}[t]{0.48\linewidth}
    \centering
    \vspace{0pt}
    \includegraphics[width=\linewidth]{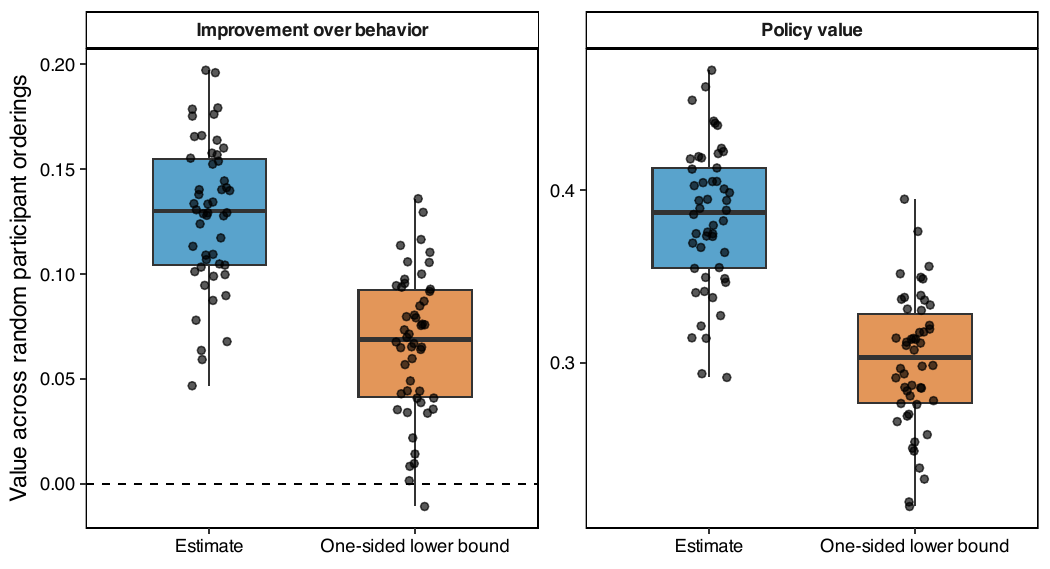}
    \par\smallskip
    {\small (b) Participant-order sensitivity.}
    \end{minipage}
    \caption{NSAVE sensitivity analyses in the Drink Less application.
    Panel (a) varies the initial training size, scale window, and refitting
    frequency; intervals show 95\% two-sided confidence intervals.  Panel
    (b) varies the random participant ordering; each point corresponds to
    one ordering.}
    \label{fig:drinkless_tuning}
    \label{fig:drinkless_order}
\end{figure}

\paragraph{NSAVE tuning sensitivity.}
We vary the initial training size over $\{70,100,130\}$, the scale window over $\{20,30,50\}$, and the refitting frequency over $\{1,5,10\}$.
Across the resulting 27 configurations, the NSAVE improvement estimates range from $0.079$ to $0.116$, and all 27 one-sided lower confidence bounds are positive.
Panel (a) of Figure~\ref{fig:drinkless_tuning} shows that the primary conclusion is not driven by the selected tuning triple.

\paragraph{Participant-order sensitivity.}
Because NSAVE evaluates predictable learned policies, the historical sequence of learned policies depends on the participant ordering.
Across 50 random orderings, the median improvement estimate is $0.130$, with 2.5th--97.5th percentiles $[0.060,0.192]$.
The median one-sided lower bound is $0.069$, with 2.5th--97.5th percentiles $[0.003,0.127]$; 49 of 50 orderings have a positive lower bound.
Panel (b) of Figure~\ref{fig:drinkless_order} reports the distribution of estimates and lower bounds.
This sensitivity supports the value-improvement conclusion while also illustrating that the identity of the selected near-tied actions should not be overinterpreted.

\begin{figure}[!htbp]
    \centering
    \begin{minipage}[t]{0.58\linewidth}
    \centering
    \vspace{0pt}
    \includegraphics[width=\linewidth]{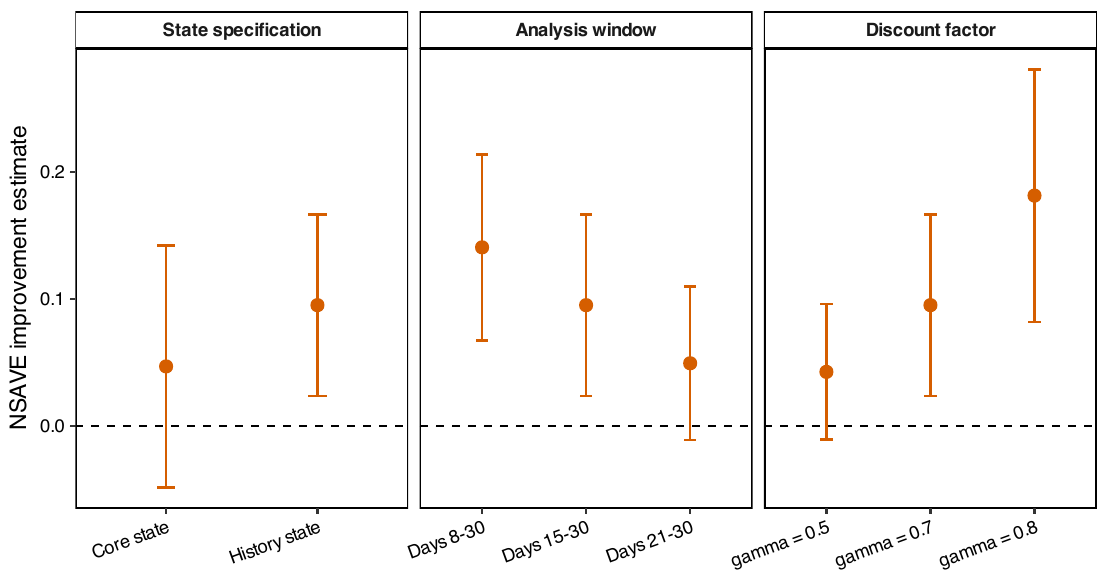}
    \par\smallskip
    {\small (a) Structural sensitivity.}
    \end{minipage}\hfill
    \begin{minipage}[t]{0.38\linewidth}
    \centering
    \vspace{0pt}
    \includegraphics[width=\linewidth]{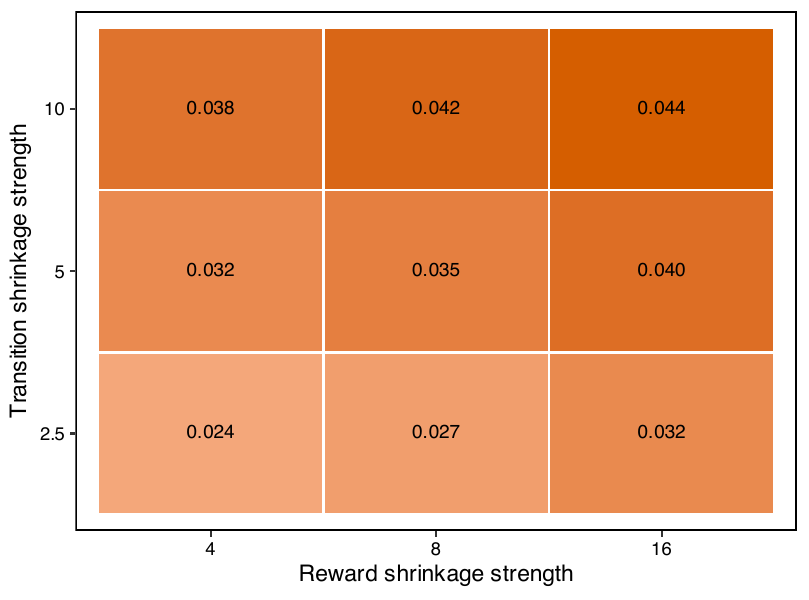}
    \par\smallskip
    {\small (b) Sparse-cell shrinkage sensitivity.}
    \end{minipage}
    \caption{Additional Drink Less sensitivity analyses.  Panel (a)
    summarizes structural sensitivity across state representations,
    analysis windows, and discount factors.  Panel (b) shows sensitivity of
    the NSAVE one-sided lower bound to transition and reward shrinkage
    strengths.}
    \label{fig:drinkless_structural}
    \label{fig:drinkless_prior}
\end{figure}

\paragraph{State, window, and discount-factor sensitivity.}
Panel (a) of Figure~\ref{fig:drinkless_structural} summarizes the main structural sensitivity analyses.
Using the coarser \texttt{core16} state gives an improvement estimate $0.047$ with one-sided lower bound $-0.033$, whereas the primary \texttt{history18} state gives $0.095$ with lower bound $0.035$.
For analysis windows days 8--30, 15--30, and 21--30, the improvement estimates are $0.141$, $0.095$, and $0.049$, with one-sided lower bounds $0.079$, $0.035$, and $-0.001$.
For $\gamma=0.5,0.7,0.8$, the corresponding improvement estimates are $0.043$, $0.095$, and $0.181$.
These results show that state construction and horizon choice are substantive modeling decisions; the main text therefore emphasizes policy-value evidence rather than a unique policy recommendation.

\paragraph{Sparse-cell shrinkage sensitivity.}
The primary analysis regularizes sparse tabular cells using transition shrinkage strength $5$ and reward shrinkage strength $8$.
We repeat NSAVE over a $3\times 3$ grid of transition strengths $\{2.5,5,10\}$ and reward strengths $\{4,8,16\}$.
The one-sided lower bound for the paired improvement ranges from $0.024$ to $0.044$ across this grid, so the positive primary lower bound is not an artifact of a single shrinkage choice.
Panel (b) of Figure~\ref{fig:drinkless_prior} reports the full grid.

\begin{figure}[!htbp]
    \centering
    \includegraphics[width=0.4\linewidth]{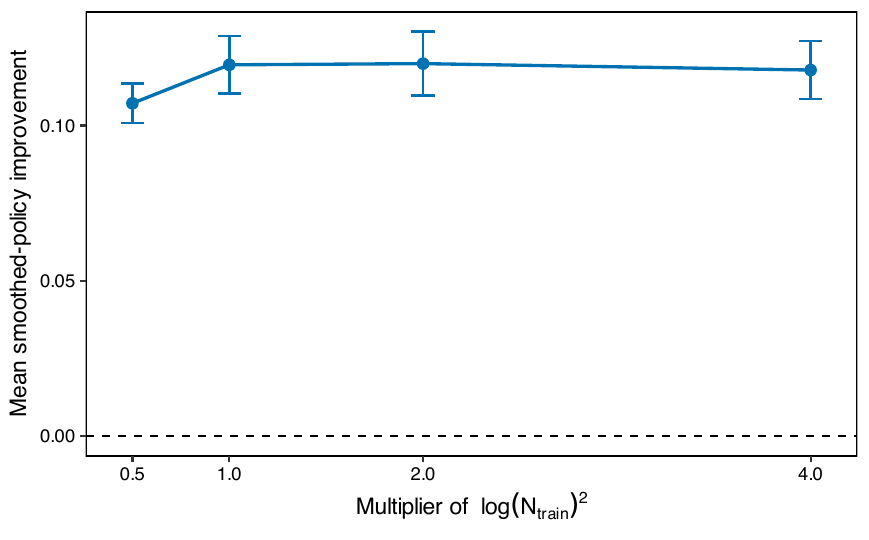}
    \caption{Sensitivity of the cross-fitted smoothed policy to the inverse-temperature multiplier.}
    \label{fig:drinkless_smoothing_beta}
\end{figure}

\paragraph{Smoothing sensitivity and near non-uniqueness.}
We also repeat the smoothed cross-fitted analysis over inverse-temperature multipliers $0.5$, $1$, $2$, and $4$ relative to the primary choice $\{\log(N_{\mathrm{train}})\}^2$.
The mean paired improvement across 20 random splits ranges from $0.107$ to $0.120$, while the split-to-split standard deviation remains substantially smaller than for hard-greedy cross-fitting.
This is the empirical pattern expected under near non-uniqueness: smoothing stabilizes the learned action probabilities without materially changing the estimated value.
Figure~\ref{fig:drinkless_smoothing_beta} reports the smoothing sensitivity.

\section{Alternative Confidence Set Constructions for Post-Selection Inference}\label{sec:psi_other_methods}

In this section we summarize several confidence set constructions for post-selection inference (PSI)
in our setting.
Let $\bm{Z}=(Z_1,\dots,Z_K)^\top \rightsquigarrow \mathcal{N}(\bm{0},\boldsymbol{R})$ with
$\boldsymbol{R}:=\mathrm{diag}(\boldsymbol{\Sigma})^{-1/2}\boldsymbol{\Sigma}\,\mathrm{diag}(\boldsymbol{\Sigma})^{-1/2}$.

\paragraph{Gaussian critical values.}
In all methods below, we compute critical values via the Gaussian approximation induced by
$\widehat{\boldsymbol{\Sigma}}$.
We assume the reported coordinates have asymptotic variances bounded
away from zero; any exactly degenerate coordinate should instead be
handled separately.
Let
\[
\widehat{\boldsymbol{R}}
:=
\mathrm{diag}(\widehat{\boldsymbol{\Sigma}})^{-1/2}\widehat{\boldsymbol{\Sigma}}\,
\mathrm{diag}(\widehat{\boldsymbol{\Sigma}})^{-1/2},
\qquad
\bm{G}\sim \mathcal{N}(\bm{0},\widehat{\boldsymbol{R}}).
\]
For any index set $S\subseteq[K]$ and level $u\in(0,1)$, define
\begin{equation}\label{eq:q_def_general}
q_u(S)
\ :=\ 
\inf\Bigl\{t\in\mathbb{R}:\ \pr\big(\max_{k\in S}|G_k|\le t \ \bigm|\ \widehat{\boldsymbol{R}}\big)\ge u\Bigr\},
\end{equation}
which can be approximated by Monte Carlo simulation from $\mathcal{N}(\bm{0},\widehat{\boldsymbol{R}})$.

\subsection{Projection (global simultaneous inference)}\label{subsec:proj}
The projection approach ignores the selection event and instead provides a uniform (simultaneous)
guarantee over all $K$ coordinates. It corresponds to taking $\delta_1 \equiv 0$ and calibrating the critical value
against the full maximum.

\paragraph{Critical value.}
Set $q^{\mathrm{proj}}_{1-\delta_2}:=q_{1-\delta_2}([K]).$

\paragraph{Confidence set.}
Define
\begin{equation*}\label{eq:CS_proj}
\mathcal{C}_{\mathrm{proj}}
:=
\bigtimes_{k\in\widehat{\mathcal{A}}_{\mathrm{opt}}}
\Bigl[
\widehat{\eta}_k \ \pm\ q^{\mathrm{proj}}_{1-\delta_2}\cdot \sqrt{\widehat{\boldsymbol{\Sigma}}_{kk}/N}
\Bigr].
\end{equation*}
This method is always valid (asymptotically) under the joint Gaussian approximation, but is typically
conservative when $K$ is large.

\subsection{Locally simultaneous inference}\label{subsec:localsim}
Locally simultaneous inference first constructs a high-probability superset of policies that could be
optimal, and then calibrates a simultaneous critical value over this smaller set.

\paragraph{Step 1 (plausible-optimal superset).}
Fix $\delta_1\in(0,\delta_2)$. Construct marginal confidence bounds
\[
\mathrm{UCB}_k:=\widehat{\eta}_k + c_{1-\delta_1}\sqrt{\widehat{\boldsymbol{\Sigma}}_{kk}/N},
\qquad
\mathrm{LCB}_k:=\widehat{\eta}_k - c_{1-\delta_1}\sqrt{\widehat{\boldsymbol{\Sigma}}_{kk}/N},
\]
where one may take the conservative \textbf{Bonferroni} choice
$c_{1-\delta_1}:=z_{1-\delta_1/(2K)}$ (alternatively one may use a Gaussian max-quantile over $[K]$ at level
$1-\delta_1$).
Define the plausible-optimal set
\begin{equation*}\label{eq:A_plus}
\widehat{\mathcal{A}}^{+}
:=
\Bigl\{k\in[K]:\ \mathrm{UCB}_k \ge \max_{\ell\in[K]}\mathrm{LCB}_\ell\Bigr\}.
\end{equation*}
Intuitively, $\widehat{\mathcal{A}}^{+}$ removes policies whose upper confidence bound lies below the
lower confidence bound of (at least) one competitor, so they cannot be optimal within the $(1-\delta_1)$
uncertainty set.

\paragraph{Step 2 (local simultaneous calibration).}
Set
$
q^{\mathrm{LS}}_{1-(\delta_2-\delta_1)}:=q_{1-(\delta_2-\delta_1)}\bigl(\widehat{\mathcal{A}}^{+}\bigr).
$
Since $\widehat{\mathcal A}^{+}$ is data dependent, this plug-in
definition is justified only when the locally simultaneous
construction verifies the joint calibration condition stated in
Step~2 of Section~\ref{sec:PSI}.  A Gaussian quantile for the realized
set, without that additional result, is not by itself a coverage
proof.  The projection critical value over $[K]$ is the conservative
fallback.

\paragraph{Confidence set.}
Define
\begin{equation*}\label{eq:CS_LS}
\mathcal{C}_{\mathrm{LS}}
:=
\bigtimes_{k\in\widehat{\mathcal{A}}_{\mathrm{opt}}}
\Bigl[
\widehat{\eta}_k \ \pm\ q^{\mathrm{LS}}_{1-(\delta_2-\delta_1)}\cdot \sqrt{\widehat{\boldsymbol{\Sigma}}_{kk}/N}
\Bigr].
\end{equation*}
When the required joint calibration condition holds and
$\widehat{\mathcal{A}}^{+}$ is substantially smaller than $[K]$,
\eqref{eq:CS_LS} can be much tighter than projection while controlling
the overall error by splitting
$\delta_2=\delta_1+(\delta_2-\delta_1)$.

\subsection{Hybrid constructions}\label{subsec:hybrid}
Hybrid methods combine a global ``safety net'' region with a sharper selective/local procedure.
We present a practical hybrid that is simple to implement and guarantees that the resulting interval
is never wider than the global projection band.

\paragraph{Hybrid by intersection.}
Fix $\delta_1\in(0,\delta_2)$.
Compute the projection critical value at level $1-\delta_1$,
$
q^{\mathrm{proj}}_{1-\delta_1}:=q_{1-\delta_1}([K]),
$
and compute a selective/local critical value, e.g.\ $q^{\mathrm{LS}}_{1-(\delta_2-\delta_1)}$ from
Section~\ref{subsec:localsim} (one may replace it by the two-step critical value in
Section~\ref{sec:PSI} if desired). Define the coordinate-wise radius
\[
r_k^{\mathrm{hyb}}
:=
\min\Bigl\{
q^{\mathrm{proj}}_{1-\delta_1},\ q^{\mathrm{LS}}_{1-(\delta_2-\delta_1)}
\Bigr\}\cdot \sqrt{\widehat{\boldsymbol{\Sigma}}_{kk}/N},
\]
and the hybrid confidence set
\begin{equation*}\label{eq:CS_hybrid}
\mathcal{C}_{\mathrm{hyb}}
:=
\bigtimes_{k\in\widehat{\mathcal{A}}_{\mathrm{opt}}}
\Bigl[
\widehat{\eta}_k \ \pm\ r_k^{\mathrm{hyb}}
\Bigr].
\end{equation*}
By construction, $\mathcal{C}_{\mathrm{hyb}}$ is never wider than the
projection band at level $1-\delta_1$.  Its asserted coverage requires
both component procedures to satisfy their respective allocated error
bounds; under those bounds, the union bound gives total error at most
$\delta_2$.  Without a valid local/selective calibration, taking the
smaller radius need not preserve coverage.

\subsection{Conditional selective inference}\label{subsec:conditional}
Conditional selective inference calibrates inference \emph{given} the selection event
$\{\widehat{\mathcal{A}}_{\mathrm{opt}}=a\}$ rather than via a worst-case bound. When
$\widehat{\mathcal{A}}_{\mathrm{opt}}=\{k^\star\}$ is a singleton (unique empirical winner), the selection event
can be written as the polyhedral constraint
\begin{equation*}\label{eq:selection_polyhedron}
\widehat{\eta}_{k^\star} \ge \widehat{\eta}_\ell,\qquad \forall \ell\neq k^\star,
\end{equation*}
and under the Gaussian approximation
$\widehat{\bm{\eta}}\approx \mathcal{N}(\bm{\eta},\boldsymbol{\Sigma}/N)$, the conditional law of
$\widehat{\bm{\eta}}$ given \eqref{eq:selection_polyhedron} is a truncated multivariate normal over a polyhedral
cone.

\paragraph{Generic test inversion.}
Let $a$ denote the realized selection outcome (e.g.\ $a=\{k^\star\}$).
For a candidate parameter vector $\bm{\eta}'$, let $\pr_{\bm{\eta}'}(\cdot\mid \widehat{\mathcal{A}}_{\mathrm{opt}}=a)$
denote the induced conditional probability under the Gaussian model. Define a family of tests
$\{\varphi_{\bm{\eta}'}\}$ with conditional size control,
\[
\pr_{\bm{\eta}'}\bigl(\varphi_{\bm{\eta}'}=1\mid \widehat{\mathcal{A}}_{\mathrm{opt}}=a\bigr)\le \delta_2,
\]
and define the conditional confidence set by inversion,
\begin{equation}\label{eq:CS_conditional_generic}
\mathcal{C}_{\mathrm{cond}}(a)
:=
\Bigl\{\bm{\eta}'\in\mathbb{R}^K:\ \varphi_{\bm{\eta}'}=0\Bigr\}.
\end{equation}
A confidence interval for the selected coordinate(s) is then obtained by projecting
$\mathcal{C}_{\mathrm{cond}}(a)$ onto $\{\eta_k:k\in a\}$.

While \eqref{eq:CS_conditional_generic} provides the conceptually sharpest adjustment, implementing it for large
$K$ (and/or non-unique winners) typically requires nontrivial computation for truncated multivariate normals and
test inversion.
If $|\widehat{\mathcal{A}}_{\mathrm{opt}}|>1$, the selection event is the union of polyhedral regions (or can be
expressed via additional constraints encoding ties), which further increases computational complexity for exact
conditional inference. The max-based procedures above naturally accommodate ties by reporting intervals for all
$k\in\widehat{\mathcal{A}}_{\mathrm{opt}}$.

\section{Technical Proofs in Section \ref{sec:IFF_for_EIF}}

\subsection{Occupancy stability and the absence of a converse bound}

Write
\[
    d_P^\pi := (1-\gamma)\sum_{t=0}^{\infty}\gamma^t
    \mathcal{L}_P^\pi(S_t)
\]
for the discounted state-occupancy probability measure.  When
$d_P^\pi$ has a density relative to the stationary state law, that
density is denoted
$\rho_P^\pi(s):=\mathrm d d_P^\pi/\mathrm d f_b(s)$.  The
state--action marginal ratio used in the estimating score is then
$\omega_P^\pi(a,s)
=\rho_P^\pi(s)\pi(a\mid s)/b(a\mid s)$.
The policy-induced state-transition kernel is
\begin{equation}\label{eq:transition_kernel}
    K_{\pi,P}(B\mid s)
    :=\int_{\mathcal A}P(S'\in B\mid s,a)\pi(\mathrm da\mid s).
\end{equation}

\begin{lemma}[One-Sided Occupancy Stability]
\label{lem:diff_polcies_upper}
    For two policies under the same law $P$,
    \[
        \|d_P^{\pi_2}-d_P^{\pi_1}\|_{\mathrm{TV}}
        \leq
        \frac{\gamma}{1-\gamma}
        \E_{S\sim d_P^{\pi_1}}
        \!\left[\operatorname{TV}
        \{\pi_2(\cdot\mid S),\pi_1(\cdot\mid S)\}\right].
    \]
    Equivalently, when the occupancy measures admit densities,
    \[
        \|\rho_P^{\pi_2}-\rho_P^{\pi_1}\|_{L_1(f_b)}
        \leq
        \frac{2\gamma}{1-\gamma}
        \E_{S\sim d_P^{\pi_1}}
        \!\left[\operatorname{TV}
        \{\pi_2(\cdot\mid S),\pi_1(\cdot\mid S)\}\right].
    \]
    If $\nu_P^\pi(\mathrm ds,\mathrm da)
    :=d_P^\pi(\mathrm ds)\pi(\mathrm da\mid s)$ denotes the discounted
    state--action occupancy measure, then
    \[
        \|\nu_P^{\pi_2}-\nu_P^{\pi_1}\|_{\mathrm{TV}}
        \leq
        \frac{1}{1-\gamma}
        \E_{S\sim d_P^{\pi_1}}
        \!\left[\operatorname{TV}
        \{\pi_2(\cdot\mid S),\pi_1(\cdot\mid S)\}\right].
    \]
\end{lemma}

\begin{proof}
    We regard probability measures as row vectors and kernels as
    operators acting on their right.  For $k=1,2$, abbreviate
    $K_k:=K_{\pi_k,P}$ and $d_k:=d_P^{\pi_k}$.  By the definition of
    discounted occupancy,
    \[
        \begin{aligned}
        d_k
        &=(1-\gamma)\sum_{t=0}^{\infty}
          \gamma^t\mu_0K_k^t\\
        &=(1-\gamma)\mu_0+
          \gamma(1-\gamma)\sum_{t=0}^{\infty}
          \gamma^t\mu_0K_k^{t+1}\\
        &=(1-\gamma)\mu_0+\gamma d_kK_k.
        \end{aligned}
    \]
    Subtract the two flow equations:
    \[
        d_2-d_1
        =\gamma(d_2K_2-d_1K_1).
    \]
    Adding and subtracting $\gamma d_1K_2$ gives
    \[
        (d_2-d_1)(I-\gamma K_2)
        =\gamma d_1(K_2-K_1).
    \]
    Since $\gamma<1$, the inverse exists as the Neumann series
    \[
        (I-\gamma K_2)^{-1}
        =\sum_{m=0}^{\infty}\gamma^mK_2^m,
    \]
    and hence
    \begin{equation}\label{eq:occupancy_resolvent_difference}
        d_2-d_1
        =\gamma d_1(K_2-K_1)
        \sum_{m=0}^{\infty}\gamma^mK_2^m.
    \end{equation}

    We next bound each factor.  We use the convention
    $\|\mu-\nu\|_{\mathrm{TV}}
      =\frac12\sup_{\|f\|_\infty\leq1}
      |\int f\,\mathrm d(\mu-\nu)|$ for differences of probability
      measures, and use the same half-variation norm for finite
      signed measures of total mass zero.

    If $\xi$ is a finite signed measure of total mass zero and $K$
    is a Markov kernel, then
    \[
        \|\xi K\|_{\mathrm{TV}}
        =\frac12\sup_{\|f\|_\infty\leq1}
          \left|\int Kf\,\mathrm d\xi\right|
        \leq\|\xi\|_{\mathrm{TV}},
    \]
    because $\|Kf\|_\infty\leq\|f\|_\infty$.  Therefore,
    \[
        \left\|\xi(I-\gamma K_2)^{-1}\right\|_{\mathrm{TV}}
        \leq\sum_{m=0}^{\infty}\gamma^m
        \|\xi K_2^m\|_{\mathrm{TV}}
        \leq\frac{\|\xi\|_{\mathrm{TV}}}{1-\gamma}.
    \]
    For each state $s$,
    \[
        K_2(\cdot\mid s)-K_1(\cdot\mid s)
        =
        \int_{\mathcal A}
        P(\cdot\mid s,a)
        \{\pi_2-\pi_1\}(\mathrm da\mid s).
    \]
    Hence, using the same dual formula,
    \[
        \begin{aligned}
        \|K_2(\cdot\mid s)-K_1(\cdot\mid s)\|_{\mathrm{TV}}
        &=
        \frac12\sup_{\|f\|_\infty\leq1}
        \left|
        \int P f(s,a)\{\pi_2-\pi_1\}(\mathrm da\mid s)
        \right|\\
        &\leq
        \operatorname{TV}\{\pi_2(\cdot\mid s),
        \pi_1(\cdot\mid s)\},
        \end{aligned}
    \]
    because $|Pf(s,a)|\leq1$.  Integrating with respect to $d_1$
    yields
    \[
        \|d_1(K_2-K_1)\|_{\mathrm{TV}}
        \leq
        \int
        \operatorname{TV}\{\pi_2(\cdot\mid s),
        \pi_1(\cdot\mid s)\}\,d_1(s).
    \]
    Applying the last two bounds to
    \eqref{eq:occupancy_resolvent_difference} proves
    \[
        \|d_2-d_1\|_{\mathrm{TV}}
        \leq\frac{\gamma}{1-\gamma}
        \E_{S\sim d_1}\!\left[
        \operatorname{TV}\{\pi_2(\cdot\mid S),
        \pi_1(\cdot\mid S)\}\right].
    \]

    If $d_k$ has density $\rho_k$ relative to $f_b$, then
    \[
        \|\rho_2-\rho_1\|_{L_1(f_b)}
        =2\|d_2-d_1\|_{\mathrm{TV}},
    \]
    which proves the density bound.

    Finally, introduce the intermediate state--action measure
    $\widetilde\nu(\mathrm ds,\mathrm da)
      :=d_1(\mathrm ds)\pi_2(\mathrm da\mid s)$.  The triangle
    inequality gives
    \[
        \|\nu_2-\nu_1\|_{\mathrm{TV}}
        \leq\|\nu_2-\widetilde\nu\|_{\mathrm{TV}}
        +\|\widetilde\nu-\nu_1\|_{\mathrm{TV}}.
    \]
    The first term equals $\|d_2-d_1\|_{\mathrm{TV}}$, because both
    measures use the same conditional action kernel $\pi_2$.  The
    second term equals
    \[
        \int
        \operatorname{TV}\{\pi_2(\cdot\mid s),
        \pi_1(\cdot\mid s)\}\,d_1(s).
    \]
    Combining these expressions with the state-occupancy bound gives
    the constant
    $\gamma/(1-\gamma)+1=1/(1-\gamma)$.
\end{proof}

\begin{lemma}[No Policy-Only Converse]\label{lem:diff_polcies_lower}
    Under Assumptions \ref{ass:data_obs} and \ref{ass:Mark}, there is
    no strictly positive lower bound on
    $\|d_P^{\pi_2}-d_P^{\pi_1}\|$ that depends only on a positive
    divergence between $\pi_2$ and $\pi_1$.
\end{lemma}

\begin{proof}
    It is enough to construct one class of counterexamples.  Let the
    transition kernel be action independent:
    \[
        P(s'\mid s,a)=P(s'\mid s)
        \quad\text{for every }(s,a,s').
    \]
    Then, for every stationary policy $\pi$,
    \[
        K_{\pi,P}(s,s')
        =\sum_a\pi(a\mid s)P(s'\mid s,a)
        =P(s'\mid s)\sum_a\pi(a\mid s)
        =P(s'\mid s).
    \]
    Thus all policies induce the same state-transition kernel, say
    $K$.  Their discounted state occupancies are consequently
    identical:
    \[
        d_P^\pi
        =(1-\gamma)\mu_0(I-\gamma K)^{-1}
        \quad\text{for every }\pi.
    \]
    Choose any state $s_0$ with positive occupancy and two distinct
    action distributions at that state.  For example, with two
    actions, let
    $\pi_1(1\mid s_0)=1/4$ and
    $\pi_2(1\mid s_0)=3/4$, and make the policies equal elsewhere.
    Then their total-variation, sup-norm, and chi-square divergences
    at $s_0$ are strictly positive, whereas
    \[
        \|d_P^{\pi_2}-d_P^{\pi_1}\|_{\mathrm{TV}}=0.
    \]
    Therefore no inequality with a strictly positive right-hand side
    depending only on policy divergence can lower-bound the
    state-occupancy distance without additional assumptions linking
    actions to distinguishable transition behavior.
\end{proof}

Lemma~\ref{lem:diff_polcies_lower} is why the proof below does not try
to infer policy stability from occupancy-ratio stability.  Only the
valid one-sided bound in Lemma~\ref{lem:diff_polcies_upper} is used.%

\begin{lemma}[Value Continuity and Argmax Stability]
\label{lem:diff_pol_to_Q_fun}
    Suppose rewards are bounded by $\overline c_R$, $\gamma<1$, and
    define
    $d_\Pi(\pi_1,\pi_2)
      :=\sup_s\operatorname{TV}\{\pi_1(\cdot\mid s),
      \pi_2(\cdot\mid s)\}$.
    Then:
    \begin{enumerate}
        \item
        \[
            \|V_P^{\pi_2}-V_P^{\pi_1}\|_\infty
            \leq
            \frac{2\overline c_R}{(1-\gamma)^2}
            d_\Pi(\pi_1,\pi_2),
            \qquad
            \|Q_P^{\pi_2}-Q_P^{\pi_1}\|_\infty
            \leq
            \frac{2\gamma\overline c_R}{(1-\gamma)^2}
            d_\Pi(\pi_1,\pi_2).
        \]
        \item If $(\Pi,d_\Pi)$ is compact, $M:\Pi\to\mathbb R$ is
        continuous with unique maximizer $\pi^*$, and
        $\sup_{\pi\in\Pi}|M_n(\pi)-M(\pi)|\to0$, then every
        $\widehat\pi_n\in\arg\max_{\pi\in\Pi}M_n(\pi)$ satisfies
        $d_\Pi(\widehat\pi_n,\pi^*)\to0$.
    \end{enumerate}
\end{lemma}

\begin{proof}
    Define the Bellman operator
    \[
        (\mathbb T_P^\pi f)(s)
        :=
        \sum_a\pi(a\mid s)
        \left\{r_P(s,a)+\gamma
        \int f(s')P(\mathrm ds'\mid s,a)\right\}.
    \]
    It is a $\gamma$-contraction in sup norm:
    \[
        \|\mathbb T_P^\pi f-\mathbb T_P^\pi g\|_\infty
        \leq\gamma\|f-g\|_\infty.
    \]
    Since $V_P^{\pi_k}=\mathbb T_P^{\pi_k}V_P^{\pi_k}$,
    add and subtract
    $\mathbb T_P^{\pi_2}V_P^{\pi_1}$ to obtain
    \[
        \begin{aligned}
        \|V_P^{\pi_2}-V_P^{\pi_1}\|_\infty
        &\leq
        \|\mathbb T_P^{\pi_2}V_P^{\pi_2}
          -\mathbb T_P^{\pi_2}V_P^{\pi_1}\|_\infty\\
        &\quad+
        \|\mathbb T_P^{\pi_2}V_P^{\pi_1}
          -\mathbb T_P^{\pi_1}V_P^{\pi_1}\|_\infty\\
        &\leq
        \gamma\|V_P^{\pi_2}-V_P^{\pi_1}\|_\infty\\
        &\quad+
        \sup_s\left|
        \sum_a\{\pi_2(a\mid s)-\pi_1(a\mid s)\}
        Q_P^{\pi_1}(a,s)\right|.
        \end{aligned}
    \]
    For each $s$, the variational characterization of total
    variation gives
    \[
        \left|
        \sum_a\{\pi_2(a\mid s)-\pi_1(a\mid s)\}
        Q_P^{\pi_1}(a,s)\right|
        \leq
        2\operatorname{TV}\{\pi_2(\cdot\mid s),
        \pi_1(\cdot\mid s)\}
        \|Q_P^{\pi_1}\|_\infty.
    \]
    Bounded rewards imply, directly from the discounted series,
    \[
        \|Q_P^{\pi_1}\|_\infty
        \leq\sum_{t=0}^{\infty}\gamma^t\overline c_R
        =\frac{\overline c_R}{1-\gamma}.
    \]
    Moving the contraction term to the left therefore yields
    \[
        (1-\gamma)
        \|V_P^{\pi_2}-V_P^{\pi_1}\|_\infty
        \leq
        2d_\Pi(\pi_1,\pi_2)\|Q_P^{\pi_1}\|_\infty
        \leq
        \frac{2\overline c_R}{1-\gamma}
        d_\Pi(\pi_1,\pi_2).
    \]
    Dividing by $1-\gamma$ proves the first displayed bound.
    Moreover,
    \[
        Q_P^{\pi_k}(a,s)
        =r_P(s,a)+\gamma
        \int V_P^{\pi_k}(s')P(\mathrm ds'\mid s,a),
    \]
    so
    \[
        \|Q_P^{\pi_2}-Q_P^{\pi_1}\|_\infty
        \leq
        \gamma\|V_P^{\pi_2}-V_P^{\pi_1}\|_\infty,
    \]
    which proves the claimed $Q$-bound.

    For the second statement, suppose for contradiction that
    $d_\Pi(\widehat\pi_n,\pi^*)\not\to0$.  Then there are
    $\varepsilon>0$ and a subsequence, still indexed by $n$, such
    that
    $d_\Pi(\widehat\pi_n,\pi^*)\geq\varepsilon$ for every $n$.
    Compactness gives a further subsequence satisfying
    $\widehat\pi_n\to\bar\pi\in\Pi$, and continuity of the metric
    implies $d_\Pi(\bar\pi,\pi^*)\geq\varepsilon$.
    For any fixed $\pi\in\Pi$, optimality of $\widehat\pi_n$ gives
    $M_n(\widehat\pi_n)\geq M_n(\pi)$.  Hence
    \[
        M(\widehat\pi_n)
        \geq M(\pi)
        -2\sup_{\widetilde\pi\in\Pi}
        |M_n(\widetilde\pi)-M(\widetilde\pi)|.
    \]
    Letting $n\to\infty$ and using continuity of $M$ gives
    $M(\bar\pi)\geq M(\pi)$ for every $\pi\in\Pi$.  Thus
    $\bar\pi$ is a maximizer of $M$.  Uniqueness forces
    $\bar\pi=\pi^*$, contradicting
    $d_\Pi(\bar\pi,\pi^*)\geq\varepsilon$.
\end{proof}

\begin{lemma}[Continuity of the Normalized Episode Score]
\label{lem:episode_score_continuity}
    Under Assumptions~\ref{ass:Mark}, \ref{ass:stationary}, and
    \ref{ass:traj_score_reg}, there is a constant $C_T<\infty$ such
    that, for any two evaluated policies,
    \[
        \|\phi_T(\cdot;\pi_2)-\phi_T(\cdot;\pi_1)\|_{P_0,2}
        \leq C_T d_\Pi(\pi_1,\pi_2).
    \]
    In particular, the normalized trajectory score is
    $L_2(P_0)$-continuous in the policy.
\end{lemma}

\begin{proof}
    Let
    $\nu^\pi(\mathrm ds,\mathrm da)
      =d_{P_0}^\pi(\mathrm ds)\pi(\mathrm da\mid s)$.
    Lemma~\ref{lem:diff_polcies_upper} and the definition of
    $d_\Pi$ give
    \[
        \|\nu^{\pi_2}-\nu^{\pi_1}\|_{\mathrm{TV}}
        \leq\frac{d_\Pi(\pi_1,\pi_2)}{1-\gamma}.
    \]
    On the finite state--action space,
    $\omega^\pi(a,s)
      =\nu^\pi(s,a)/\{f_b(s)b(a\mid s)\}$.  Moreover, on a finite
    state--action space,
    $\max_{s,a}|\nu^{\pi_2}(s,a)-\nu^{\pi_1}(s,a)|
    \leq 2\|\nu^{\pi_2}-\nu^{\pi_1}\|_{\mathrm{TV}}$, and the
    denominator $f_b(s)b(a\mid s)$ is bounded below on the relevant
    support by the overlap condition.
    The lower overlap bounds therefore imply
    \[
        \|\omega^{\pi_2}-\omega^{\pi_1}\|_\infty
        \leq C_\omega d_\Pi(\pi_1,\pi_2),
        \qquad
        \sup_{\pi}\|\omega^\pi\|_\infty<\infty.
    \]
    Lemma~\ref{lem:diff_pol_to_Q_fun} similarly gives
    \[
        \|V^{\pi_2}-V^{\pi_1}\|_\infty
        +\|Q^{\pi_2}-Q^{\pi_1}\|_\infty
        +|\eta(\pi_2)-\eta(\pi_1)|
        \leq C_Vd_\Pi(\pi_1,\pi_2).
    \]

    For $k=1,2$, define
    $\varepsilon_t^{\pi_k}
      :=R_t+\gamma V^{\pi_k}(S_{t+1})
        -Q^{\pi_k}(A_t,S_t)$.
    Bounded rewards and values imply
    $\sup_{t,k}|\varepsilon_t^{\pi_k}|\leq C_\varepsilon$, and
    \[
        |\varepsilon_t^{\pi_2}
          -\varepsilon_t^{\pi_1}|
        \leq
        \gamma\|V^{\pi_2}-V^{\pi_1}\|_\infty
        +\|Q^{\pi_2}-Q^{\pi_1}\|_\infty.
    \]
    For each transition,
    \[
        \begin{aligned}
        |\omega^{\pi_2}\varepsilon_t^{\pi_2}
          -\omega^{\pi_1}\varepsilon_t^{\pi_1}|
        &\leq
        |\omega^{\pi_2}-\omega^{\pi_1}|
        |\varepsilon_t^{\pi_2}|\\
        &\quad+
        |\omega^{\pi_1}|
        |\varepsilon_t^{\pi_2}-\varepsilon_t^{\pi_1}|\\
        &\leq C d_\Pi(\pi_1,\pi_2).
        \end{aligned}
    \]
    Substituting this bound into the finite weighted sum defining
    $\phi_T$, and applying the corresponding bounds to its initial
    value and centering terms, yields the pointwise inequality
    \[
        |\phi_T(O_{0:T};\pi_2)-\phi_T(O_{0:T};\pi_1)|
        \leq C_Td_\Pi(\pi_1,\pi_2).
    \]
    Taking the $L_2(P_0)$ norm proves the result.
\end{proof}

\begin{lemma}[Uniform Fixed-Policy Expansion in the Finite MDP]
\label{lem:finite_mdp_uniform_expansion}
    Under Assumption~\ref{ass:regularity}, define
    \[
        r_{\epsilon,\pi}(s)
        :=\sum_a\pi(a\mid s)r_\epsilon(s,a),
        \qquad
        K_{\epsilon,\pi}(s,s')
        :=\sum_a\pi(a\mid s)p_\epsilon(s'\mid s,a),
    \]
    with
    $\dot r_\pi(s):=\sum_a\pi(a\mid s)\dot r(s,a)$ and
    $\dot K_\pi(s,s')
      :=\sum_a\pi(a\mid s)\dot p(s'\mid s,a)$,
    and $G_{0,\pi}:=(I-\gamma K_{0,\pi})^{-1}$.  Then
    \[
        V_{P_\epsilon}^{\pi}
        =V_{P_0}^{\pi}
        +\epsilon G_{0,\pi}
        \{\dot r_\pi+\gamma\dot K_\pi V_{P_0}^{\pi}\}
        +o(|\epsilon|)
    \]
    uniformly over $\pi\in\Pi$, where the remainder is in sup norm.
    Consequently,
    \begin{equation}\label{eq:finite_mdp_uniform_expansion}
        \sup_{\pi\in\Pi}
        \left|
        \Psi(P_\epsilon;\pi)-\Psi(P_0;\pi)
        -\epsilon\dot\Psi_{P_0}^{\pi}
        \right|
        =o(|\epsilon|),
    \end{equation}
    with
    \[
        \dot\Psi_{P_0}^{\pi}
        =
        \dot\mu^\top V_{P_0}^{\pi}
        +
        \mu_0^\top G_{0,\pi}
        \{\dot r_\pi+\gamma\dot K_\pi V_{P_0}^{\pi}\}.
    \]
    In the auxiliary experiment $\mathcal M^\dagger$,
    $\dot\Psi_{P_0}^{\pi}
      =\E_{P_0^\dagger}
      [D_{P_0}^{\dagger,\pi}\dot\ell^\dagger]$.
\end{lemma}

\begin{proof}
    Throughout this proof, the matrix norm is the maximum absolute
    row-sum norm, and the vector norm is the sup norm.  Because
    $K_{\epsilon,\pi}$ is stochastic,
    \[
        \sup_{\epsilon,\pi}
        \|(I-\gamma K_{\epsilon,\pi})^{-1}\|_\infty
        \leq
        \sum_{m=0}^{\infty}\gamma^m
        \|K_{\epsilon,\pi}^m\|_\infty
        =\frac1{1-\gamma}.
    \]
    Write
    $G_{\epsilon,\pi}:=(I-\gamma K_{\epsilon,\pi})^{-1}$,
    $\Delta r_{\epsilon,\pi}:=
      r_{\epsilon,\pi}-r_{0,\pi}$, and
    $\Delta K_{\epsilon,\pi}:=
      K_{\epsilon,\pi}-K_{0,\pi}$.
    The componentwise expansions in
    Assumption~\ref{ass:regularity} remain uniform after averaging
    with any policy, because all policy weights are nonnegative and
    sum to one.  More explicitly,
    \[
        \begin{aligned}
        \sup_{\pi}\|
        \Delta r_{\epsilon,\pi}-\epsilon\dot r_\pi
        \|_\infty
        &\leq
        \max_{s,a}|r_\epsilon(s,a)-r_0(s,a)
        -\epsilon\dot r(s,a)|
        =o(|\epsilon|),\\
        \sup_{\pi}\|
        \Delta K_{\epsilon,\pi}-\epsilon\dot K_\pi
        \|_\infty
        &\leq
        \max_{s,a}\sum_{s'}
        |p_\epsilon(s'\mid s,a)-p_0(s'\mid s,a)
        -\epsilon\dot p(s'\mid s,a)|
        =o(|\epsilon|).
        \end{aligned}
    \]
    In particular,
    $\sup_\pi\|\Delta r_{\epsilon,\pi}\|_\infty=O(|\epsilon|)$
    and
    $\sup_\pi\|\Delta K_{\epsilon,\pi}\|_\infty=O(|\epsilon|)$.

    The resolvent identity follows by multiplying out:
    \[
        \begin{aligned}
        G_{\epsilon,\pi}-G_{0,\pi}
        &=
        G_{\epsilon,\pi}
        \{(I-\gamma K_{0,\pi})
          -(I-\gamma K_{\epsilon,\pi})\}
        G_{0,\pi}\\
        &=
        \gamma G_{\epsilon,\pi}
        (K_{\epsilon,\pi}-K_{0,\pi})G_{0,\pi}
        =
        \gamma G_{\epsilon,\pi}
        \Delta K_{\epsilon,\pi}G_{0,\pi}.
        \end{aligned}
    \]
    Consequently,
    \[
        \sup_\pi\|G_{\epsilon,\pi}-G_{0,\pi}\|_\infty
        \leq
        \frac{\gamma}{(1-\gamma)^2}
        \sup_\pi\|\Delta K_{\epsilon,\pi}\|_\infty
        =O(|\epsilon|).
    \]

    We now derive the value expansion directly from the two Bellman
    equations.  Since
    $(I-\gamma K_{\epsilon,\pi})V_{P_\epsilon}^\pi
      =r_{\epsilon,\pi}$ and
    $(I-\gamma K_{0,\pi})V_{P_0}^\pi=r_{0,\pi}$,
    \[
        \begin{aligned}
        &(I-\gamma K_{\epsilon,\pi})
        (V_{P_\epsilon}^\pi-V_{P_0}^\pi)\\
        &\quad=
        r_{\epsilon,\pi}
        -(I-\gamma K_{\epsilon,\pi})V_{P_0}^\pi\\
        &\quad=
        \Delta r_{\epsilon,\pi}
        +\gamma\Delta K_{\epsilon,\pi}V_{P_0}^\pi.
        \end{aligned}
    \]
    Therefore the following identity is exact:
    \begin{equation}\label{eq:exact_value_perturbation}
        V_{P_\epsilon}^\pi-V_{P_0}^\pi
        =
        G_{\epsilon,\pi}
        \{\Delta r_{\epsilon,\pi}
        +\gamma\Delta K_{\epsilon,\pi}V_{P_0}^\pi\}.
    \end{equation}
    Bounded rewards and the resolvent bound imply
    $\sup_\pi\|V_{P_0}^\pi\|_\infty
      \leq\overline c_R/(1-\gamma)$.
    Define
    \[
        h_\pi:=\dot r_\pi+\gamma\dot K_\pi V_{P_0}^\pi.
    \]
    The primitive expansions and the uniform bound on $V_{P_0}^\pi$
    give
    \[
        \sup_\pi
        \|\Delta r_{\epsilon,\pi}
        +\gamma\Delta K_{\epsilon,\pi}V_{P_0}^\pi
        -\epsilon h_\pi\|_\infty=o(|\epsilon|),
        \qquad
        \sup_\pi\|h_\pi\|_\infty<\infty.
    \]
    Subtracting $\epsilon G_{0,\pi}h_\pi$ from
    \eqref{eq:exact_value_perturbation} yields
    \[
        \begin{aligned}
        &V_{P_\epsilon}^\pi-V_{P_0}^\pi
        -\epsilon G_{0,\pi}h_\pi\\
        &\quad=
        \epsilon(G_{\epsilon,\pi}-G_{0,\pi})h_\pi\\
        &\qquad+
        G_{\epsilon,\pi}
        \{\Delta r_{\epsilon,\pi}
        +\gamma\Delta K_{\epsilon,\pi}V_{P_0}^\pi
        -\epsilon h_\pi\}.
        \end{aligned}
    \]
    The first term is uniformly $O(\epsilon^2)$ and the second is
    uniformly $o(|\epsilon|)$.  This proves the asserted uniform
    expansion of $V^\pi$.

    Next expand the target value without suppressing the cross term:
    \[
        \begin{aligned}
        \Psi(P_\epsilon;\pi)-\Psi(P_0;\pi)
        &=(\mu_\epsilon-\mu_0)^\top V_{P_0}^\pi\\
        &\quad+\mu_0^\top
        (V_{P_\epsilon}^\pi-V_{P_0}^\pi)\\
        &\quad+(\mu_\epsilon-\mu_0)^\top
        (V_{P_\epsilon}^\pi-V_{P_0}^\pi).
        \end{aligned}
    \]
    The first two lines equal
    \[
        \epsilon\dot\mu^\top V_{P_0}^\pi
        +\epsilon\mu_0^\top G_{0,\pi}h_\pi
        +o(|\epsilon|)
    \]
    uniformly in $\pi$.  The final cross term is
    $O(\epsilon^2)$ uniformly, because both factors are
    $O(|\epsilon|)$.  This proves
    \eqref{eq:finite_mdp_uniform_expansion} and the displayed formula
    for $\dot\Psi_{P_0}^\pi$.

    It remains to identify the derivative in
    $\mathcal M^\dagger$.  Its product-factorized density can be
    written as
    \[
        p^\dagger(s_0,s,a,r,s')
        =\mu(s_0)g(s)b(a\mid s)c(r,s'\mid a,s),
    \]
    so a regular score decomposes as
    \[
        \dot\ell^\dagger
        =\dot\ell_\mu(S_0)+\dot\ell_g(S)
        +\dot\ell_b(A\mid S)
        +\dot\ell_c(R,S'\mid A,S),
    \]
    where $g$ is the transition-sampling state marginal.  The
    conditional score components satisfy
    $\E(\dot\ell_\mu)=0$,
    $\E(\dot\ell_g)=0$,
    $\E(\dot\ell_b\mid S)=0$, and
    $\E(\dot\ell_c\mid S,A)=0$.
    Differentiating the primitive laws gives
    \[
        \begin{aligned}
        \dot\mu^\top V_{P_0}^\pi
        &=\E\!\left[
        \{V_{P_0}^\pi(S_0)-\Psi(P_0;\pi)\}
        \dot\ell_\mu(S_0)\right],\\
        \dot r(s,a)
        &=\E\!\left[
        \{R-r_0(s,a)\}\dot\ell_c(R,S'\mid a,s)
        \mid S=s,A=a\right],\\
        \sum_{s'}\dot p(s'\mid s,a)V_{P_0}^\pi(s')
        &=\E\!\left[
        \{V_{P_0}^\pi(S')
        -P_0V_{P_0}^\pi(s,a)\}
        \dot\ell_c(R,S'\mid a,s)
        \mid S=s,A=a\right].
        \end{aligned}
    \]
    The centering terms in the last two lines may be omitted because
    the conditional score has mean zero.  By the Bellman equation,
    $Q_{P_0}^\pi(s,a)
      =r_0(s,a)+\gamma P_0V_{P_0}^\pi(s,a)$, and hence
    \[
        h_\pi(s)
        =
        \sum_a\pi(a\mid s)
        \E\!\left[
        \{R+\gamma V_{P_0}^\pi(S')
        -Q_{P_0}^\pi(A,S)\}\dot\ell_c
        \mid S=s,A=a\right].
    \]

    The discounted occupancy identity follows from the Neumann
    series:
    \[
        \begin{aligned}
        \mu_0^\top G_{0,\pi}h
        &=\sum_{t=0}^{\infty}\gamma^t
          \mu_0^\top K_{0,\pi}^t h\\
        &=\frac1{1-\gamma}
          \E_{S\sim d_{P_0}^\pi}\{h(S)\}.
        \end{aligned}
    \]
    Changing measure from the discounted target occupancy to the
    behavior transition marginal uses
    \[
        \omega^\pi(a,s)
        =\frac{d_{P_0}^\pi(s)\pi(a\mid s)}
        {f_b(s)b(a\mid s)}.
    \]
    Thus the second term of $\dot\Psi_{P_0}^\pi$ is
    \[
        \frac1{1-\gamma}
        \E_{P_{0,b}^{\mathrm{tr}}}\!\left[
        \omega^\pi(A,S)
        \{R+\gamma V_{P_0}^\pi(S')
        -Q_{P_0}^\pi(A,S)\}
        \dot\ell_c(R,S'\mid A,S)\right].
    \]
    Let
    $\varepsilon^\pi
      :=R+\gamma V_{P_0}^\pi(S')-Q_{P_0}^\pi(A,S)$.
    The Bellman equation gives
    $\E(\varepsilon^\pi\mid S,A)=0$.  Consequently,
    \[
        \E\{\omega^\pi\varepsilon^\pi\dot\ell_g(S)\}=0,
        \qquad
        \E\{\omega^\pi\varepsilon^\pi
        \dot\ell_b(A\mid S)\}=0.
    \]
    The initial-state component is centered and, in the product
    auxiliary experiment, is orthogonal to the transition-sampling
    score components.  Combining these orthogonality relations with
    the preceding derivative calculation gives
    \[
        \dot\Psi_{P_0}^{\pi}
        =\E_{P_0^\dagger}\!\left[
        D_{P_0}^{\dagger,\pi}(O^\dagger)
        \dot\ell^\dagger(O^\dagger)\right].
    \]
    Each summand of $D_{P_0}^{\dagger,\pi}$ is mean zero and belongs
    to the corresponding component of the full nonparametric product
    tangent space.  Therefore the displayed gradient belongs to
    $\mathcal T^\dagger$ and is the canonical gradient.

    Finally, in the finite MDP,
    $\pi\mapsto K_{0,\pi}$ and $\pi\mapsto r_{0,\pi}$ are continuous.
    The resolvent identity then makes
    $\pi\mapsto V_{P_0}^\pi$ and
    $\pi\mapsto Q_{P_0}^\pi$ continuous in sup norm.
    Lemma~\ref{lem:diff_polcies_upper}, the finite-state lower overlap
    bounds, and
    $\omega^\pi(a,s)
      =d^\pi(s)\pi(a\mid s)/\{f_b(s)b(a\mid s)\}$ make
    $\pi\mapsto\omega^\pi$ continuous in $L_1$ and hence in $L_2$
    on this finite, uniformly bounded space.  All factors in
    $D_{P_0}^{\dagger,\pi}$ are uniformly bounded, so
    $\pi\mapsto D_{P_0}^{\dagger,\pi}$ is continuous in
    $L_2(P_0^\dagger)$.
\end{proof}

\begin{lemma}[Exact Fixed-Policy Double-Robustness Identity]
\label{lem:exact_dr_identity}
    Fix a policy $\pi$.  Let $\bar Q$ and $\bar\omega$ be arbitrary
    square-integrable functions and set
    $\bar V(s)=\sum_a\pi(a\mid s)\bar Q(a,s)$.  Under Assumptions
    \ref{ass:Mark} and \ref{ass:stationary},
    \[
        \begin{aligned}
        &\E_{P_0^\dagger}\!
        \left[\EstFunPoi_{\eta(\pi)}
        (O^\dagger;\bar Q,\bar\omega,\bar V,\pi)\right]-\eta(\pi)\\
        &\quad =
        \frac{1}{1-\gamma}
        \E_{P_{0,b}^{\mathrm{tr}}}\!\left[
        \{\bar\omega-\omega^\pi\}(A,S)
        \{\gamma(\bar V-V^\pi)(S')
        -(\bar Q-Q^\pi)(A,S)\}\right].
        \end{aligned}
    \]
    The normalized trajectory score in
    \eqref{eq:debiased_terms_for_fixed_policy} has the same
    expectation.  Thus the score is unbiased if either
    $\bar Q=Q^\pi$ or $\bar\omega=\omega^\pi$.
\end{lemma}

\begin{proof}
    Put
    \[
        \delta Q:=\bar Q-Q^\pi,
        \qquad
        \delta V:=\bar V-V^\pi.
    \]
    Because
    $\E\{R+\gamma V^\pi(S')-Q^\pi(A,S)\mid S,A\}=0$,
    multiplying the true Bellman residual by any square-integrable
    function of $(S,A)$ still gives mean zero.  Therefore, after
    subtracting the score at the true nuisances,
    \[
        \begin{aligned}
        &\E_{P_0^\dagger}
        [\EstFunPoi_{\eta(\pi)}
        (O^\dagger;\bar Q,\bar\omega,\bar V,\pi)]
        -\eta(\pi)\\
        &\quad=
        \frac1{1-\gamma}
        \E_{P_{0,b}^{\mathrm{tr}}}\!
        \left[\bar\omega(A,S)
        \{\gamma\delta V(S')-\delta Q(A,S)\}\right]
        +\E_{\mu_0}\{\delta V(S_0)\}.
        \end{aligned}
    \]
    We now prove the cancellation involving the true occupancy ratio.
    By the definition of $\omega^\pi$,
    \[
        \begin{aligned}
        \E_{P_{0,b}^{\mathrm{tr}}}
        \{\omega^\pi(A,S)\delta Q(A,S)\}
        &=
        \E_{S\sim d_{P_0}^\pi}
        \sum_a\pi(a\mid S)\delta Q(a,S)\\
        &=\E_{S\sim d_{P_0}^\pi}\{\delta V(S)\},
        \end{aligned}
    \]
    while
    \[
        \E_{P_{0,b}^{\mathrm{tr}}}
        \{\omega^\pi(A,S)\delta V(S')\}
        =
        \E_{S\sim d_{P_0}^\pi K_{\pi,P_0}}
        \{\delta V(S)\}.
    \]
    The discounted occupancy flow equation is
    \[
        d_{P_0}^\pi
        =(1-\gamma)\mu_0
        +\gamma d_{P_0}^\pi K_{\pi,P_0}.
    \]
    Integrating $\delta V$ against both sides gives the exact
    identity
    \[
        \E_{P_{0,b}^{\mathrm{tr}}}
        \left[
        \omega^\pi(A,S)
        \{\delta Q(A,S)-\gamma\delta V(S')\}\right]
        =(1-\gamma)\E_{\mu_0}\{\delta V(S_0)\}.
    \]
    Equivalently,
    \[
        \frac1{1-\gamma}
        \E\!\left[
        \omega^\pi\{\gamma\delta V(S')-\delta Q\}\right]
        +\E_{\mu_0}(\delta V)=0.
    \]
    Subtract this zero quantity from the previous score expansion.
    The result is
    \[
        \frac1{1-\gamma}
        \E\!\left[
        \{\bar\omega-\omega^\pi\}
        \{\gamma\delta V(S')-\delta Q(A,S)\}\right],
    \]
    which is the claimed product remainder.

    Under Assumption~\ref{ass:stationary}, every transition tuple
    $(S_t,A_t,R_t,S_{t+1})$ has marginal law
    $P_{0,b}^{\mathrm{tr}}$.  If $c$ denotes the expectation of the
    unscaled transition term, then the transition part of the
    normalized trajectory score has expectation
    \[
        \frac1{1-\gamma^{T+1}}
        \sum_{t=0}^T\gamma^t c
        =
        \frac1{1-\gamma^{T+1}}
        \frac{1-\gamma^{T+1}}{1-\gamma}c
        =\frac{c}{1-\gamma}.
    \]
    Its initial-state term has the same expectation as above.
    Hence the trajectory and point-level scores have the same
    expectation and the same double-robustness identity.
\end{proof}

\subsection{Envelope proof of Theorems
\ref{thm:suff_cond_for_EIF_exp} and \ref{thm:nec_cond_for_EIF}}

Fix a regular parametric submodel of $\mathcal M^\dagger$ through
$P_0^\dagger$ with score $\dot\ell^\dagger$, and write
\[
    L_\pi^\dagger(\dot\ell^\dagger)
    :=\E_{P_0^\dagger}\!\left[
    D_{P_0}^{\dagger,\pi}(O^\dagger)
    \dot\ell^\dagger(O^\dagger)\right].
\]
Write $m_\epsilon(\pi):=\Psi(P_\epsilon;\pi)$,
$v_\epsilon:=\max_{\pi\in\Pi}m_\epsilon(\pi)$, and
$\mathcal S:=\Pi^*(P_0)$.  Lemma
\ref{lem:finite_mdp_uniform_expansion} gives
\begin{equation}\label{eq:uniform_envelope_remainder}
    m_\epsilon(\pi)
    =m_0(\pi)+\epsilon
    L_\pi^\dagger(\dot\ell^\dagger)+r_\epsilon(\pi),
    \qquad
    \sup_{\pi\in\Pi}
    \frac{|r_\epsilon(\pi)|}{|\epsilon|}\longrightarrow0.
\end{equation}
Because $\Pi$ is a closed subset of the finite product of action
simplexes, it is compact; Lemma~\ref{lem:diff_pol_to_Q_fun} gives
continuity of each $m_\epsilon$.  Hence all maxima below are attained.
We now prove the one-sided envelope formulas
\begin{equation}\label{eq:optimal_value_envelope}
    \begin{aligned}
        \lim_{\epsilon\downarrow0}
        \frac{\Psi^*(P_\epsilon)-\Psi^*(P_0)}{\epsilon}
        &=
        \max_{\pi\in\Pi^*(P_0)}
        L_\pi^\dagger(\dot\ell^\dagger),\\
        \lim_{\epsilon\uparrow0}
        \frac{\Psi^*(P_\epsilon)-\Psi^*(P_0)}{\epsilon}
        &=
        \min_{\pi\in\Pi^*(P_0)}
        L_\pi^\dagger(\dot\ell^\dagger).
    \end{aligned}
\end{equation}
Fix any $\pi\in\mathcal S$.  For $\epsilon>0$,
$v_\epsilon\geq m_\epsilon(\pi)$ and $v_0=m_0(\pi)$, so
\[
    \frac{v_\epsilon-v_0}{\epsilon}
    \geq
    L_\pi^\dagger(\dot\ell^\dagger)
    +\frac{r_\epsilon(\pi)}{\epsilon}.
\]
Taking the lower limit and maximizing over $\pi\in\mathcal S$ gives
\begin{equation}\label{eq:envelope_lower}
    \liminf_{\epsilon\downarrow0}
    \frac{v_\epsilon-v_0}{\epsilon}
    \geq
    \max_{\pi\in\mathcal S}
    L_\pi^\dagger(\dot\ell^\dagger).
\end{equation}

For the reverse inequality, choose
$\pi_\epsilon\in\arg\max_\pi m_\epsilon(\pi)$.  Since
$m_0(\pi_\epsilon)\leq v_0$,
\[
    \begin{aligned}
    v_\epsilon-v_0
    &=m_\epsilon(\pi_\epsilon)-v_0\\
    &\leq m_\epsilon(\pi_\epsilon)-m_0(\pi_\epsilon)\\
    &=\epsilon L_{\pi_\epsilon}^\dagger
      (\dot\ell^\dagger)+r_\epsilon(\pi_\epsilon).
    \end{aligned}
\]
Uniform convergence of $m_\epsilon$ to $m_0$ implies
\[
    |v_\epsilon-v_0|
    \leq\sup_{\pi\in\Pi}|m_\epsilon(\pi)-m_0(\pi)|
    \longrightarrow0.
\]
It also implies
\[
    \begin{aligned}
    0\leq v_0-m_0(\pi_\epsilon)
    &\leq |v_0-v_\epsilon|
      +|m_\epsilon(\pi_\epsilon)-m_0(\pi_\epsilon)|\\
    &\leq |v_0-v_\epsilon|
      +\sup_{\pi\in\Pi}|m_\epsilon(\pi)-m_0(\pi)|
      \longrightarrow0.
    \end{aligned}
\]
Thus every limit point of $\pi_\epsilon$ belongs to $\mathcal S$.
Lemma~\ref{lem:finite_mdp_uniform_expansion} also gives continuity of
$\pi\mapsto L_\pi^\dagger(\dot\ell^\dagger)$.  Dividing the preceding
upper bound by $\epsilon>0$ and taking the upper limit yields
\[
    \limsup_{\epsilon\downarrow0}
    \frac{v_\epsilon-v_0}{\epsilon}
    \leq
    \max_{\pi\in\mathcal S}
    L_\pi^\dagger(\dot\ell^\dagger).
\]
Together with \eqref{eq:envelope_lower}, this proves the first line of
\eqref{eq:optimal_value_envelope}.

For the left derivative, write $\epsilon=-t$ with $t\downarrow0$.
The expansion along $P_{-t}$ has linear coefficient
$-L_\pi^\dagger(\dot\ell^\dagger)$.  Applying the right-derivative
result to $t\mapsto P_{-t}$ gives
\[
    \lim_{t\downarrow0}\frac{v_{-t}-v_0}{t}
    =\max_{\pi\in\mathcal S}
    \{-L_\pi^\dagger(\dot\ell^\dagger)\}
    =-\min_{\pi\in\mathcal S}
    L_\pi^\dagger(\dot\ell^\dagger).
\]
Multiplication by $-1$ proves the second line of
\eqref{eq:optimal_value_envelope}.

\begin{proof}[Proof of Theorem \ref{thm:suff_cond_for_EIF_exp}]
    Under Assumption~\ref{ass:deterministic_policy},
    $\Pi^*(P_0)=\{\pi_0^*\}$, so the maximum and minimum in
    \eqref{eq:optimal_value_envelope} coincide for every score:
    \[
        \left.\frac{\mathrm d}{\mathrm d\epsilon}
        \Psi^*(P_\epsilon)\right|_{\epsilon=0}
        =
        \E_{P_0^\dagger}\!\left[
        D_{P_0}^{\dagger,\pi_0^*}(O^\dagger)
        \dot\ell^\dagger(O^\dagger)\right].
    \]
    Thus the right and left derivatives agree along every regular
    submodel.  The derivative is the continuous linear functional
    \[
        \dot\ell^\dagger\longmapsto
        \langle D_{P_0}^{\dagger,\pi_0^*},
        \dot\ell^\dagger\rangle_{P_0^\dagger},
    \]
    whose absolute value is at most
    $\|D_{P_0}^{\dagger,\pi_0^*}\|_2
      \|\dot\ell^\dagger\|_2$ by Cauchy--Schwarz.
    Therefore $\Psi^*$ is pathwise differentiable with gradient
    $D_{P_0}^{\dagger,\pi_0^*}$.  Lemma
    \ref{lem:finite_mdp_uniform_expansion} shows that this gradient
    belongs to the tangent space, so it is the canonical gradient and
    hence the EIF of the optimized value in
    $\mathcal M^\dagger$.
\end{proof}

\begin{proof}[Proof of Corollary
\ref{cor:equivalent_tied_optima}]
    If every optimizer has the same gradient
    $D_{P_0}^{\dagger,*}$, the maximum and minimum in
    \eqref{eq:optimal_value_envelope} both equal
    \[
        \E_{P_0^\dagger}\!\left[
        D_{P_0}^{\dagger,*}(O^\dagger)
        \dot\ell^\dagger(O^\dagger)\right]
    \]
    for every regular score $\dot\ell^\dagger$.  Hence the right and
    left derivatives agree and are represented by a single
    continuous linear functional of the score.  Because the common
    gradient belongs to the tangent space of
    $\mathcal M^\dagger$, it is the canonical gradient and therefore
    the EIF of $\Psi^*$ in that experiment.
\end{proof}

\begin{proof}[Proof of Theorem \ref{thm:nec_cond_for_EIF}]
    Under Assumption~\ref{ass:unres_rules}, choose
    $\pi_1,\pi_2\in\Pi^*(P_0)$ with distinct fixed-policy EIFs.
    Let
    $g=D_{P_0}^{\dagger,\pi_1}
       -D_{P_0}^{\dagger,\pi_2}$.
    Both canonical gradients belong to the tangent space
    $\mathcal T^\dagger$, so
    $g\in\mathcal T^\dagger$ and
    $\|g\|_{P_0^\dagger,2}>0$.  Because $\mathcal T^\dagger$ is
    the $L_2(P_0^\dagger)$ closure of regular-submodel scores, some
    regular score $h$ must satisfy
    $\E_{P_0^\dagger}(gh)\neq0$.  Indeed, if
    $\E(gh)=0$ for every regular score, continuity of the inner
    product would make $g$ orthogonal to their closure
    $\mathcal T^\dagger$.  Since $g\in\mathcal T^\dagger$, this would
    imply $\E(g^2)=0$, contradicting $\|g\|_2>0$.

    For this score $h$,
    \[
        \begin{aligned}
        L_{\pi_1}^\dagger(h)-L_{\pi_2}^\dagger(h)
        &=
        \E_{P_0^\dagger}\!\left[
        \{D_{P_0}^{\dagger,\pi_1}
        -D_{P_0}^{\dagger,\pi_2}\}h\right]\\
        &=\E_{P_0^\dagger}(gh)\neq0.
        \end{aligned}
    \]
    Thus at least two values in
    $\{L_\pi^\dagger(h):\pi\in\Pi^*(P_0)\}$ differ, so the maximum is
    strictly larger than the minimum.  By
    \eqref{eq:optimal_value_envelope}, the right and left derivatives
    along this regular submodel are unequal.  An influence function
    would have to represent one two-sided linear derivative
    $\dot\ell^\dagger\mapsto\E(D\dot\ell^\dagger)$ along every
    regular submodel.  Such a representation is impossible here;
    hence no influence function exists.
\end{proof}

\section{Technical Proofs in Section \ref{sec:inference}}

\begin{lemma}[Greedy-Policy Regret under the Margin Condition]
\label{lem:greedy_regret_margin}
    Let $\widehat Q^*$ be any estimate of $Q_{P_0}^*$, let
    $\widehat\pi$ be the deterministic policy obtained by the fixed
    tie-breaking greedy rule applied to $\widehat Q^*$, and put
    $\delta:=\|\widehat Q^*-Q_{P_0}^*\|_\infty$.  Under
    Assumptions~\ref{ass:traj_score_reg},
    \ref{ass:bellman_complete}, and \ref{ass:margin_con}, there is a
    finite constant
    $C_{\mathrm{reg}}$ such that, for every $\delta\geq0$,
    \[
        0\leq\eta^*-\eta(\widehat\pi)
        \leq C_{\mathrm{reg}}\delta^{1+\alpha},
    \]
    where $C_{\mathrm{reg}}$ does not depend on
    $\widehat Q^*$ or $\widehat\pi$.
\end{lemma}

\begin{proof}
    Write
    $\Delta_{\min}(s):=
      \min_{a\in\mathcal A_{\mathrm{sub-opt}}(s)}
      \Delta^*(a,s)$, with the maintained convention
    $\Delta_{\min}(s)=+\infty$ if the set is empty.
    Fix a state $s$, let $a^*(s)$ be any optimal action, and let
    $\widehat a(s)$ be the action selected by the greedy rule.
    Greediness gives
    $\widehat Q^*(\widehat a(s),s)
      \geq\widehat Q^*(a^*(s),s)$.  Therefore
    \[
        \begin{aligned}
        \Delta^*(\widehat a(s),s)
        &=Q_{P_0}^*(a^*(s),s)
          -Q_{P_0}^*(\widehat a(s),s)\\
        &=
        \{Q_{P_0}^*(a^*(s),s)
          -\widehat Q^*(a^*(s),s)\}\\
        &\quad+
        \{\widehat Q^*(a^*(s),s)
          -\widehat Q^*(\widehat a(s),s)\}\\
        &\quad+
        \{\widehat Q^*(\widehat a(s),s)
          -Q_{P_0}^*(\widehat a(s),s)\}\\
        &\leq\delta+0+\delta=2\delta.
        \end{aligned}
    \]
    If $\widehat a(s)$ is suboptimal, then the smallest positive gap
    at $s$ is no larger than its gap.  Thus
    \[
        \Delta^*(\widehat a(s),s)
        \leq
        2\delta\,
        \mathds{1}\!\left\{
        0<
        \min_{a\in\mathcal A_{\mathrm{sub-opt}}(s)}
        \Delta^*(a,s)
        \leq2\delta\right\}.
    \]
    Under Assumption~\ref{ass:bellman_complete}, choose an optimal
    greedy policy $\pi^*$ for which
    $V^{\pi^*}=V^*$ and $Q^{\pi^*}=Q^*$.  Applying the
    performance-difference identity with reference policy $\pi^*$
    and occupancy of $\widehat\pi$ gives
    \[
        \eta^*-\eta(\widehat\pi)
        =
        \frac1{1-\gamma}
        \E_{S\sim d_{P_0}^{\widehat\pi}}
        \{\Delta^*(\widehat a(S),S)\}.
    \]
    The preceding pointwise inequality and the uniform margin
    condition yield
    \[
        \begin{aligned}
        0\leq\eta^*-\eta(\widehat\pi)
        &\leq
        \frac{2\delta}{1-\gamma}
        \Pr_{S\sim d_{P_0}^{\widehat\pi}}
        \left\{0<\Delta_{\min}(S)\leq2\delta\right\}\\
        &\leq
        \frac{2\delta}{1-\gamma}C(2\delta)^\alpha
        =\frac{2^{1+\alpha}C}{1-\gamma}
        \delta^{1+\alpha}.
        \end{aligned}
    \]
    This proves the claim when $2\delta$ lies in the range where the
    margin condition holds.  To extend it globally, choose
    $\bar\delta>0$ so that the margin condition applies whenever
    $2\delta\leq2\bar\delta$.  Bounded rewards give
    \[
        0\leq\eta^*-\eta(\widehat\pi)
        \leq\frac{2\overline c_R}{1-\gamma}.
    \]
    If $\delta>\bar\delta$, the last bound is at most
    \[
        \frac{2\overline c_R}
        {(1-\gamma)\bar\delta^{1+\alpha}}
        \delta^{1+\alpha}.
    \]
    Taking $C_{\mathrm{reg}}$ to be the larger of this constant and
    $2^{1+\alpha}C/(1-\gamma)$ proves the stated bound for all
    $\delta\geq0$.
\end{proof}

\subsection{Proof of Theorem \ref{thm:CLT_for_R1_CLB} and
Corollary \ref{cor:lower_bound_method}}\label{sec:lower_bound_proofs_new}

Write $n_N=N-\ell_N$ and let
$\widehat w_j=\widehat\sigma_{\tau(j-1)}^{-1}$ and
$w_j^0=\widetilde\sigma_{\tau(j-1)}^{-1}$.  With
$\widehat\phi_j$ and $M_j$ as defined above, put
\[
    B_j:=\E(\widehat\phi_j\mid\mathcal F_{j-1})
    -\eta(\widehat\pi_j).
\]
Then
\[
    \begin{aligned}
    R_{\mathrm{CLB},1N}
    &=
    \left(\sum_{j=\ell_N+1}^N\widehat w_j\right)^{-1}
    \sum_{j=\ell_N+1}^N
    \widehat w_j\{\widehat\phi_j-\eta(\widehat\pi_j)\}\\
    &=
    \left(\sum_{j=\ell_N+1}^N\widehat w_j\right)^{-1}
    \sum_{j=\ell_N+1}^N\widehat w_j(M_j+B_j).
    \end{aligned}
\]
Conditional on $\mathcal F_{j-1}$, the policy and nuisance estimates
are fixed.  Lemma~\ref{lem:exact_dr_identity} gives
\[
    B_j=
    \frac1{1-\gamma}
    \E\!\left[
    \left.
    (\widehat\omega_j-\omega^{\widehat\pi_j})
    \{\gamma(\widehat V_j-V^{\widehat\pi_j})(S')
    -(\widehat Q_j-Q^{\widehat\pi_j})(A,S)\}
    \right|\mathcal F_{j-1}\right].
\]
To verify the needed norm comparison, write
$e_Q=\widehat Q_j-Q^{\widehat\pi_j}$ and
$e_V(s)=\sum_a\widehat\pi_j(a\mid s)e_Q(a,s)$.
Jensen's inequality and behavior overlap give
\[
    \begin{aligned}
    \E_{S\sim f_b}\{e_V^2(S)\}
    &\leq
    \E_{S\sim f_b}
    \sum_a\widehat\pi_j(a\mid S)e_Q^2(a,S)\\
    &\leq
    \underline b^{-1}
    \E_{S\sim f_b,A\sim b(\cdot\mid S)}
    \{e_Q^2(A,S)\},
    \end{aligned}
\]
where $\underline b>0$ is the overlap lower bound.  Stationarity
gives $S'\sim f_b$, so the triangle inequality implies
\[
    \|\gamma(\widehat V_j-V^{\widehat\pi_j})(S')
      -(\widehat Q_j-Q^{\widehat\pi_j})(A,S)\|_{P_0,2}
    \leq(1+\gamma\underline b^{-1/2})
    \|\widehat Q_j-Q^{\widehat\pi_j}\|_{P_0,2}.
\]
Thus Cauchy--Schwarz yields
\[
    |B_j|
    \leq C
    \|\widehat Q_j-Q^{\widehat\pi_j}\|_{P_0,2}
    \|\widehat\omega_j-\omega^{\widehat\pi_j}\|_{P_0,2}.
\]
On the scale-bounded event in Assumption~\ref{ass:non_zero_var},
$\widehat w_j\leq\underline\sigma^{-1}$.  Hence Assumption
\ref{ass:nui_rates} gives
\[
    \begin{aligned}
    \left|
    \frac1{\sqrt{n_N}}
    \sum_{j=\ell_N+1}^N\widehat w_jB_j
    \right|
    &\leq
    \frac{C}{\underline\sigma\sqrt{n_N}}
    \sum_{j=\ell_N+1}^N
    \|\widehat Q_j-Q^{\widehat\pi_j}\|_{P_0,2}
    \|\widehat\omega_j-\omega^{\widehat\pi_j}\|_{P_0,2}\\
    &=o_{P_0}(1).
    \end{aligned}
\]

Put
$e_j:=\widehat\sigma_{\tau(j-1)}
/\widetilde\sigma_{\tau(j-1)}-1$.  On the same event,
\[
    |\widehat w_j-w_j^0|
    \widetilde\sigma_{\tau(j-1)}
    =
    \left|\frac{e_j}{1+e_j}\right|
    \leq C|e_j|.
\]
Assumption~\ref{ass:est_var} and Cauchy--Schwarz therefore imply
\[
    \begin{aligned}
    \frac1{n_N}\sum_j
    (\widehat w_j-w_j^0)^2
    \widetilde\sigma_{\tau(j-1)}^2
    &=o_{P_0}(1),\\
    \frac1{n_N}\sum_j|\widehat w_j-w_j^0|
    &=o_{P_0}(1).
    \end{aligned}
\]
Because
$\overline\sigma^{-1}\leq n_N^{-1}\sum_jw_j^0
\leq\underline\sigma^{-1}$, it follows that
\[
    \frac{\sum_j\widehat w_j}{\sum_jw_j^0}\overset{P_0}{\to}1,
    \qquad
    \frac{\widehat\sigma_{R_{1N}}}{\sigma_{R_{1N}}}
    =\frac{\sum_jw_j^0}{\sum_j\widehat w_j}
    \overset{P_0}{\to}1.
\]

Since $\widehat w_j-w_j^0$ is predictable,
\[
    D_N:=\frac1{\sqrt{n_N}}
    \sum_{j=\ell_N+1}^N(\widehat w_j-w_j^0)M_j
\]
is a martingale with predictable quadratic variation
\[
    \begin{aligned}
    \langle D_N\rangle
    &=
    \frac1{n_N}\sum_{j=\ell_N+1}^N
    (\widehat w_j-w_j^0)^2
    \E(M_j^2\mid\mathcal F_{j-1})\\
    &=
    \frac1{n_N}\sum_{j=\ell_N+1}^N
    (\widehat w_j-w_j^0)^2
    \widetilde\sigma_{\tau(j-1)}^2
    =o_{P_0}(1).
    \end{aligned}
\]
More explicitly, Lenglart's inequality gives, for any
$\varepsilon,\delta>0$,
\[
    \Pr(|D_N|>\varepsilon)
    \leq\frac{\delta}{\varepsilon^2}
    +\Pr\{\langle D_N\rangle>\delta\}.
\]
First let $N\to\infty$ and then $\delta\downarrow0$ to obtain
$D_N=o_{P_0}(1)$.

Using
$\sigma_{R_{1N}}=\sqrt{n_N}/\sum_jw_j^0$, we now obtain
\[
    \begin{aligned}
    \sigma_{R_{1N}}^{-1}R_{\mathrm{CLB},1N}
    &=
    \frac{\sum_jw_j^0}{\sum_j\widehat w_j}
    \frac1{\sqrt{n_N}}
    \sum_j\widehat w_j(M_j+B_j)\\
    &=
    \frac1{\sqrt{n_N}}\sum_{j=\ell_N+1}^Nw_j^0M_j
    +o_{P_0}(1)\\
    &=
    \frac{1}{\sqrt{n_N}}
    \sum_{j=\ell_N+1}^N
    \frac{M_j}{\widetilde\sigma_{\tau(j-1)}}
    +o_{P_0}(1).
    \end{aligned}
\]
Define
$X_{N,j}:=M_j/
\{\sqrt{n_N}\widetilde\sigma_{\tau(j-1)}\}$.
Then
\[
    \E(X_{N,j}\mid\mathcal F_{j-1})=0,
    \qquad
    \sum_{j=\ell_N+1}^N
    \E(X_{N,j}^2\mid\mathcal F_{j-1})=1.
\]
Assumption~\ref{ass:Lindeberg} is exactly the conditional Lindeberg
condition for this array.  The martingale CLT gives
$\sum_jX_{N,j}\rightsquigarrow\mathcal N(0,1)$
\citep{hall2014martingale}.
The preceding expansion proves Theorem~\ref{thm:CLT_for_R1_CLB};
the estimated-scale conclusion follows from Slutsky's theorem.

Finally,
$\overline\eta_w(\widehat\pi)\leq\eta^*$ because it is a convex
combination of feasible policy values.  On the event
$R_{\mathrm{CLB},1N}
\leq z_{1-\alpha}\widehat\sigma_{R_{1N}}$,
\[
    \eta^*
    \geq\overline\eta_w(\widehat\pi)
    =\widehat\eta_{\mathrm{NSAVE}}-R_{\mathrm{CLB},1N}
    \geq
    \widehat\eta_{\mathrm{NSAVE}}
    -z_{1-\alpha}\widehat\sigma_{R_{1N}}.
\]
The studentized CLT makes the probability of the first event converge
to $1-\alpha$.  This proves
Corollary~\ref{cor:lower_bound_method}.

\subsection{Proof of Theorem \ref{thm:CLT_for_R1_TCI}}

Decompose
\[
    R_{\mathrm{TCI},1N}
    =
    R_{\mathrm{CLB},1N}
    +
    \left(\sum_j\widehat w_j\right)^{-1}
    \sum_j\widehat w_j
    \{\eta(\widehat\pi_j)-\eta(\widehat\pi_{\mathrm{fin}})\}.
\]
Let
$\delta_j=\|\widehat Q_j^*-Q_{P_0}^*\|_\infty$ and
$r=\kappa_*(1+\alpha)>1/2$.
By Assumption~\ref{ass:est_pol}, for every $\varepsilon>0$ there is
$C_\varepsilon<\infty$ such that, with probability at least
$1-\varepsilon$ for all sufficiently large $N$,
\[
    \delta_j\leq
    C_\varepsilon(j-\ell_N)^{-\kappa_*}
    \quad\text{simultaneously for }\ell_N<j\leq N.
\]
On this event, Lemma~\ref{lem:greedy_regret_margin} gives
\[
    0\leq\eta^*-\eta(\widehat\pi_j)
    \leq C_\varepsilon'(j-\ell_N)^{-r}
\]
simultaneously over the same indices.  The separate final-fit rate in
Assumption~\ref{ass:est_pol} and the same lemma also give
\[
    \eta^*-\eta(\widehat\pi_{\mathrm{fin}})
    =O_{P_0}(n_N^{-r}).
\]
The triangle inequality gives
\[
    \begin{aligned}
    &\left|
    \left(\sum_j\widehat w_j\right)^{-1}
    \sum_j\widehat w_j
    \{\eta(\widehat\pi_j)-\eta(\widehat\pi_{\mathrm{fin}})\}
    \right|\\
    &\quad\leq
    \left(\sum_j\widehat w_j\right)^{-1}
    \sum_j\widehat w_j
    \{\eta^*-\eta(\widehat\pi_j)\}
    +\{\eta^*-\eta(\widehat\pi_{\mathrm{fin}})\}.
    \end{aligned}
\]
The scale bounds imply
$\widehat w_j/\sum_i\widehat w_i\leq C/n_N$ uniformly.
Consequently, the last display is at most
\[
    O_{P_0}(1)
    \left\{
    \frac1{n_N}\sum_{k=1}^{n_N}k^{-r}
    +n_N^{-r}\right\}
\!.
\]
If $1/2<r<1$, the braces are $O(n_N^{-r})$; if $r=1$, they are
$O(n_N^{-1}\log n_N)$; and if $r>1$, they are $O(n_N^{-1})$.
Every case is $o(n_N^{-1/2})$.
Since the scale bounds also give
$\sigma_{R_{1N}}\asymp n_N^{-1/2}$, the difference between
$R_{\mathrm{TCI},1N}$ and $R_{\mathrm{CLB},1N}$ is
$o_{P_0}(\sigma_{R_{1N}})$.  Theorem
\ref{thm:CLT_for_R1_CLB}, scale consistency, and Slutsky's theorem
prove both conclusions of Theorem~\ref{thm:CLT_for_R1_TCI}.

\subsection{Proof of Theorem \ref{thm:o_P_for_R2_TCI} and
Corollary \ref{cor:CLT_for_NSAVE}}

Put
$\delta_{\mathrm{fin}}
  =\|\widehat Q_{\mathrm{fin}}^*-Q_{P_0}^*\|_\infty$.
Assumption~\ref{ass:est_pol} gives
$\delta_{\mathrm{fin}}=O_{P_0}(n_N^{-\kappa_*})$.
Lemma~\ref{lem:greedy_regret_margin} therefore yields
\[
    \begin{aligned}
    0\leq\eta^*-\eta(\widehat\pi_{\mathrm{fin}})
    &\leq C_{\mathrm{reg}}\delta_{\mathrm{fin}}^{1+\alpha}\\
    &=O_{P_0}\!
    \left(n_N^{-\kappa_*(1+\alpha)}\right)
    =o_{P_0}(n_N^{-1/2}),
    \end{aligned}
\]
because $\kappa_*(1+\alpha)>1/2$.  Since
$R_{\mathrm{TCI},2N}
  =\eta(\widehat\pi_{\mathrm{fin}})-\eta^*$, this proves
Theorem~\ref{thm:o_P_for_R2_TCI}.
Finally,
\[
    \widehat\eta_{\mathrm{NSAVE}}-\eta^*
    =R_{\mathrm{TCI},1N}+R_{\mathrm{TCI},2N}.
\]
The first term has the studentized Gaussian limit from Theorem
\ref{thm:CLT_for_R1_TCI}.  By the scale bounds and scale consistency,
the second is
$o_{P_0}(n_N^{-1/2})
  =o_{P_0}(\widehat\sigma_{R_{1N}})$.
Slutsky's theorem proves
Corollary~\ref{cor:CLT_for_NSAVE}.

\subsection{Proof of Theorem \ref{thm:DR_under_policy_alignment}}

Using the notation from the CLB proof, the following decomposition is
exact:
\[
    \widehat\eta_{\mathrm{NSAVE}}
    -\overline\eta_w(\widehat\pi)
    =
    \frac{\sum_{j=\ell_N+1}^N\widehat w_jM_j}
    {\sum_{j=\ell_N+1}^N\widehat w_j}
    +
    \frac{\sum_{j=\ell_N+1}^N\widehat w_jB_j}
    {\sum_{j=\ell_N+1}^N\widehat w_j}.
\]
Lemma~\ref{lem:exact_dr_identity} and the norm comparison used in the
CLB proof give
\[
    |B_j|
    \leq C
    \|\widehat Q_j-Q^{\widehat\pi_j}\|_{P_0,2}
    \|\widehat\omega_j-\omega^{\widehat\pi_j}\|_{P_0,2}.
\]
On the event in Assumption~\ref{ass:non_zero_var},
\[
    \sum_j\widehat w_j\geq
    \frac{n_N}{\overline\sigma},
    \qquad
    \max_j\widehat w_j\leq\frac1{\underline\sigma}.
\]
Consequently,
\[
    \begin{aligned}
    \left|
    \frac{\sum_j\widehat w_jB_j}
    {\sum_j\widehat w_j}\right|
    &\leq
    \frac{C\overline\sigma}{\underline\sigma n_N}
    \sum_{j=\ell_N+1}^N
    \|\widehat Q_j-Q^{\widehat\pi_j}\|_{P_0,2}
    \|\widehat\omega_j-\omega^{\widehat\pi_j}\|_{P_0,2}\\
    &=o_{P_0}(1)
    \end{aligned}
\]
by the assumed average product condition.

The numerator of the martingale term is a martingale because
$\widehat w_j$ is $\mathcal F_{j-1}$-measurable.  Its predictable
quadratic variation satisfies
\[
    \begin{aligned}
    \sum_{j=\ell_N+1}^N
    \widehat w_j^2
    \E(M_j^2\mid\mathcal F_{j-1})
    &=
    \sum_{j=\ell_N+1}^N
    \widehat w_j^2
    \widetilde\sigma_{\tau(j-1)}^2\\
    &\leq
    n_N\frac{\overline\sigma^2}{\underline\sigma^2}
    =O(n_N).
    \end{aligned}
\]
Applying the same Lenglart inequality to
$n_N^{-1/2}\sum_j\widehat w_jM_j$ shows that this scaled martingale is
$O_{P_0}(1)$.  Hence
$\sum_j\widehat w_jM_j=O_{P_0}(\sqrt{n_N})$.
Dividing by
$\sum_j\widehat w_j\geq n_N/\overline\sigma$ shows that the
martingale term is $O_{P_0}(n_N^{-1/2})=o_{P_0}(1)$.
We have proved
\[
    \widehat\eta_{\mathrm{NSAVE}}
    -\overline\eta_w(\widehat\pi)=o_{P_0}(1).
\]
The assumed value convergence now gives
\[
    \widehat\eta_{\mathrm{NSAVE}}-\eta^*
    =
    \{\widehat\eta_{\mathrm{NSAVE}}
      -\overline\eta_w(\widehat\pi)\}
    +\{\overline\eta_w(\widehat\pi)-\eta^*\}
    \overset{P_0}{\to}0.
\]

For the sufficient double-robustness condition, suppose, for example,
\[
    \max_j\|\widehat Q_j-Q^{\widehat\pi_j}\|_{P_0,2}
    =o_{P_0}(1),
    \qquad
    \max_j\|\widehat\omega_j-\omega^{\widehat\pi_j}\|_{P_0,2}
    =O_{P_0}(1).
\]
Then the average product is bounded by the product of these two
maxima and is therefore $o_{P_0}(1)$.  The same argument applies with
the roles of $Q$ and $\omega$ interchanged.

For completeness, Assumptions~\ref{ass:bellman_complete}--\ref{ass:margin_con}
also provide one sufficient route to the policy-value condition in
the theorem.  The scale bounds and the sequential regret calculation
give
\[
    \begin{aligned}
    0\leq
    \eta^*-\overline\eta_w(\widehat\pi)
    &=
    \frac{\sum_j\widehat w_j
    \{\eta^*-\eta(\widehat\pi_j)\}}
    {\sum_j\widehat w_j}\\
    &\leq
    O_{P_0}(1)\frac1{n_N}
    \sum_{k=1}^{n_N}k^{-\kappa_*(1+\alpha)}
    =o_{P_0}(1).
    \end{aligned}
\]
Thus the assumed convergence of
$\overline\eta_w(\widehat\pi)$ is automatic under those stronger
policy-learning conditions, although it can hold under weaker ones.

\subsection{Proof of Theorem \ref{thm:traj_eif} and
Corollaries \ref{cor:eff_est}--\ref{cor:eff_opt}}
\label{sec:traj_eff_proofs}

\emph{Tangent space of the i.i.d.-episode model.}
The likelihood of one episode factorizes as
\[
    p(O_{0:T})
    =\mu_0(S_0)\prod_{t=0}^{T}b(A_t\mid S_t)
     \prod_{t=0}^{T}g(R_t,S_{t+1}\mid S_t,A_t),
\]
with the \emph{same} $b$ and $g$ at every $t$
(Assumption~\ref{ass:Mark}).  Hence a regular parametric submodel
$\{P_\epsilon\}$ through $P_0$ has score
\[
    \dot\ell(O_{0:T})
    =a(S_0)+\sum_{t=0}^{T}\beta(A_t,S_t)
    +\sum_{t=0}^{T}\psi(R_t,S_{t+1},A_t,S_t),
\]
where $\E\{a(S_0)\}=0$, $\E\{\beta(A,S)\mid S\}=0$, and
$\E\{\psi(R,S',A,S)\mid S,A\}=0$.  To keep the martingale
conditioning explicit, write
\[
    \mathcal H_t^-:=
    \sigma(S_0,A_0,R_0,\ldots,S_{t-1},A_{t-1},R_{t-1},S_t),
    \qquad
    \mathcal H_t:=\mathcal H_t^-\vee\sigma(A_t),
\]
with the convention $\mathcal H_0^-=\sigma(S_0)$.  Then
$(R_t,S_{t+1})$ is generated after $\mathcal H_t$, and
$(R_t,S_{t+1})$ is included in $\mathcal H_{t+1}^-$.  The tangent space is the
orthogonal sum
$\mathcal T=\mathcal T_\mu\oplus\mathcal T_b\oplus\mathcal T_g$ of the
$L_2(P_0)$-closures of the three displayed blocks; orthogonality follows
from the tower property, since
$\E\{\beta(A_t,S_t)\mid\mathcal H_t^{-}\}=0$ and
$\E\{\psi(\cdot_t)\mid\mathcal H_t\}=0$ make cross-products of distinct
blocks vanish.

\emph{Pathwise derivative.}
Because $\eta(\pi)=\E_{\mu_0}\{V^\pi(S_0)\}$ depends on $(\mu_0,g)$ but
not on $b$,
\[
    \frac{\mathrm d}{\mathrm d\epsilon}\eta(\pi)\Big|_{0}
    =\E\{(V^\pi(S_0)-\eta(\pi))\,a(S_0)\}
    +\E_{\mu_0}\Big\{\tfrac{\mathrm d}{\mathrm d\epsilon}
    V^\pi_\epsilon(S_0)\Big|_0\Big\}.
\]
For the second term, the value-perturbation identity of
Lemma~\ref{lem:finite_mdp_uniform_expansion} gives
$\tfrac{\mathrm d}{\mathrm d\epsilon}V^\pi_\epsilon
=G_{0,\pi}\{\dot r_\pi+\gamma\dot K_\pi V^\pi\}$, and the Neumann series
$\mu_0^\top G_{0,\pi}=\frac1{1-\gamma}\E_{S\sim d_{P_0}^\pi}$ together
with the reward/transition score identities used in that lemma yield
\[
    \E_{\mu_0}\Big\{\tfrac{\mathrm d}{\mathrm d\epsilon}
    V^\pi_\epsilon(S_0)\Big\}
    =\frac1{1-\gamma}\E_{S\sim d_{P_0}^\pi}
    \Big\{\sum_a\pi(a\mid S)\,\E[\varepsilon^\pi\psi\mid S,a]\Big\}
    =\frac1{1-\gamma}
    \E_{P_{0,b}^{\mathrm{tr}}}\{\omega^\pi(A,S)\varepsilon^\pi\psi\},
\]
the last equality changing measure through
$\omega^\pi(a,s)=d_{P_0}^\pi(s)\pi(a\mid s)/\{f_b(s)b(a\mid s)\}$.

\emph{The canonical gradient.}  Set
\[
    \phi^{\mathrm{eff}}_T
    =\underbrace{V^\pi(S_0)-\eta(\pi)}_{\in\,\mathcal T_\mu}
    +\underbrace{\frac{1}{(1-\gamma)(T+1)}
    \sum_{t=0}^T\omega^\pi(A_t,S_t)\varepsilon^\pi_t}
    _{=\,\sum_t\psi^\star(\cdot_t)},
    \qquad
    \psi^\star:=\frac{\omega^\pi\varepsilon^\pi}{(1-\gamma)(T+1)} .
\]
Since $\E\{\psi^\star\mid S,A\}=0$, the second block lies in
$\mathcal T_g$, so $\phi^{\mathrm{eff}}_T\in\mathcal T$ with no
$\mathcal T_b$ component.  We check
$\langle\phi^{\mathrm{eff}}_T,\dot\ell\rangle
=\frac{\mathrm d}{\mathrm d\epsilon}\eta(\pi)$ for every score.  By
block orthogonality only matching blocks survive:
$\langle V^\pi(S_0)-\eta,\,a\rangle=\E\{(V^\pi(S_0)-\eta)a\}$, and
\[
    \Big\langle\sum_t\psi^\star(\cdot_t),
    \sum_{t'}\psi(\cdot_{t'})\Big\rangle
    =\sum_{t}\E\{\psi^\star(\cdot_t)\psi(\cdot_t)\}
    =\frac{1}{1-\gamma}
    \E_{P_{0,b}^{\mathrm{tr}}}\{\omega^\pi\varepsilon^\pi\psi\},
\]
because the off-diagonal terms vanish: for $t'<t$,
$\psi(\cdot_{t'})\in\mathcal H_t$ while
$\E\{\varepsilon^\pi_t\mid\mathcal H_t\}=0$; for $t'>t$,
$\omega^\pi(A_t,S_t)\varepsilon^\pi_t\in\mathcal H_{t'}$ while
$\E\{\psi(\cdot_{t'})\mid\mathcal H_{t'}\}=0$; and on the diagonal
stationarity makes each of the $T+1$ summands equal to
$\E\{\omega^\pi\varepsilon^\pi\psi\}$.  The remaining cross terms
$\langle\sum_t\psi^\star,\sum\beta\rangle$ and
$\langle\sum_t\psi^\star,a\rangle$ vanish by the same martingale
identities, consistent with $\eta(\pi)$ not depending on $b$.  Summing
the two surviving inner products reproduces the pathwise derivative.
Therefore $\phi^{\mathrm{eff}}_T$ is a gradient that lies in
$\mathcal T$, hence it is the canonical gradient (efficient influence
function), and $\eta(\pi)$ is pathwise differentiable.

\emph{Variance.}  The terms $\omega^\pi(A_t,S_t)\varepsilon^\pi_t$ are
martingale increments in the sense that they are
$\mathcal H_{t+1}^-$-measurable and have conditional mean zero given
$\mathcal H_t$:
$\E\{\omega^\pi(A_t,S_t)\varepsilon^\pi_t\mid\mathcal H_t\}=0$.  Hence
they are pairwise uncorrelated by the tower property, and
$V^\pi(S_0)$ is uncorrelated with each of them by the same argument.
Writing
$v:=\var_{P_{0,b}^{\mathrm{tr}}}\{\omega^\pi\varepsilon^\pi\}$ (equal
across $t$ by stationarity),
\[
    \var\{\phi^{\mathrm{eff}}_T\}
    =\var\{V^\pi(S_0)\}
    +\frac{(T+1)\,v}{(1-\gamma)^2(T+1)^2}
    =\var\{V^\pi(S_0)\}+\frac{v}{(1-\gamma)^2(T+1)}
    =\sigma^2_{\mathrm{eff},T}(\pi),
\]
which proves Theorem~\ref{thm:traj_eif}.  The same computation applied
to the discounted weights $\gamma^t/(1-\gamma^{T+1})$ gives transition
variance $v(1+\gamma^{T+1})/\{(1-\gamma^2)(1-\gamma^{T+1})\}$, yielding
the comparison stated after Corollary~\ref{cor:eff_opt}.

\emph{Efficient estimator (Corollary~\ref{cor:eff_est}).}
The map $(Q,\omega,V)\mapsto\Phi^{\mathrm{eff}}_T$ is the uniformly
weighted analogue of the score analyzed in
Theorem~\ref{thm:smoothing_RAL}.  We spell out the corresponding
product remainder.  Write
$\delta Q=\widehat Q-Q^\pi$, $\delta V=\widehat V-V^\pi$, and
$\delta\omega=\widehat\omega-\omega^\pi$.  By stationarity, every
transition in the uniform score has the same marginal
$P_{0,b}^{\mathrm{tr}}$, and the uniform weights sum to
$1/(1-\gamma)$.  Therefore
\[
    \begin{aligned}
    &P_0\{\Phi^{\mathrm{eff}}_T
      (O_{0:T};\widehat Q,\widehat\omega,\widehat V,\pi)\}
      -\eta(\pi)\\
    &\qquad
      =
      \frac1{1-\gamma}
      \E_{P_{0,b}^{\mathrm{tr}}}
      \!\left[
      \widehat\omega(A,S)
      \{\varepsilon^\pi+\gamma\delta V(S')-\delta Q(A,S)\}
      \right]
      +\E_{\mu_0}\{\delta V(S_0)\}.
    \end{aligned}
\]
The term involving $\widehat\omega\,\varepsilon^\pi$ is zero because
$\E(\varepsilon^\pi\mid S,A)=0$ and $\widehat\omega$ is a function of
$(S,A)$.  The true occupancy ratio also satisfies
\[
    \frac1{1-\gamma}
    \E_{P_{0,b}^{\mathrm{tr}}}
    \{\omega^\pi(A,S)[\gamma\delta V(S')-\delta Q(A,S)]\}
    +\E_{\mu_0}\{\delta V(S_0)\}=0,
\]
the Bellman identity used in
Lemma~\ref{lem:exact_dr_identity}.  Subtracting this zero quantity
leaves exactly
\[
    \frac1{1-\gamma}\E_{P_{0,b}^{\mathrm{tr}}}
    \{(\widehat\omega-\omega^\pi)
    [\gamma(\widehat V-V^\pi)(S')-(\widehat Q-Q^\pi)(A,S)]\}.
\]
Under cross-fitting the empirical-process term is
$o_{P_0}(N^{-1/2})$ by conditional Chebyshev (Step~1 of the proof of
Theorem~\ref{thm:smoothing_RAL}), and the product remainder is
$o_{P_0}(N^{-1/2})$ by assumption; hence
$\sqrt N\{\widehat\eta_{\mathrm{eff}}(\pi)-\eta(\pi)\}
=\frac1{\sqrt N}\sum_{i=1}^N\phi^{\mathrm{eff}}_T(O_i;\pi)+o_{P_0}(1)
\rightsquigarrow\mathcal N(0,\sigma^2_{\mathrm{eff},T}(\pi))$.  The
limiting variance equals the bound of Theorem~\ref{thm:traj_eif}, so the
estimator is efficient.

\emph{Uniqueness (Corollary~\ref{cor:eff_opt}).}
Under Assumption~\ref{ass:deterministic_policy} the envelope identity
\eqref{eq:optimal_value_envelope} holds with the linear coefficient
$L_\pi(\dot\ell)=\langle\phi^{\mathrm{eff}}_T(\cdot;\pi),\dot\ell\rangle$,
because the value expansion of
Lemma~\ref{lem:finite_mdp_uniform_expansion} is intrinsic to the MDP and
its episode-level Riesz representation is exactly the inner product
computed above; the required continuity of $\pi\mapsto
L_\pi(\dot\ell)$ follows from
Lemma~\ref{lem:episode_score_continuity}, whose proof applies verbatim
to the uniform weights.  With a singleton optimal set the one-sided
derivatives coincide, so $\Psi^*$ is pathwise differentiable in the
i.i.d.-episode model with canonical gradient
$\phi^{\mathrm{eff}}_T(\cdot;\pi^*_0)$ and bound
$\sigma^2_{\mathrm{eff},T}(\pi^*_0)$.  Finally, with
$\widehat\pi_{\mathrm{fin}}$ trained by sample splitting,
\[
    \sqrt N\{\widehat\eta_{\mathrm{eff}}(\widehat\pi_{\mathrm{fin}})
    -\Psi^*(P_0)\}
    =\sqrt N\{\widehat\eta_{\mathrm{eff}}(\widehat\pi_{\mathrm{fin}})
    -\eta(\widehat\pi_{\mathrm{fin}})\}
    +\sqrt N\{\eta(\widehat\pi_{\mathrm{fin}})-\Psi^*(P_0)\}.
\]
The first term obeys Corollary~\ref{cor:eff_est} conditionally on the
training split, with
$\|\phi^{\mathrm{eff}}_T(\cdot;\widehat\pi_{\mathrm{fin}})
-\phi^{\mathrm{eff}}_T(\cdot;\pi^*_0)\|_{P_0,2}=o_{P_0}(1)$ by
Lemma~\ref{lem:episode_score_continuity} and
$d_\Pi(\widehat\pi_{\mathrm{fin}},\pi^*_0)\overset{P_0}{\to}0$; the second is
$o_{P_0}(1)$ after scaling because
$0\le\Psi^*(P_0)-\eta(\widehat\pi_{\mathrm{fin}})
\le C_{\mathrm{reg}}\|\widehat Q^*_{\mathrm{fin}}-Q^*_{P_0}\|_\infty
^{1+\alpha}=o_{P_0}(N^{-1/2})$ by
Lemma~\ref{lem:greedy_regret_margin} and
Assumption~\ref{ass:est_pol} (as in Theorem~\ref{thm:o_P_for_R2_TCI}).
Slutsky's theorem gives Corollary~\ref{cor:eff_opt}.

\section{Technical Proofs in Section \ref{sec:alternative}}

\subsection{Proof of Theorem \ref{thm:smoothing_RAL}}\label{sec:smoothing_proofs}

\begin{proof}
    We divide the proof into five steps.

    \emph{Step 1: conditional one-step expansion.}
    Conditional on the training sample,
    $\widehat\pi_{\beta_N}$ and all policy-evaluation nuisance
    estimates are fixed on the evaluation fold.  Let
    $\widehat\psi_N$ denote the uncentered normalized trajectory score
    using these estimated nuisances and policy, and let
    $\psi_N^0$ denote the same score using
    $(Q^{\widehat\pi_{\beta_N}},
      \omega^{\widehat\pi_{\beta_N}},
      V^{\widehat\pi_{\beta_N}})$.
    Then
    $P_0\psi_N^0=\eta(\widehat\pi_{\beta_N})$ and
    \[
        \begin{aligned}
        \widehat\eta_{\beta_N}
        -\eta(\widehat\pi_{\beta_N})
        &=(\mathbb P_N-P_0)
        \{\psi_N^0-\eta(\widehat\pi_{\beta_N})\}\\
        &\quad+
        (\mathbb P_N-P_0)(\widehat\psi_N-\psi_N^0)
        +P_0(\widehat\psi_N-\psi_N^0).
        \end{aligned}
    \]
    The first term is
    $(\mathbb P_N-P_0)
      \phi_T(\cdot;\widehat\pi_{\beta_N})$.
    Lemma~\ref{lem:exact_dr_identity} and the finite-window
    stationarity calculation in its proof imply
    \[
        |P_0(\widehat\psi_N-\psi_N^0)|
        \leq C
        \|\widehat Q_{\mathrm{eval}}
        -Q^{\widehat\pi_{\beta_N}}\|_{P_0,2}
        \|\widehat\omega_{\mathrm{eval}}
        -\omega^{\widehat\pi_{\beta_N}}\|_{P_0,2}
        =o_{P_0}(N^{-1/2}).
    \]
    To verify the empirical-process norm explicitly, write
    $e_\omega=\widehat\omega_{\mathrm{eval}}
      -\omega^{\widehat\pi_{\beta_N}}$,
    $e_Q=\widehat Q_{\mathrm{eval}}
      -Q^{\widehat\pi_{\beta_N}}$, and
    $e_V(s)=\sum_a\widehat\pi_{\beta_N}(a\mid s)e_Q(a,s)$.
    At every transition, the difference between the estimated and
    true transition terms is
    \[
        e_\omega(A,S)
        \{R+\gamma V^{\widehat\pi_{\beta_N}}(S')
          -Q^{\widehat\pi_{\beta_N}}(A,S)\}
        +
        \widehat\omega_{\mathrm{eval}}(A,S)
        \{\gamma e_V(S')-e_Q(A,S)\}.
    \]
    The true residual and estimated ratio are uniformly bounded, and
    the same Jensen--overlap calculation used in the CLB proof gives
    $\|e_V(S')\|_{P_0,2}\leq C\|e_Q\|_{P_0,2}$.
    Summing the finitely many normalized transition terms and adding
    the initial term $e_V(S_0)$ therefore yields
    \[
        \|\widehat\psi_N-\psi_N^0\|_{P_0,2}
        \leq C_T\{\|e_\omega\|_{P_0,2}
          +\|e_Q\|_{P_0,2}\}
        =o_{P_0}(1).
    \]
    Conditional independence of the evaluation fold therefore yields
    \[
        \E\!\left[
        \left.
        \left\{\sqrt N(\mathbb P_N-P_0)
        (\widehat\psi_N-\psi_N^0)\right\}^2
        \right|\text{training data}\right]
        \leq
        \|\widehat\psi_N-\psi_N^0\|_{P_0,2}^2
        =o_{P_0}(1).
    \]
    Conditional Chebyshev's inequality proves
    $(\mathbb P_N-P_0)(\widehat\psi_N-\psi_N^0)
      =o_{P_0}(N^{-1/2})$.  We have thus shown
    \[
        \widehat\eta_{\beta_N}
        -\eta(\widehat\pi_{\beta_N})
        =
        (\mathbb P_N-P_0)
        \phi_T(\cdot;\widehat\pi_{\beta_N})
        +o_{P_0}(N^{-1/2}).
    \]

    \emph{Step 2: error from the estimated softmax policy.}
    We first verify the softmax Lipschitz bound.  For vectors $q,q'$
    and $h=q'-q$, define
    $p_t(a):=\exp\{\beta(q_a+th_a)\}/
      \sum_b\exp\{\beta(q_b+th_b)\}$.  Direct differentiation gives
    \[
        \frac{\mathrm d}{\mathrm dt}p_t(a)
        =
        \beta p_t(a)
        \left\{h_a-\sum_bp_t(b)h_b\right\}.
    \]
    Therefore
    \[
        \sum_a\left|
        \frac{\mathrm d}{\mathrm dt}p_t(a)\right|
        \leq
        \beta\sum_ap_t(a)
        \left(|h_a|+\sum_bp_t(b)|h_b|\right)
        \leq2\beta\|h\|_\infty.
    \]
    Integrating from $t=0$ to $t=1$ and dividing the $L_1$ distance
    by two gives
    \[
        d_\Pi(\widehat\pi_{\beta_N},\pi_{\beta_N})
        \leq
        \beta_N
        \|\widehat Q_{\mathrm{opt}}-Q_{P_0}^*\|_\infty
        =\beta_N r_N.
    \]
    Lemma~\ref{lem:diff_pol_to_Q_fun}, integrated over $\mu_0$,
    therefore implies
    \[
        |\eta(\widehat\pi_{\beta_N})-\eta(\pi_{\beta_N})|
        \leq C\beta_Nr_N=o_{P_0}(N^{-1/2}).
    \]

    \emph{Step 3: deterministic population softmax bias.}
    Let $\Delta_{\min}(s)$ denote the smallest positive
    optimal-action gap, with $\Delta_{\min}(s)=+\infty$ if every
    action is optimal.  Write
    $\Delta_a(s):=V_{P_0}^*(s)-Q_{P_0}^*(a,s)$.  Cancelling the common
    factor $\exp\{\beta_N V_{P_0}^*(s)\}$ from the softmax
    probabilities gives
    \[
        \pi_{\beta_N}(a\mid s)
        =
        \frac{\exp\{-\beta_N\Delta_a(s)\}}
        {\sum_b\exp\{-\beta_N\Delta_b(s)\}}.
    \]
    The denominator is at least one because at least one action has
    zero gap.  Thus
    \[
        \begin{aligned}
        \sum_a\pi_{\beta_N}(a\mid s)
        \Delta_a(s)
        &\leq
        \sum_{a:\Delta_a(s)>0}
        \Delta_a(s)e^{-\beta_N\Delta_a(s)}\\
        &\leq
        \frac{2|\mathcal A|}{e\beta_N}
        e^{-\beta_N\Delta_{\min}(s)/2}.
        \end{aligned}
    \]
    The last inequality uses
    $xe^{-\beta x/2}\leq2/(e\beta)$ and
    $e^{-\beta x/2}\leq
      e^{-\beta\Delta_{\min}(s)/2}$ for every positive gap $x$.
    Using an optimal greedy reference policy with
    $Q^{\pi^*}=Q^*$, the performance-difference lemma gives
    \[
        \begin{aligned}
        0\leq\eta^*-\eta(\pi_{\beta_N})
        &=
        \frac1{1-\gamma}
        \E_{S\sim d_{P_0}^{\pi_{\beta_N}}}
        \sum_a\pi_{\beta_N}(a\mid S)\Delta_a(S)\\
        &\leq
        \frac{C}{\beta_N}
        \E_{S\sim d_{P_0}^{\pi_{\beta_N}}}
        e^{-\beta_N\Delta_{\min}(S)/2}.
        \end{aligned}
    \]

    \emph{Step 4: use of the margin condition.}
    Put $\lambda=\beta_N/2$ and let
    $F_\pi(x):=
      \Pr_{S\sim d_{P_0}^\pi}\{0<\Delta_{\min}(S)\leq x\}$.
    Choose $\delta_0>0$ such that
    $\sup_{\pi\in\Pi}F_\pi(x)\leq Cx^\alpha$ for
    $0<x\leq\delta_0$.  Splitting at $\delta_0$ and integrating by
    parts,
    \[
        \begin{aligned}
        \E_{d^\pi}
        \{e^{-\lambda\Delta_{\min}(S)}
        \mathds{1}(\Delta_{\min}(S)<\infty)\}
        &\leq e^{-\lambda\delta_0}
        +\int_{(0,\delta_0]}e^{-\lambda x}\,\mathrm dF_\pi(x)\\
        &\leq
        2e^{-\lambda\delta_0}
        +\lambda\int_0^{\delta_0}e^{-\lambda x}F_\pi(x)\,\mathrm dx\\
        &\leq
        2e^{-\lambda\delta_0}
        +C\lambda\int_0^\infty e^{-\lambda x}x^\alpha\,\mathrm dx\\
        &=
        2e^{-\lambda\delta_0}
        +C\Gamma(1+\alpha)\lambda^{-\alpha}
        =O(\lambda^{-\alpha}),
        \end{aligned}
    \]
    uniformly in $\pi$.  Applying this bound to
    $\pi=\pi_{\beta_N}$ in the preceding display yields
    \[
        0\leq\eta^*-\eta(\pi_{\beta_N})
        \leq C\beta_N^{-(1+\alpha)}
        =o(N^{-1/2}).
    \]
    Together with Step 2,
    \[
        \eta(\widehat\pi_{\beta_N})-\eta^*
        =o_{P_0}(N^{-1/2}).
    \]

    \emph{Step 5: stabilization of the random score.}
    For each state, the population softmax converges to equal mass on
    the maximizing actions and zero mass on suboptimal actions.
    Since the state and action spaces are finite,
    $d_\Pi(\pi_{\beta_N},\pi_\infty)\to0$.  Step 2 and
    $\beta_Nr_N=o_{P_0}(N^{-1/2})$ imply
    \[
        d_\Pi(\widehat\pi_{\beta_N},\pi_\infty)
        \leq
        d_\Pi(\widehat\pi_{\beta_N},\pi_{\beta_N})
        +d_\Pi(\pi_{\beta_N},\pi_\infty)
        \overset{P_0}{\to}0.
    \]
    By Lemma~\ref{lem:episode_score_continuity},
    \[
        \|\phi_T(\cdot;\widehat\pi_{\beta_N})
          -\phi_T(\cdot;\pi_\infty)\|_{P_0,2}
        =o_{P_0}(1).
    \]
    Conditional on the training data,
    \[
        \begin{aligned}
        &\E\!\left[
        \left.
        \left\{\sqrt N(\mathbb P_N-P_0)
        [\phi_T(\cdot;\widehat\pi_{\beta_N})
        -\phi_T(\cdot;\pi_\infty)]\right\}^2
        \right|\text{training data}\right]\\
        &\qquad\leq
        \|\phi_T(\cdot;\widehat\pi_{\beta_N})
          -\phi_T(\cdot;\pi_\infty)\|_{P_0,2}^2
        =o_{P_0}(1).
        \end{aligned}
    \]
    Conditional Chebyshev's inequality gives
    \[
        (\mathbb P_N-P_0)
        \{\phi_T(\cdot;\widehat\pi_{\beta_N})
        -\phi_T(\cdot;\pi_\infty)\}
        =o_{P_0}(N^{-1/2}).
    \]

    Combining Steps 1--5,
    \[
        \begin{aligned}
        \sqrt N(\widehat\eta_{\beta_N}-\eta^*)
        &=
        \sqrt N(\mathbb P_N-P_0)
        \phi_T(\cdot;\pi_\infty)+o_{P_0}(1)\\
        &=
        \frac1{\sqrt N}\sum_{i=1}^N
        \phi_T(O_i;\pi_\infty)+o_{P_0}(1).
        \end{aligned}
    \]
    The summands in the last line are i.i.d., centered, and
    square-integrable.  The ordinary central limit theorem therefore
    gives the asserted pointwise Gaussian limit with variance
    $\var\{\phi_T(O;\pi_\infty)\}$.
\end{proof}

\subsection{Proof of Corollary \ref{cor:PSI}}

Let
\[
    \begin{aligned}
    \mathcal E_1
    &:=
    \{\bm\eta\in
    \mathcal C_\eta(\widehat{\bm\eta};\delta_1)\},\\
    \mathcal E_2
    &:=
    \left\{
    \max_{k\in\widehat{\mathcal A}^{+}}|Z_k|
    \leq
    q_N(\widehat{\mathcal A}^{+})\right\}.
    \end{aligned}
\]
By construction of the plausible set and the fact that the centered
confidence region contains $\widehat{\bm\eta}$,
\[
    \widehat{\mathcal A}_{\mathrm{opt}}
    =\arg\max_k\widehat\eta(\pi_k)
    \subseteq\widehat{\mathcal A}^{+}.
\]
On $\mathcal E_2$, every $k\in\widehat{\mathcal A}^{+}$ satisfies
\[
    \sqrt N\,
    \widehat{\boldsymbol\Sigma}_{kk}^{-1/2}
    |\widehat\eta(\pi_k)-\eta(\pi_k)|
    \leq
    q_N(\widehat{\mathcal A}^{+}).
\]
Equivalently,
\[
    \eta(\pi_k)\in
    \left[
    \widehat\eta(\pi_k)\pm
    q_N(\widehat{\mathcal A}^{+})
    \sqrt{\widehat{\boldsymbol\Sigma}_{kk}/N}
    \right].
\]
Since every empirically selected coordinate belongs to
$\widehat{\mathcal A}^{+}$, all factors in
$\mathcal C_{\mathrm{PSI}}$ cover their corresponding true values on
$\mathcal E_1\cap\mathcal E_2$.

The two advertised probability guarantees and the union bound give
\[
    \begin{aligned}
    \liminf_{N\to\infty}\Pr(\mathcal E_1\cap\mathcal E_2)
    &\geq
    1-\limsup_{N\to\infty}\Pr(\mathcal E_1^c)
      -\limsup_{N\to\infty}
       \Pr(\mathcal E_1\cap\mathcal E_2^c)\\
    &\geq
    1-\delta_1-(\delta_2-\delta_1)
    =1-\delta_2.
    \end{aligned}
\]
The coverage event contains $\mathcal E_1\cap\mathcal E_2$, proving
the corollary.

\section{Auxiliary Lemmata}

\begin{lemma}[Performance Difference Lemma, see \cite{kakade2002approximately}]\label{AuLem:PerDiffLem}
    Suppose that Assumptions \ref{ass:data_obs} and \ref{ass:Mark}
    hold, $\gamma<1$, and rewards are uniformly bounded.  Let
    \[
        d_{P,s_0}^{\pi_2}
        :=(1-\gamma)\sum_{t=0}^{\infty}\gamma^t
        \mathcal L_P^{\pi_2}(S_t\mid S_0=s_0).
    \]
    Then
    \[
        V_P^{\pi_2}(s_0)-V_P^{\pi_1}(s_0)
        =
        \frac{1}{1-\gamma}
        \E_{S\sim d_{P,s_0}^{\pi_2}}
        \left[
        \E_{A\sim\pi_2(\cdot\mid S)}
        \mathbb A_P^{\pi_1}(A,S)\right],
    \]
    where
    $\mathbb A_P^\pi(a,s):=Q_P^\pi(a,s)-V_P^\pi(s)$.
    Integrating over $S_0\sim\mu_0$ gives the corresponding identity
    for $\eta(\pi_2)-\eta(\pi_1)$ with occupancy
    $d_P^{\pi_2}$.
\end{lemma}

\begin{proof}
    Generate a trajectory under $\pi_2$ starting from $S_0=s_0$.
    By the definitions of $Q_P^{\pi_1}$ and the advantage,
    \[
        \begin{aligned}
        \E\{\mathbb A_P^{\pi_1}(A_t,S_t)\}
        &=
        \E\{Q_P^{\pi_1}(A_t,S_t)-V_P^{\pi_1}(S_t)\}\\
        &=
        \E\{R_t+\gamma V_P^{\pi_1}(S_{t+1})
        -V_P^{\pi_1}(S_t)\}.
        \end{aligned}
    \]
    Multiply by $\gamma^t$ and sum from $t=0$ to $m$:
    \[
        \begin{aligned}
        \sum_{t=0}^m\gamma^t
        \E\{\mathbb A_P^{\pi_1}(A_t,S_t)\}
        &=
        \sum_{t=0}^m\gamma^t\E(R_t)\\
        &\quad+
        \sum_{t=0}^m\gamma^{t+1}
        \E\{V_P^{\pi_1}(S_{t+1})\}\\
        &\quad-
        \sum_{t=0}^m\gamma^t
        \E\{V_P^{\pi_1}(S_t)\}\\
        &=
        \sum_{t=0}^m\gamma^t\E(R_t)
        -V_P^{\pi_1}(s_0)\\
        &\quad+
        \gamma^{m+1}
        \E\{V_P^{\pi_1}(S_{m+1})\}.
        \end{aligned}
    \]
    Bounded rewards imply bounded $V_P^{\pi_1}$, so the final term
    tends to zero as $m\to\infty$.  The reward sum tends to
    $V_P^{\pi_2}(s_0)$.  Hence
    \[
        V_P^{\pi_2}(s_0)-V_P^{\pi_1}(s_0)
        =
        \sum_{t=0}^{\infty}\gamma^t
        \E\{\mathbb A_P^{\pi_1}(A_t,S_t)\}.
    \]
    By the definition of $d_{P,s_0}^{\pi_2}$, the right-hand side is
    \[
        \frac1{1-\gamma}
        \E_{S\sim d_{P,s_0}^{\pi_2}}
        \E_{A\sim\pi_2(\cdot\mid S)}
        \{\mathbb A_P^{\pi_1}(A,S)\}.
    \]
    Integrating both sides with respect to $\mu_0$ proves the final
    statement.
\end{proof}







\end{document}